\documentclass[11pt]{article}
\usepackage{amsmath}
\usepackage{amssymb}
\usepackage{amsthm}
\usepackage{amsfonts}
\usepackage{cite}
\usepackage{calc}
\usepackage{epsfig,afterpage}
\advance\textwidth by +1.4in
\advance\textheight by +1.8in
\advance\oddsidemargin by -0.5in
\advance\evensidemargin by -1.0in
\advance\topmargin by -0.5in
\catcode`\@=11

\@addtoreset{equation}{section}
\catcode`\@=12

\def\be{\begin{equation}}
\def\ee{\end{equation}}
\def\bal{\begin{aligned}}
\def\eal{\end{aligned}}
\def\bq{\begin{eqnarray}}
\def\eq{\end{eqnarray}}
\def\beq{\begin{eqnarray*}}
\def\eeq{\end{eqnarray*}}

\def\ba{\begin{array}}
\def\ea{\end{array}}
\def\bth{\begin{theorem}}
\def\eth{\end{theorem}}
\def\blm{\begin{lemma}}
\def\elm{\end{lemma}}
\def\bdf{\begin{definition}}
\def\edf{\end{definition}}
\def\bpr{\begin{proposition}}
\def\epr{\end{proposition}}
\def\brm{\begin{remark}}
\def\erm{\end{remark}}
\def\bnot{\begin{notation}}
\def\enot{\end{notation}}
\def\bobs{\begin{observation}}
\def\eobs{\end{observation}}
\def\bcrl{\begin{corrolary}}
\def\ecrl{\end{corrolary}}

\newcommand{\ab}{\mbox{\boldmath $a$}}
\newcommand{\bb}{\mbox{\boldmath $b$}}

\newcommand{\SSS}{{\bf S}}
\newcommand{\QSL}{{\bf Q\bf S\bf L}}

\newcommand{\QS}{{\bf Q\bf S}}
\newcommand{\IL}{{\bf I\bf L}}
\newcommand{\D}{{\bf D}}

\newcommand{\PP}{{\mathbb P}}
\newcommand{\R}{\mathbb{R}}
\newcommand{\C}{\mathbb{C}}

\newcommand{\sign}{\mbox{\rm sign\,}}
\newcommand{\Res}{\mbox{\rm Res\,}}
\newcommand{\Discriminant}{\mbox{\rm Discriminant\,}}
\newcommand{\Jacob}{\mbox{\rm Jacob\,}}
\newcommand{\Hess}{\mbox{\rm Hess\,}}
\newcommand{\Supp}{\mbox{\rm Supp\,}}

\newtheorem{theorem}{Theorem}[section]
\newtheorem{lemma}{Lemma}[section]
\newtheorem{definition}{Definition}[section]
\newtheorem{proposition}[lemma]{Proposition}
\newtheorem{remark}{Remark}[section]
\newtheorem{notation}[remark]{Notation}
\newtheorem{observation}[remark]{Observation}
\newtheorem{corrolary}[remark]{Corrolary}
\newcommand{\BProof}{\noindent{\it Proof:}\ \,}
\newcommand{\EProof}{\ \hfill\rule[-1mm]{1.5mm}{2.0mm}}

\begin{document}

\setcounter{secnumdepth}{5}

\title{Planar quadratic differential systems
         with invariant straight lines of the total multiplicity 4
 }

\author{Dana SCHLOMIUK\thanks{Work supported by NSERC and
        by the Quebec Education Ministry}\\
        D\'epartement de Math\'ematiques et de Statistiques\\
        Universit\'e de Montr\'eal
 \and
        Nicolae VULPE\thanks{Partially supported by  NSERC}\\
        Institute of Mathematics and Computer Science\\
        Academy of Science of Moldova }
\date{}
\maketitle
\begin{abstract}

In this article we  consider the action of affine group and time
rescaling on planar quadratic differential systems. We construct a
system of representatives of the orbits of systems with four
invariant lines, including the line at infinity and including
multiplicities.  For each orbit we  exhibit  its configuration. We
characterize in terms of algebraic invariants and comitants and
also geometrically, using divisors of the complex projective
plane, the class of quadratic differential systems with four
invariant lines. These conditions are such that no matter how a
system may be presented, one can verify by using them whether the
system has exactly  four invariant lines  including
multiplicities, and if it is so, to check to which orbit (or
family of orbits) it belongs.
\end{abstract}

\section{Introduction } We consider here real planar  differential systems of the form
\be\label{il1}
 (S)\qquad \frac {dx}{dt}= p(x,y),\qquad
  \frac {dy}{dt}= q(x,y),\hphantom{--  }
\ee
where $p,\,q\in \R[x,y]$, i.e. $p,\ q$ are polynomials in $x,\ y$
over $\R$, and their associated vector fields
 \be
     \tilde D=p(x,y)\frac{\partial}{\partial x} + q(x,y)\frac{\partial}{\partial y}.
     \label{il2}
\ee
Each such system generates a complex differential vector field when the
variables range over $\C$. To the complex systems we can apply the work of
Darboux on integrability via invariant algebraic curves (cf.\cite{Darb}).

For the systems (\ref{il1}) we can use the following definition.
 \bdf \label{df1}
 An affine algebraic invariant  curve of a polynomial system {\rm(\ref{il1})}
(or an algebraic particular integral) is a curve $f(x,y)=0$ where
$f\in \C[x,y]$, $\deg(f)\ge1$, such that there exists $k(x,y)\in
\C[x,y]$ satisfying $\tilde Df=fk$ in $\C[x,y]$. We call $k$ the
cofactor of $f$ with respect to the system. \edf

We are interested in polynomial systems (\ref{il1}) possessing
algebraic invariant curves. The presence of a sufficient number of
such curves implies integrability of the system via the geometric
method of integration of Darboux (cf.\cite{Darb}). For a brief
introduction to the work of Darboux we refer to the survey article
\cite {Dana1}.  Some applications of the work of Darboux in
connection with the problem of the center are given in
\cite{Dana2}. There is a growing literatures on problems related
to the work of Darboux on invariant algebraic curves of
differential equations. In particular we mention here the recent
work of C. Christopher, J.V. Perreira and J. Llibre
\cite{Cris_Llib_Per} on the notion of multiplicity of an invariant
algebraic curve of a differential system.

 In this article we shall
consider the simplest kind of such a structure, i.e. quadratic systems
(\ref{il1}) possessing invariant lines. Some references on this topic are:
\cite{Sib3,Druzhkova,Art_Llib2,Lyubim1, Lyubim2,Popa_Sib1,Popa2,Sokulski,ZX}.

To a line $f(x,y)=ux+vy+w=0$, $(u,v)\ne(0,0)$ we associate its
projective completion $F(X,Y,Z)=uX+vY+wZ=0$ under the embedding
$\C^2\hookrightarrow \PP_2(\C)$, $(x,y)\mapsto [x:y:1]$. The line
$Z=0$ is   called the line at infinity of the system (\ref{il1}).
It follows from the work of Darboux that each system of
differential equations of the form (\ref{il1}) yields a
differential equation on the complex projective plane which is the
compactification of the complex system (\ref{il1}) on $\PP_2(\C)$.
The line $Z=0$ is an  invariant manifold of this complex
differential equation.
\bnot
  Let us denote by
\beq
  \QS &=& \left\{\ S\ \left|\ba{l}\ S\ \hbox {is a system {\rm (\ref{il1})} such that}\
       \gcd(p(x,y),q(x,y))=1\\ \ \ \hbox{and}\quad \max\big(\deg(p(x,y)),\deg(q(x,y))\big)=2
            \ea \right.\right\}; \\
  \QSL &=& \left\{\,S\in \QS \left|\ba{l} \ \hbox{ $S$ possesses at
                                       least one invariant affine line or}\\
                                     \ \, \hbox{the line at infinity with multiplicity
                                       at least two}
                       \ea \right.\right\}.
\eeq
 \enot
\noindent
 For the multiplicity of the line at infinity see
\cite{Dana_Vlp1}.

We shall call \textit{degenerate quadratic differential system} a
system (\ref{il1}) with $\deg \gcd(p(x,y),q(x,y))\ge1$ and
$\max\big(\deg(p(x,y)),\deg(q(x,y))\big)=2$.
\bpr\label{pr:m_il} \mbox{\rm\!\!\cite{Art_Llib2}} The
maximum number of invariant lines (including the line at infinity
and including multiplicities) which a quadratic system could have
is six.\epr

 \bnot
 To a quadratic system {\rm(\ref{il1})} we can associate  a point in
 $ \R^{12}$, the
ordered tuple of the coefficients of $p(x,y)$, $q(x,y)$  and this
correspondence is an injection
\be\label{bigection}
\bal
{\cal B}:\quad \QS &\hookrightarrow \R^{12}\\
S\ \ & \mapsto\ \ \ab = {\cal B}(S).
 \eal
 \ee
The topology of $\R^{12}$ yields an induced topology on {\bf QS}.
 \enot

We associate to each system in {\bf QSL} its {\it configuration }
of invariant lines, i.e. the set of its invariant lines together
with the singular points of the systems located on the union of
these lines. In analogous manner to how we view the phase
portraits of the systems on the Poincar\'e disc (see e.g.
\cite{Lib_DS}), we can also view the configurations of real lines
on the disc. To help imagining the full configurations, we
complete the picture by drawing dashed lines whenever these are
complex.

On the class of quadratic systems acts the group of affine
transformations and time rescaling. Since quadratic systems depend
on 12 parameters and since this group depends on 7 parameters, the
class of quadratic systems modulo this group action, actually
depends on five parameters.

It is clear that the configuration of invariant lines of a system is an affine
invariant.

\bdf\label{def:multipl} We say that an invariant straight line ${\cal
L}(x,y)=ux+vy+w=0$ for a quadratic vector field $\tilde D$ has multiplicity $m$
if there exists a sequence of quadratic vector fields $\tilde D_k$ converging
to $\tilde D$, such that each $\tilde D_k$ has $m$ distinct invariant straight
lines ${\cal L}^1_k=0,\ldots, {\cal L}^m_k=0$, converging to ${\cal L}=0$ as
$k\to\infty$ (with the topology of their coefficients), and this does not occur
for $m+1$.
 \edf

The notion of multiplicity thus defined is invariant under the group action,
i.e. if a quadratic system $(S)$ has an invariant line $l$ of multiplicity $m$,
then each system $(\tilde S)$ in the orbit of $(S)$ under the group action has
an invariant line $l~$ of the same multiplicity $m$.

In this article  we continue the work initiated in
\cite{Dana_Vlp1} and  consider the case when the system
(\ref{il1}) has exactly four invariant lines considered with their
multiplicities.

The problems which we solve in this article are the following:

I) Construct a system of representatives of the orbits of systems
with exactly four invariant lines, including the line at infinity
and including multiplicities.
 For each orbit  exhibit  its configuration.

II) Characterize in terms of algebraic invariants and comitants
and also geometrically, using divisors of the complex projective
plane, the class of quadratic differential systems with four
invariant lines. These conditions should be such that no matter
how a system may be presented to us, we should be able to verify
by using them whether the system has or does not have four
invariant lines and   to check to which orbit or perhaps family of
orbits it belongs.

Our main results are formulated in Theorem \ref{th_mil_4}. Theorem
\ref{th_mil_4} gives a
  a complete list of
representatives of the orbits of systems with exactly four
invariant lines including the line at infinity and including
multiplicities.  These representatives are classified in 12
two-parameters families, 28  one-parameter families and 6 concrete
systems. We characterize each one of these 40 families in terms of
algebraic invariants or comitants and also geometrically. As the
calculation of invariants and comitants can be implemented on a
computer, this verification can be done by a computer.

The invariants and comitants of differential equations used in the
classification Theorem \ref{th_mil_4}  are obtained following the
theory established by K.Sibirsky and his disciples (cf.
\cite{Sib1}, \cite{Sib2}, \cite{Vlp1}, \cite{Popa4}).

\section{ Divisors associated to invariant lines configurations }
\label{divisors} Consider real quadratic systems, i.e. systems  of
the form:
\be\label{2l1}
(S)\qquad \left\{ \ba{ll}
\displaystyle  \frac {dx}{dt}&=p_0+ p_1(x,y)+\,p_2(x,y)\equiv p(x,y), \\[2mm]
\displaystyle  \frac {dx}{dt}&=q_0+ q_1(x,y)+\,q_2(x,y)\equiv
q(x,y)
\ea\right.
\ee
with  $\max(\deg(p),\deg(q))=2$, $\gcd(p,q)=1$  and
$$
\bal
&p_0=a_{00},\quad p_1(x,y)=  a_{10}x+ a_{01}y,\quad  p_2(x,y)=
a_{20}x^2 +2
a_{11}xy + a_{02}y^2,\\
&q_0=b_{00},\quad q_1(x,y)=  b_{10}x+ b_{01}y,\quad\  q_2(x,y)=
b_{20}x^2 +2
b_{11}xy + b_{02}y^2.\\
\eal
$$
 Let $
a=(a_{00},a_{10},a_{01},a_{20},a_{11},a_{02},b_{00},b_{10},b_{01},b_{20},b_{11},b_{02})$
be the 12-tuple of the coefficients of system \eqref{2l1} and
denote $\R[
a,x,y]=\R[a_{00},a_{10},a_{01},a_{20},a_{11},a_{02},b_{00},b_{10},b_{01},b_{20},b_{11},b_{02},x,y]$.
\bnot Let us denote by $\ab=(\ab_{00},\ab_{10}\ldots,\bb_{02})$ a
point in $\R^{12}$. Each  particular system \eqref{2l1} yields an
ordered 12-tuple $\ab$ of its coefficients.
\enot
\bnot Let
$$
\bal
  P(X,Y,Z)=& p_0(\ab)Z^2+ p_1(\ab,X,Y)Z+\,p_2(\ab,X,Y)=0,\\
  Q(X,Y,Z)=& q_0(\ab)Z^2+ q_1(\ab,X,Y)Z+\,q_2(\ab,X,Y)=0.\\
\eal
$$
We denote $\quad \sigma(P,Q)= \{w\in\PP_2(\C)\ |\ P(w)= Q(w)=0\}$.
\enot

\bdf\label{df3_1a} A formal expression of the form $\D = \sum
n(w)w$\ where $w \in\PP_2(\C)$, $n(w)$ is an integer and only a
finite number of the numbers $n(w)$ are not zero, will be called a
zero-cycle of $\PP_2(\C)$ and if $w$ only belongs to the line
$Z=0$ will be called a divisor of this line. We call  degree of
the expression  $\D$ the integer $\deg(\D) = \sum n(w)$. We call
support of   $\D$ the set $\Supp(\D)$ of points $w$ such that
$n(w)\ne0$. \edf

\bdf Let $C(X,Y,Z)=YP(X,Y,Z) - XQ(X,Y,Z)$.
\beq
    \D_{{}_S}(P,Q) &=& \sum_{w\in\sigma(P,Q)}I_w(P,Q)w;\\
    \D_{{}_S}(C,Z) &=& \sum_{w\in\{ Z = 0\}}I_w(C,Z)wa \quad \mbox{if}\quad Z\nmid C(X,Y,Z);\\
    \D_{{}_S}(P,Q;Z) &=& \sum_{w\in\{ Z = 0\}}I_w(P,Q)w;\\
    \widehat\D_{{}_S}(P,Q,Z) &=& \sum_{w\in\{ Z = 0\}}\Big(I_w(C,Z),\,I_w(P,Q)\Big)w,\\
\eeq
where $I_w(F,G)$ is the intersection number (see, \cite{Fult}) of
the  curves defined by homogeneous polynomials  $F,\ G\in
\C[X,Y,Z]$\ and $\deg(F),\deg(G)\ge1$.

\edf
\bnot
\be\label{invar1}
\bal
  n_{{}_\R}^{\!{}^\infty} =&\# \{ w\in Supp\, \D_{{}_S}(C,Z)\,\big|\, w\in \PP_2(\R)\};\\
  d_{\sigma}^{{}^\infty} =&\deg \D_{{}_S}(P,Q;Z).
\eal
\ee
\enot

A complex projective line $uX+vY+wZ=0$ is invariant for the system
$(S)$ if either it coincides with $Z=0$ or is the projective
completion of an invariant affine line $ux+vy+w=0$.

\bnot
 Let  $S\in \QSL$. Let us  denote
$$
\bal
    \IL(S)=&\left\{\ \ l\ \ \left|\ba{ll} & l\  \hbox {is a line in $\PP_2(\C)$ such }\\
                                           & \hbox{that}\  l\ \hbox{is invariant for}\ (S)\\
                      \ea\ \right\}\right.;\\
     M(l)=& \ \hbox{the multiplicity of the invariant line $l$ of
     $(S)$}.
\eal
$$
\enot
\brm
We note that the line  $l_\infty: Z=0$ is included  in $\IL(S)$
for any $S\in \QSL$.
\erm
 Let $l_i\,:\ f_i(x,y)=0$, $i=1,\ldots,k$,  be all the distinct invariant
affine lines (real or complex) of a system $S\in \QSL$. Let
$l'_i\,:\ {\cal F}_i(X,Y,Z)=0 $ be the complex projective
completion of $l_i$.
\bnot We denote
\beq
    && {\cal G}\,:\quad \ \prod_{i}{\cal F}_i(X,Y,Z)\,Z=0;\quad
        Sing\, {\cal G}=\left\{w\in {\cal G}|\ w\ \mbox{is a singular point of}\ {\cal G}\right\};  \\
    &&\nu(w)= \ \hbox{the multiplicity of the point $w$, as a point of}\
               {\cal G}.
\eeq
\enot
\bdf
\beq
    \D_{{}_\IL}(S)&=&\sum_{l\in\IL(S)} M(l)l,\quad (S)\in\QSL;\\
    Supp\, \D_{{}_\IL}(S)& =& \{\, l\ |\ l\in \IL(S)\}.
\eeq
\edf
\bnot
\be\label{invar}
\bal
  M_{{}_\IL}=&\deg\D_{{}_\IL}(S);\\
  N_{{}_\C} =&\# Supp\, \D_{{}_\IL};\\
  N_{{}_\R} =&\# \{   l\in Supp\, \D_{{}_\IL}\,\big|\,l\in  \PP_2(\R) \};\\
 n^{{}^\R}_{{}_{{\cal G},\,\sigma}}=&\# \{\omega\in Supp\, \D_{{}_S}(P,Q)\, |\,
                    \omega\in  {\cal G}\raisebox{-0.3em}[0pt][0pt]{$\big|_{\R^2}$}\};\\
 d^{\,{}^\R}_{{}_{{\cal G},\,\sigma}}=&\sum_{\omega\in{\cal G}\raisebox{-0.2em}[0pt][0pt]{$|_{\R^2}$}}I_\omega(P,Q);\\
  m_{{}_{\cal G}}=& \max\{\nu(\omega)\, |\, \omega\in Sing\, {\cal G}\};\\
  m_{{}_{\cal G}}^{\!{}^\infty}=& \max\{\nu(\omega)\, |\, \omega\in Sing\, {\cal
  G}\cap\{Z=0\}\}.
\eal
\ee
\enot
\section{The main $T$-comitants associated to configurations of
invariant lines}\label{T-comit}
 It is known that on the set {\QS} of all
quadratic differential systems \eqref{2l1} acts the group
$Aff(2,\R)$ of affine transformation on the plane \mbox{(cf. \cite
{Dana_Vlp1})}.
  For every  subgroup
$G\subseteq Aff(2,\R)$ we have an  induced action of $G$ on {\QS}.
We can identify the set {\QS} of systems \eqref{2l1} with a subset
of $\R^{12}$ via the map  \QS $\longrightarrow \R^{12}$  which
associates to
 each system \eqref{2l1} the 12-tuple $\ab=(\ab_{00},\ldots,\bb_{02})$ of
 its coefficients.

For the definitions of  an affine $GL$-comitant and invariant as
well as for the definition of a $T$-comitant and $CT$-comitant we
refer reader to the paper  \cite {Dana_Vlp1}. Here we  shall only
construct the necessary $T$-comitants associated to configurations
of invariant lines for the class of quadratic systems with exactly
four invariant lines including the line at infinity and including
multiplicities.

 Let us consider the polynomials
\beq
   C_i(a,x,y)&=&yp_i(a,x,y)-xq_i(a,x,y)\in \R[a,x,y],\ i=0,1,2, \\
  D_i(a,x,y)&=&\frac{\partial}{\partial x}p_i(a,x,y)+
        \frac{\partial}{\partial y}q_i(a,x,y)\in \R[a,x,y],\ i=1,2.
\eeq
As it was shown in \cite{Sib1} the polynomials
\be\label{C_i:D_i}
\big\{\ C_0(a,x,y),\quad C_1(a,x,y),\quad C_2(a,x,y),\quad D_1(a),
\quad D_2(a,x,y)\ \big\}
\ee
of degree one in the coefficients of systems \eqref{2l1} are
$GL$-comitants   of these systems.
\bnot Let $f,$
$g\in$ $\R[a,x,y]$ and
\be \label{trsv}
  (f,g)^{(k)}=
   \sum_{h=0}^k (-1)^h {k\choose h}
   \frac{\partial^k f}{ \partial x^{k-h}\partial y^h}\
   \frac{\partial^k g}{ \partial x^h\partial y^{k-h}}.
\ee
$(f,g)^{(k)}\in \R[a,x,y] $ is called the transvectant of index $k$ of $(f,g)$
(cf. {\rm\cite{Gr_Yng}, \cite{Olver}})
\enot

\bth\label{th:Vlp} \mbox{\rm \cite{Vlp1}}   Any $GL$-comitant  of systems \eqref{2l1} can
be constructed from the elements of the set \eqref{C_i:D_i} by
using the operations: $+,\, -,\,\times,$   and by applying the
differential operation $(f,g)^{(k)}$.
\eth
\bnot\label{not:1}  Consider the
polynomial  $\Phi_{\alpha,\beta}=\alpha P+\beta Q\in
\R[a,X,Y,Z,\alpha,\beta]$ where $P=Z^2p(X/Z,Y/Z),$
$Q=Z^2q(X/Z,Y/Z)$, $p,$ $q\in \R[a,x,y]$ and $\max
(\deg_{(x,y)}p,\deg_{(x,y)}q)=2$. Then
$$
\bal
&\Phi_{\alpha,\beta}=\ c_{11}(\alpha,\beta)X^2  +2
c_{12}(\alpha,\beta)XY
       + c_{22}(\alpha,\beta)Y^2+
   2c_{13}(\alpha,\beta)XZ  +2c_{23}(\alpha,\beta)YZ\\
&\qquad\qquad  +c_{33}(\alpha,\beta)Z^2,\qquad\quad
\Delta(a,\alpha,\beta) =\ \det\left|\left|c_{ij}(\alpha,\beta)
    \right|\right|_{i,j\in\{1,2,3\}},\\
&  D(a,\alpha,\beta) = 4\Delta(a,-\beta,\alpha),\qquad
H(a,\alpha,\beta)  = 4\big[\det\left|\left|c_{ij}(-\beta,\alpha)
    \right|\right|_{i,j\in\{1,2\}}\big].\\
\eal
$$
\enot
\blm\label{lem:H}\mbox{\rm\!\!\cite{Dana_Vlp1}} Consider two parallel invariant affine lines
\mbox{${\cal L}_i(x,y)\equiv ux+vy+w_i=0$}, ${\cal L}_i(x,y)\in
\C[x,y],$ $ (i=1,2)$ of a quadratic system $S$ of coefficients
$\ab$. Then \mbox{$H(\ab,\!-v,u)\!=\!0$}, i.e. the T-comitant
$H(a,x,y)$ captures the directions of parallel invariant lines of
systems
\eqref{2l1}.
\elm
We construct the  following $T$-comitants:
\bnot\label{not1}
\be\label{Comit:Bi}
\bal
B_3(a,x,y)&=(C_2,D)^{(1)}=Jacob\left( C_2,D\right),\\
B_2(a,x,y)&=\left(B_3,B_3\right)^{(2)} - 6B_3(C_2,D)^{(3)},\\
B_1(a)&=\Res_x\left( C_2,D\right)/y^9=-2^{-9}3^{-8}\left(B_2,B_3\right)^{(4)}.\\
\eal
\ee
\enot
\blm\label{lm:BGI}\mbox{\rm\!\!\cite{Dana_Vlp1}}  For the existence of an
invariant straight line in one (respectively 2, 3 distinct )
directions in the affine plane it is necessary that $B_1=0$
(respectively $B_2=0$, $B_3=0$).
\elm

Let us apply a translation $x=x'+x_0$, $y=y'+y_0$ to the polynomials $p(a,x,y)$
and $q(a,x,y)$. We obtain $ \tilde p(\tilde a(a,x_0,y_0),x',y')=p(a, x'+x_0,
y'+y_0),$ $\quad \tilde q(\tilde a(a,x_0,y_0),x',y')=q(a, x'+x_0, y'+y_0).$ Let
us construct the following polynomials
$$
\bal
\Gamma_i(a,x_0,y_0)& \equiv  \Res_{x'}
    \Big(C_i\big(\tilde a(a,x_0,y_0),x',y'\big),C_0\big(\tilde
    a(a,x_0,y_0),x',y'\big)\Big)/(y')^{i+1},\\
             &   \Gamma_i(a,x_0,y_0) \in \R[a,x_0,y_0],\ (i=1,2).\\
\eal
$$
\bnot\label{not2}
\be\label{2l4a}
   \tilde{\cal E}_i(a,x,y)=\left.\Gamma_i(a,x_0,y_0)\right|_{\{x_0=x,\ y_0=y\}}\in \R[a,x,y]
    \ \ (i=1,2).
\ee
\enot
\bobs\label{obs_} We  note that the constructed polynomials
$\tilde{\cal E}_1(a,x,y)$ and $\tilde{\cal E}_2(a,x,y) $ are
affine comitants of systems \eqref{2l1} and are homogeneous
polynomials in the coefficients $a_{00},\ldots, b_{02}$ and
non-homogeneous in $x,y$ and $\
  \deg_{a} \tilde{\cal E}_1=3,\  \deg_{\,(x,y)} \tilde{\cal
  E}_1=5,\ \
    \deg_{a} \tilde{\cal E}_2=4,\  \deg_{\,(x,y)} \tilde{\cal E}_2=6.
$ \eobs
\bnot\label{GCD:Ei}
Let  ${\cal E}_i(a,X,Y,Z)$ $(i=1,2)$ be the homogenization of $\tilde{\cal
E}_i(a,x,y)$, i.e.
$$
{\cal E}_1(a,X,Y,Z)=Z^5\tilde {\cal E}_1(a,X/Z,Y/Z),\qquad {\cal
E}_2(a,X,Y,Z)=Z^6\tilde {\cal E}_1(a,X/Z,Y/Z)
$$
and $ \qquad
 {\cal H}(a,X,Y,Z)=\gcd\Big({\cal E}_1(a,X,Y,Z),\
 {\cal E}_2(a,X,Y,Z)\Big)$ in $\R[a,X,Y,Z]$.
\enot
\indent The geometrical meaning of  these affine comitants is
given by  the two following lemmas (see \cite{Dana_Vlp1}):
\blm\label{lm2}
 The straight line ${\cal L}(x,y)\equiv ux+vy+w=0$,\ $u,v,w\in \C$, $(u,v)\ne(0,0)$
is an invariant line for a quadratic system \eqref{2l1} if and
only if the polynomial ${\cal L}(x,y)$ is a common factor of  the
polynomials $\tilde{\cal E}_1(\ab,x,y)$ and $\tilde{\cal
E}_2(\ab,x,y)$ over $\C$,
 i.e.
$$
\tilde{\cal E}_i(\ab,x,y)=(ux+vy+w)\widetilde W_i(x,y)\quad
(i=1,2),
$$
where   $\widetilde W_i(x,y)\in \C[x,y].$
\elm
\blm\label{lm3} If ${\cal L}(x,y)\equiv ux+vy+w=0$, $u,v,w\in \C$, $(u,v)\ne(0,0)$ is an invariant
straight line of multiplicity $k$ for a quadratic  system
\eqref{2l1} then $[{\cal L}(x,y)]^k\mid \gcd(\tilde {\cal
E}_1,\tilde {\cal E}_2)$ in  $\R[x,y]$, i.e. there exist
$W_i(\ab,x,y)\in \C[x,y]$\ $(i=1,2)$ such that
\be
\tilde {\cal E}_i(\ab,x,y)= (ux+vy+w)^k W_i(\ab,x,y),\quad i=1,2.
\label{2l5}
\ee
\elm
\bcrl\label{Mult:Z=0}
 If the line $l_\infty:Z=0$ is of multiplicity $k>1$ then
$Z^{k-1}\mid \gcd({\cal E}_1, {\cal E}_2)$.
\ecrl

   Let us consider the
following $GL$-comitants of systems \eqref{2l1}:
\bnot\label{not3}
$$
\ba{ll}
   M(a,x,y) = 2\,\Hess\big(C_2(x,y)\big), &
    \eta(a) = \Discriminant\big(C_2(x,y)\big),\\
   K(a,x,y) = \Jacob\big(p_2(x,y),q_2(x,y)\big),\qquad &
   \mu(a) =  \Discriminant\big(K(a,x,y)\big),\\
   N(a,x,y) =  K(a,x,y) + H(a,x,y), &
   \theta(a)  =   \Discriminant\big(N(a,x,y)\big),
\ea
$$
\enot
\noindent the geometrical meaning of which is revealed by the next
3 lemmas (see \cite{Dana_Vlp1}).
\blm\label{rm2_} Let $S\in \QS$ and let $\ab\in \R^{12}$ be its
12-tuple of coefficients.  The common points of $P=0$ and $Q=0$ on
the line $Z=0$ are given by the common linear factors over $\C$ of
$p_2$ and $q_2$. Moreover,
\beq
\deg\gcd(p_2(x,y),q_2(x,y)) =\left\{\begin{array}{lcl}
           0 & iff &\mu(\ab)\ne0;\\
           1 & iff &\mu(\ab)=0,\ K(\ab,x,y)\not=0;\\
           2 & iff &K(\ab,x,y)=0.
                         \end{array}\right.
\eeq
\elm
\blm\label{lm4}
 A necessary condition for  the existence of one
couple (respectively, two couples) of parallel invariant straight
lines of a systems \eqref{2l1} corresponding to $\ab\in\R^{12}$ is
the condition $\theta(\ab) =0$ (respectively, $N(\ab,x,y)=0$).
\elm
\blm\label{lm_3:2}
The type of the divisor $D_S(C,Z)$ for systems (\ref{il1}) is
determined by the corresponding conditions indicated in Table 1,
where we write $\omega_1^c+\omega_2^c+\omega_3$ if two of the
points, i.e. $\omega_1^c, \omega_2^c$, are complex but not real.
        Moreover, for each type of the divisor $D_S(C,Z)$ given
by Table 1 the quadratic systems (\ref{il1}) can be brought via a
linear transformation to one of the following  canonical systems
$(\SSS_{I})-(\SSS_{V})$ corresponding to their behavior at
infi\-ni\-ty.
\elm
\smallskip

In order to determine the existence of a common factor of the
polynomials ${\cal E}_1(\ab,X,Y,Z)$ and ${\cal E}_2(\ab,X,Y,Z)$ we
shall use the notion of the  resultant   of two polynomials with
respect to a given indeterminate (see for instance,
\cite{Walker}).

Let us  consider two polynomials $f,g\in R[x_1,x_2,\ldots,x_r]$
where $R$ is a unique factorization domain. Then we can regard the
polynomials $f$ and  $g$  as polynomials in $x_r$ over the ring
$R[x_1,x_2,\ldots,x_{r-1}]$, i.e.
$$
\bal
&f(x_1,x_2,\ldots,x_r)=a_0+a_1x_r+\ldots+a_nx_r^n,\\
&g(x_1,x_2,\ldots,x_r)=b_0+a_1x_r+\ldots+b_mx_r^m.
\eal
$$
\blm\label{Trudi:2}{\rm \cite{Walker}} Assuming $a_nb_m\ne0$ and
$n,m>0$, the resultant $\Res_{x_r}(f,g)$ of the polynomials $f$
and $g$ with respect to $x_r$ is a polynomial in
$R[x_1,x_2,\ldots,x_{r-1}]$ which is zero if and only if $f$ and
$g$ have a common factor involving $x_r$.
 \elm

\begin{table}[!htb]
\begin{center}
\begin{tabular}{|c|c|c|c|}
\multicolumn{3}{r}{\bf Table  1}\\[1mm]
\hline
  \raisebox{-0.7em}[0pt][0pt]{Case}  & \raisebox{-0.7em}[0pt][0pt]{Type of $D_S(C,Z)$}
      & Necessary and sufficient   \\[-1mm]
         & & conditions on the comitants \\
 \hline\hline
 \rule{0pt}{1.2em} $1$ & $\omega_1+\omega_2+\omega_3 $ &  $\eta>0 $ \\[1mm]
\hline
 \rule{0pt}{1.2em}$2$  & $\omega_1^c+\omega_2^c+\omega_3 $ &  $\eta<0$ \\[1mm]
\hline
 \rule{0pt}{1.2em}  $3$ & $2\omega_1+\omega_2 $ &  $\eta=0,\quad M\ne0$ \\[1mm]
\hline
 \rule{0pt}{1.2em} $4$ & $3\omega $ &  $ M=0,\quad C_2\ne0$ \\[1mm]
\hline
 \rule{0pt}{1.2em} $5$ & $D_S(C,Z)$ undefined  &  $ C_2=0$ \\[1mm]
\hline
\end{tabular}
\end{center}
\end{table}

\begin{table}[!htb]
$$
\bal
&\left\{\ba{rcl}
 \displaystyle \frac{dx}{dt}&=&k+cx+dy+gx^2+(h-1)xy,\\[2mm]
 \displaystyle \frac{dy}{dt}&=& l+ex+fy+(g-1)xy+hy^2;
\ea\right. &\qquad (\SSS_I)\\[3mm]
&\left\{\ba{rcl}
 \displaystyle \frac{dx}{dt}&=&k+cx+dy+gx^2+(h+1)xy,\\[2mm]
 \displaystyle \frac{dy}{dt}&=& l+ex+fy-x^2+gxy+hy^2;
\ea\right.& \hspace{2cm}(\SSS_{I\!I})\\[3mm]
&\left\{\ba{rcl}
 \displaystyle \frac{dx}{dt}&=&k+cx+dy+gx^2+hxy,\\[2mm]
 \displaystyle \frac{dy}{dt}&=& l+ex+fy+(g-1)xy+hy^2;
\ea\right.&\qquad (\SSS_{I\!I\!I})\\[3mm]
&\left\{\ba{rcl}
 \displaystyle \frac{dx}{dt}&=&k+cx+dy+gx^2+hxy,\\[2mm]
 \displaystyle \frac{dy}{dt}&=& l+ex+fy-x^2+gxy+hy^2,
\ea\right.&\qquad (\SSS_{I\!V})\\[3mm]
&\left\{\ba{rcl}
 \displaystyle \frac{dx}{dt}&=&k+cx+dy+x^2,\\[2mm]
 \displaystyle \frac{dy}{dt}&=& l+ex+fy+xy.
\ea\right.&\qquad (\SSS_{V})
\eal
$$
\end{table}


\newpage

\section{ The Main Theorem } \label{Sec:m_il:5}

\bnot
We denote by $\QSL_{\bf4}$  the class of  all  quadratic
differential systems
\eqref{2l1} with $(p,q)=1$  possessing a configuration of 4 invariant straight
lines including the line at infinity and including possible
multiplicities and the line at infinity does not consist entirely
of singularities.
\enot
\bobs The case when the line at infinity  is a union of
singularities will be discussed in a forthcoming paper.
 \eobs

\blm\label{lm_NB2} If  a quadratic system $(S)$ corresponding to $\ab\in\R^{12}$
 belongs to the class $\QSL_{\bf4}$, then for this system
one of the  following two conditions are satisfied in $\R[x,y]$:\\[-5mm]
$$
\bal
&(i)\quad \theta\ne0,\ B_3(\ab,x,y)=0;\qquad (ii)\quad
\theta(\ab)=0= B_2(\ab,x,y).
\eal
$$
\elm
\BProof Indeed, if for a system \eqref{2l1} the condition $M_{{}_\IL}=4$ is
satisfied then taking into account the Definition \ref{def:multipl} we conclude
that there exists a perturbation  of the coefficients of the
 system~\eqref{2l1} within the class of quadratic systems such that
the perturbed systems have $4$ distinct invariant lines (real or
complex, including the  line $Z=0$). Hence, the   perturbed
systems must possess either three affine lines with  distinct
directions or one couple of parallel lines and another  line in a
different direction.  Then, by continuity and according to Lemmas
\ref{lm:BGI} and \ref{lm4} we respectively have
 either conditions $(i)$ or $(ii)$. \EProof

 By Lemmas \ref{lm2}, \ref{lm3}  and  \ref{lm4} we
obtain the following result:
\blm \label{gcd:3}
(a) If for a system $(S)$ of coefficients $\ab\in \R^{12}$,\
$M_{{}_\IL}=4$ then\\\hphantom{m} $\deg\gcd\big({\cal
E}_1(\ab,X,Y,Z), {\cal E}_2(\ab,X,Y,Z)\big)=3;$ \hphantom{mm} (b)
If $\theta(\ab)\ne0$ then $M_{{}_\IL}\le4$.
\elm
We shall use here the following $T$-comitants constructed in
\cite{Dana_Vlp1}
\beq
&& H_1(a)= -\big((C_2,C_2)^{(2)},C_2)^{(1)},D\big)^{(3)},\\
&& H_2(a,x,y)=\left(C_1,\ 2H-N\right)^{(1)}-2D_1N,\\
&& H_3(a,x,y)=(C_2,D)^{(2)},\\
&&H_4(a)=\big((C_2,D)^{(1)},(C_2,D_2)^{(1)}\big)^{(2)},\\
&&H_5(a)=\big((C_2,C_2)^{(2)},(D,D)^{(2)}\big)^{(2)}+
    8\big((C_2,D)^{(2)},(D,D_2)^{(1)}\big)^{(2)},\\
&& H_6(a,x,y)= 16N^2(C_2,D)^{(2)}+ H_2^2(C_2,C_2)^{(2)}
\eeq
and $CT$-comitants
\beq
&&N_1(a,x,y)=C_1(C_2,C_2)^{(2)} -2C_2(C_1,C_2)^{(2)} ,\\
&&N_2(a,x,y)=D_1(C_1,C_2)^{(2)}-\left((C_2,C_2)^{(2)},C_0\right)^{(1)},\\
&&N_3=(a,x,y)   \left(C_2,C_1\right)^{(1)},\\
&& N_4(a,x,y)=   4\left(C_2,C_0\right)^{(1)}- 3C_1D_1,\\
&& N_5(a,x,y)= \big((D_2,C_1)^{(1)} + D_1D_2\big)^2
-4\big(C_2,C_2\big)^{(2)}\big(C_0,D_2\big)^{(1)},\\
&&N_6(a,x,y)= 8D+C_2\left[8(C_0,D_2)^{(1)}-3(C_1,C_1)^{(2)}+2D_1^2\right].\\
\eeq
We  shall also use the following remark:
\brm\label{rem:transf} Assume $s,\, \gamma\in \R$, $\gamma>0$. Then the transformation
$x=\gamma^{s}x_1$, $y=\gamma^{s}y_1$ and \mbox{$t=\gamma^{-s}t_1$}
does not change the coefficients of the  quadratic part of a
quadratic system, whereas each coefficient of the linear
(respectively,
 constant ) part will be multiplied by  $\gamma^{-s}$ (respectively,
by $\gamma^{-2s}$).
\erm
\bth\label{th_mil_4}
(i) The class $\QSL_{\bf4}$ splits into 46 distinct subclasses
indicated in {\bf Diagram 1} with the corresponding Configurations
4.1-4.46 where the complex invariant straight lines are indicated
by dashed lines. If an invariant straight line has multiplicity
$k>1$, then the number $k$ appears near the corresponding straight
line and this line is in bold face. We indicate next to the
singular points their multiplicities as follows:
$\left(I_\omega(p,q)\right)$ if $\omega$ is a finite singularity,
$\left(I_\omega(C,Z),\ I_\omega(P,Q)\right)$ if $\omega$ is an
infinite singularity with $I_w(P,Q)\ne0$ and
$\left(I_\omega(C,Z)\right)$ if $\omega$ is an infinite
singularity with $I_\omega(P,Q)=0$.

\smallskip
(ii) We consider the orbits of the class $\QSL_{\bf4}$ under  the
action of the affine group and time rescaling. The systems of the
form {\sl(IV.1)} up to form {\sl(IV.46)} from the Table 2 form a
system of representatives of these orbits under this action. A
differential system $(S)$ in $\QSL_{\bf4}$ is in the orbit of   a
system belonging to $(IV.i)$ if and only if the corresponding
conditions in the middle column (where the polynomials $H_i$
$(i=7,\ldots, 11)$ are   $T$-comitants  to be introduced below) is
verified for this system $(S)$.  The conditions indicated in the
middle column are affinely invariant.

\smallskip

Wherever we have a case with invariant straight lines of
multiplicity $>1$  we indicate the corresponding perturbed systems
in the Table 3.
\eth

\begin{figure}
\centerline{\hfill{\bf Diagram 1
}\hphantom{999999999999999999999999}}
\centerline{\psfig{figure=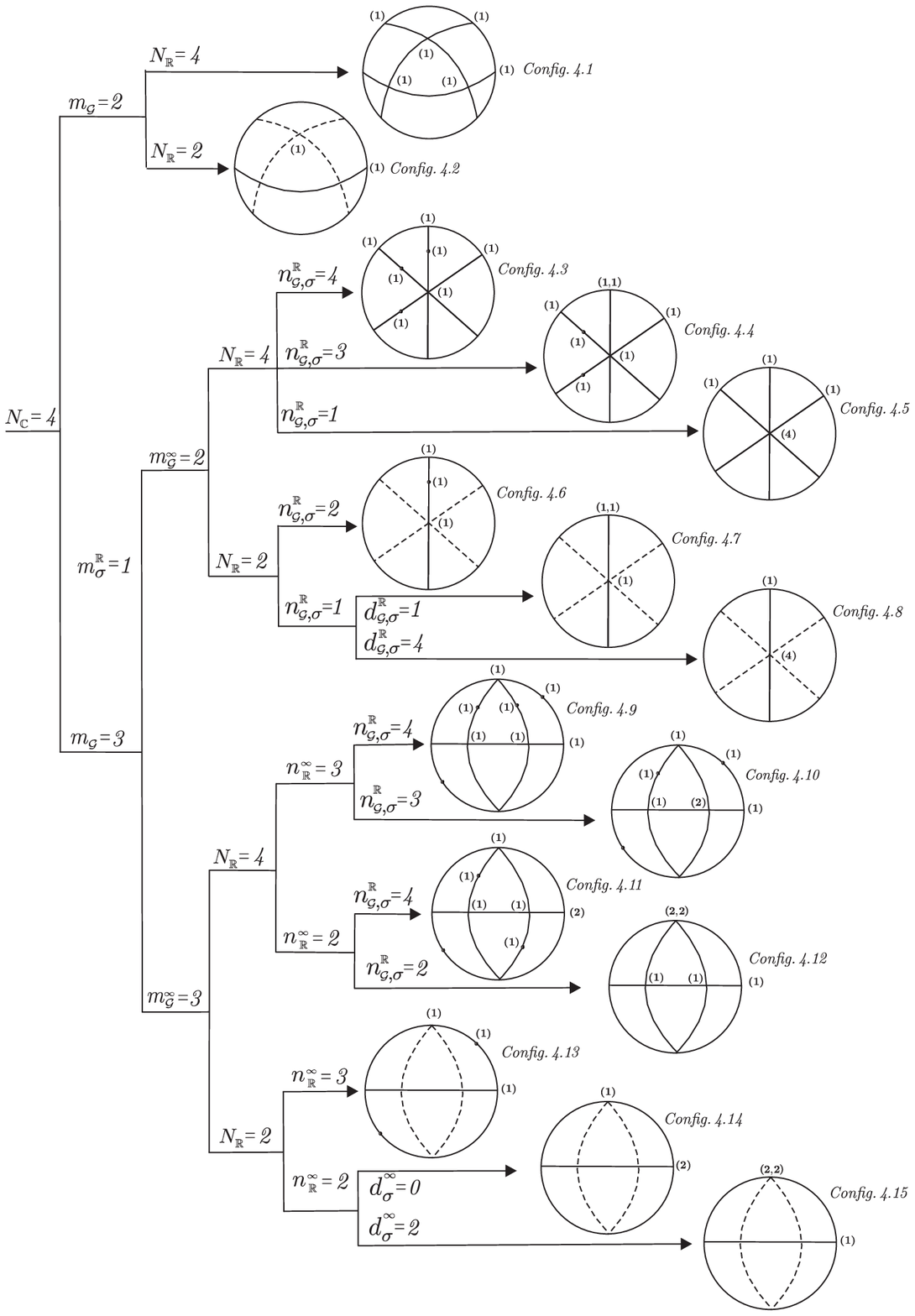}}
\end{figure}
\begin{figure}
\centerline{\hfill{\bf Diagram 1 } {\it
(continued)}\hphantom{999999999 }}
\vspace{3mm}
\centerline{ \psfig{figure=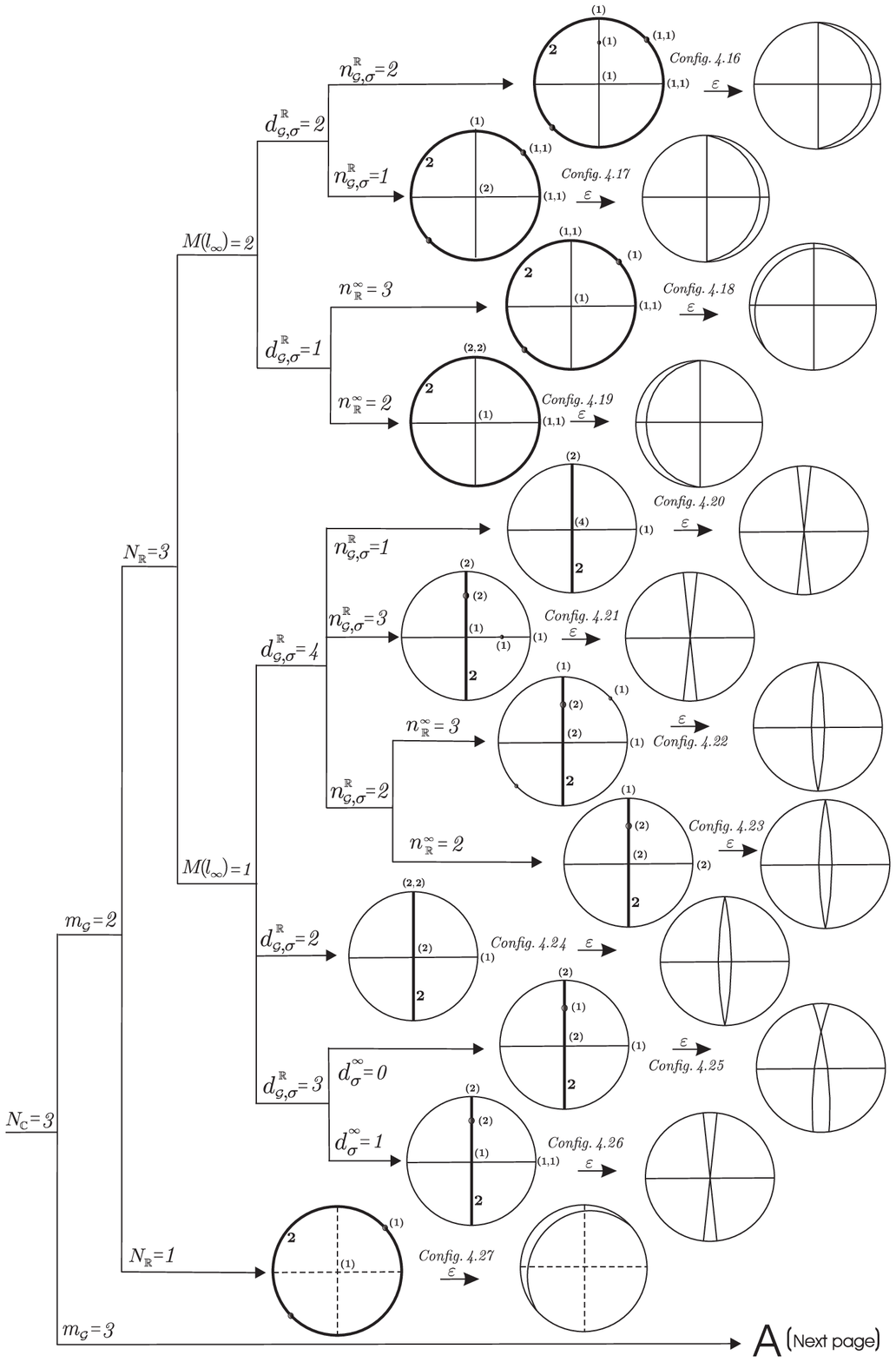}}
\end{figure}
\begin{figure}
\centerline{\hfill{\bf Diagram 1} {\it
(continued)}\hphantom{99999999999}}
\vspace{3mm}
\centerline{\psfig{figure=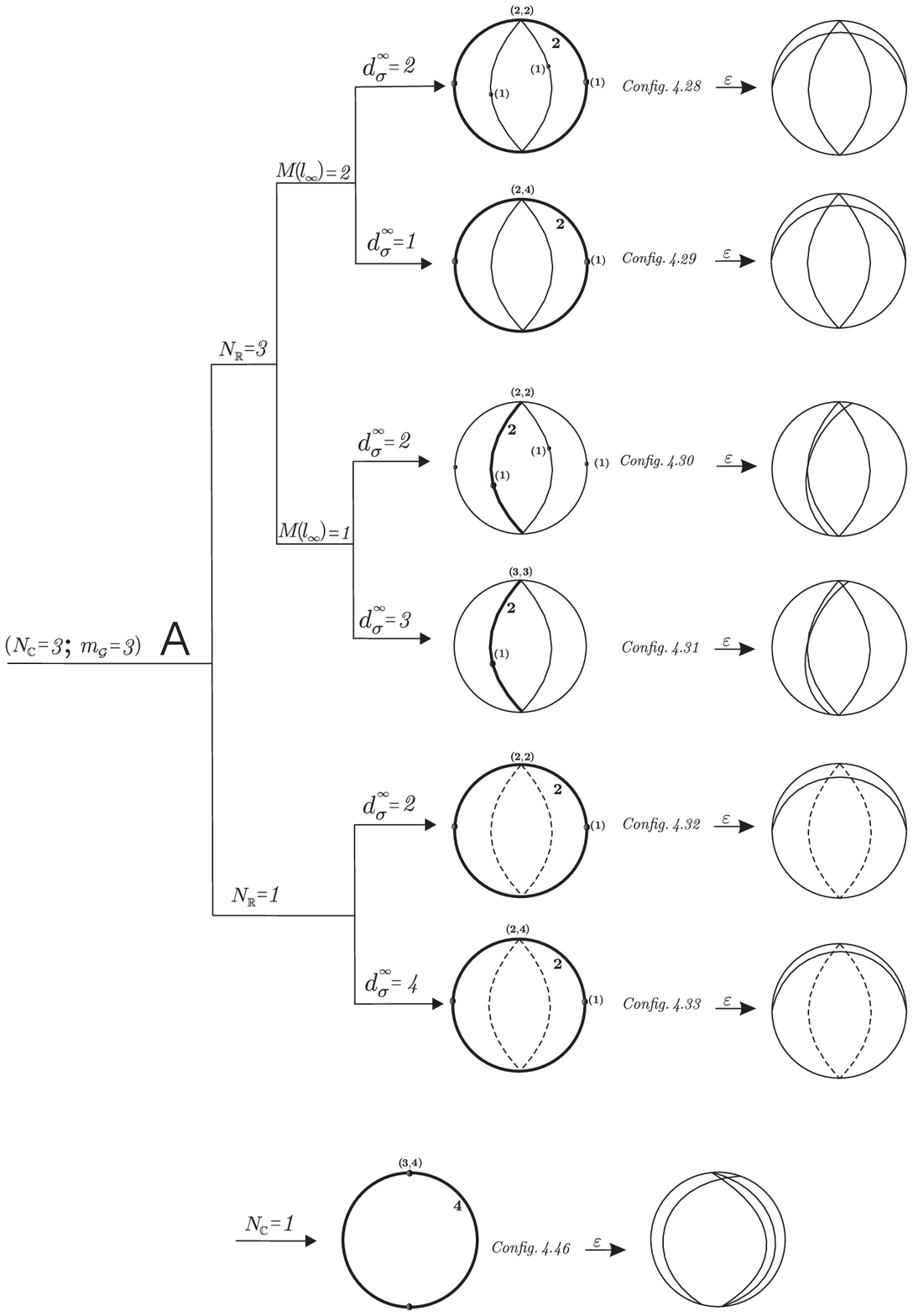}\hphantom{9999}}
\end{figure}
\begin{figure}
\centerline{\hfill{\bf Diagram 1} {\it
(continued)}\hphantom{99999999999999999}}
\vspace{3mm}
\centerline{\psfig{figure=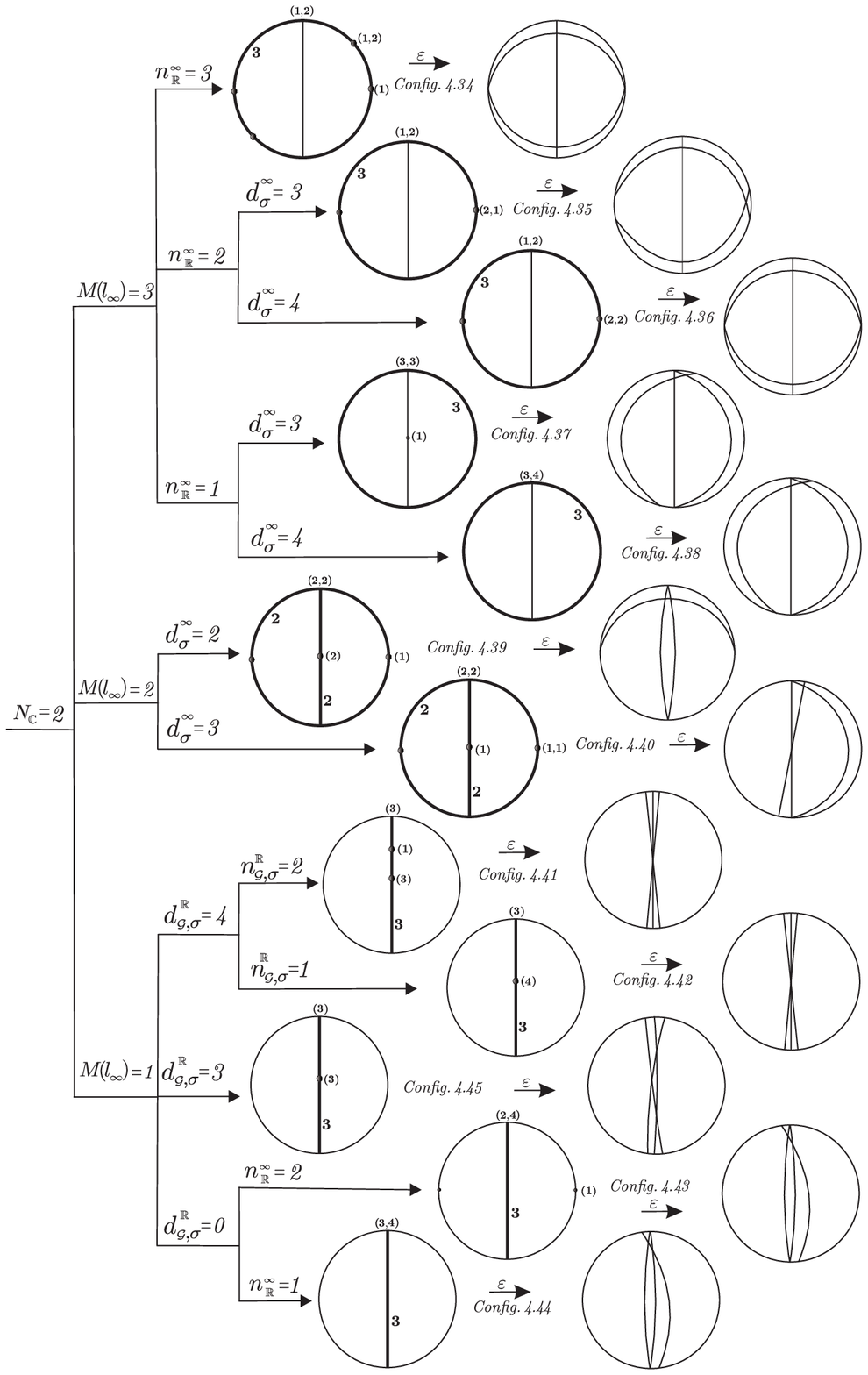}\hphantom{9999}}
\end{figure}

\begin{table}[!htb]
\begin{center}
\begin{tabular}{|l|c|c|c|}
\multicolumn{4}{r}{\bf Table 2}\\[1mm]
\hline \raisebox{-0.7em}[0pt][0pt]{\qquad Orbit representative}  &
\multicolumn{2}{c|}{Necessary and sufficient}
           & \!\!\raisebox{-0.7em}[0pt][0pt]{Configuration}\!\! \\[-1mm]
                             & \multicolumn{2}{c|}{conditions} & \\
 \hline\hline\rule[0mm]{0mm}{10.0mm}
 \!\!\!({\it IV\!.1})\ $\left\{\!\!\ba{l} \dot x= lx+lx^2+(b-1)xy, \\[-0.8mm]
                \dot y=-by+(l\!-\!1)xy+by^2,\\[-0.8mm]
                b,l\in\R,\ bl(b+l-1)\ne0,  \\[-0.8mm]
                (b-1)(l-1)(b+l)\ne0 \ea\!\!\right.$ & \multicolumn{2}{c|}{ $ \ba{c}\eta>0,\
                \theta\ne0,\  B_3=0,\\  H_7\ne0\ea $} & Config.\ 4.1\\
\hline\rule[0mm]{0mm}{10.0mm}
 \!\!\!({\it IV\!.2})\ $\left\{\!\!\ba{l} \dot x= lx^2+(b+1)xy,\ b,l\in\!\R  \\[-0.8mm]
                \dot y=b[l^2+(b+1)^2] +2lby \\[-0.8mm]
                +(l^2\!+\!1\!-\!h^2)x\! -\!x^2\!+ \!lxy\!+ \!by^2\\[-0.8mm]
                 b(b+1)\big[l^2+(b-1)^2\big]\ne0\ea\!\!\right.$ & \multicolumn{2}{c|}{$ \ba{c}\eta<0,\
                \theta\ne0,\  B_3=0,\\   H_7\ne0\ea $} & Config.\ 4.2\\
\hline\rule[0mm]{0mm}{10.0mm}
 \!\!\!({\it IV\!.3})\ $\left\{\!\!\ba{l} \dot x= x+lx^2+(b-1)xy, \\[-0.8mm]
                \dot y=y+(l-1)xy+by^2,\\[-0.8mm]
                b,l\in\R,\ bl(b+l-1)\ne0,  \\[-0.8mm]
                 (b-1)(l-1)(b+l)\ne0  \ea\!\!\right. $ &  \multicolumn{2}{c|}{$ \ba{c}\eta>0,\
                \theta\ne0,\, B_3=0,\\  H_7=0,\,H_1\ne0,\,\mu\ne0\ea $} & Config.\ 4.3\\
\hline\rule[0mm]{0mm}{8.0mm}
 \!\!\!({\it IV\!.4})\ $\left\{\!\!\ba{l} \dot x=x+ lx^2-xy, \\[-0.8mm]
                \dot y=y+(l-1)xy,\\[-0.8mm]
                l\in \R,\ l(l-1)\ne0  \ea\!\!\right. $ &   \multicolumn{2}{c|}{$ \ba{c}\eta>0,\
                \theta\ne0,\, B_3=0,\\  H_7=0,\,H_1\ne0,\,\mu=0\ea $} & Config.\ 4.4\\
\hline\rule[0mm]{0mm}{10.0mm}
 \!\!\!({\it IV\!.5})\ $\left\{\!\!\ba{l} \dot x= lx^2+(b-1)xy, \\[-0.8mm]
                \dot y=(l-1)xy+by^2,\\[-0.8mm]
                 b,l\in\R,\ bl(b+l-1)\ne0,  \\[-0.8mm]
                 (b-1)(l-1)(b+l)\ne0\ea\!\!\right.$ &  \multicolumn{2}{c|}{$ \ba{c}\eta>0,\
                \theta\ne0,\, B_3=0,\\  H_7=0,\,H_1=0\ea $} & Config.\ 4.5\\
\hline\rule[0mm]{0mm}{10.0mm}
 \!\!\!({\it IV\!.6})\ $\left\{\!\!\ba{l} \dot x= lx^2+(b+1)xy,\ b,l\in\!\R  \\[-0.8mm]
                \dot y=-1 +lx +(b-1)y \\[-0.8mm]
                -x^2+ lxy+ by^2\\[-0.8mm]
                 b(b+1)\big[l^2+(b-1)^2\big]\ne0\ea\!\!\right.$ & \multicolumn{2}{c|}{$ \ba{c}\eta<0,\
                \theta\ne0,\  B_3=0,\\ H_7=0,\, \mu\ne0,\, H_{9}\ne0\ea $} & Config.\ 4.6\\
\hline\rule[0mm]{0mm}{6.0mm}
 \!\!\!({\it IV\!.7})\ $\left\{\!\!\ba{l} \dot x= lx^2+ xy,\quad l\in\!\R  \\[-0.8mm]
                \dot y= -1+lx-y-x^2+ lxy\ea\!\!\right.$ & \multicolumn{2}{c|}{$ \ba{c}\eta<0,\
                \theta\ne0,\  B_3=0,\\ H_7=0,\, \mu=0\ea $} & Config.\ 4.7\\
\hline\rule[0mm]{0mm}{8.0mm}
 \!\!\!({\it IV\!.8})\ $\left\{\!\!\ba{l} \dot x= lx^2+(b+1)xy,\ b,l\in\!\R  \\[-0.8mm]
                \dot y= -x^2+ lxy+ by^2\\[-0.8mm]
                 b(b+1)\big[l^2+(b-1)^2\big]\ne0\ea\!\!\right.$ & \multicolumn{2}{c|}{$ \ba{c}\eta<0,\
                \theta\ne0,\  B_3=0,\\ H_7=0,\, \mu\ne0,\, H_{9}=0\ea $} & Config.\ 4.8\\
\hline\rule[0mm]{0mm}{10.0mm}
  \!\!\!({\it IV\!.9})\ $\left\{\!\!\ba{l} \dot x= l(x^2-1),\ l,b\in\R, \\[-0.8mm]
                \dot y=(y\!+\!b)[y\!+\!(l\!-\!1)x\!-\!b],\\[-0.0mm]
                 l(l\!-\!1)\!\ne\!0,\  (l\!-\!1)^2\ne4b^2,\\[-0.8mm]
                 (l+1)^2\ne4b^2 \ea\!\!\right.$  &\!\!\!\! $
                 \ba{c}\eta>0,\,
                 \theta=B_2=0\\ \mu\, B_3\,H_4\ne0,\\H_7=0,H_9\ne0\ea $\!\!\!\!
                 &\!\!\!$\ba{c}   H_{10}N>0,\\[-0.8mm]
                 \!\!\!\rule[0.5mm]{0.16 \textwidth}{0.15mm}\!\!\!\!\\ \!\!\!N\!=\!0,H_8\!>\!0 \!\!\!\ea $\!\!\!
                      & Config.\ 4.9 \\
\hline\rule[0mm]{0mm}{9.0mm}
  \!\!\!\!({\it IV\!.10})\ $\left\{\!\!\ba{l} \dot x= (2l+1)(x^2\!-\!1), \\[-0.8mm]
                \dot y=(y+l)(y+2lx-l),\\[-0.0mm]
                 l\!\in\!\R,\ l(2l+1)\ne0 \ea\!\!\right.$  &\!\!\!\! $
                 \ba{c}\eta>0,\,
                 \theta=B_2=0\\ \mu\, B_3\,H_4\ne0,\\H_7=H_9=0\ea\!\!\!\! $
                 &\!\!\!$\ba{c}   H_{10}N>0,\\[-0.9mm]
                 \!\!\!\rule[0.5mm]{0.16 \textwidth}{0.15mm}\!\!\!\!\\[-0.8mm] \!\!\!N\!=\!0,H_8\!>\!0 \!\!\!\ea $\!\!\!
                      & Config.\ 4.10 \\
\hline\rule[0mm]{0mm}{6.0mm}
 \!\!\!({\it IV\!.11})\ $\left\{\!\!\ba{l} \dot x= x^2+ xy,\ l\in\!\R,\ l\ne\pm1  \\[-0.8mm]
                \dot y= (y+l)^2-1\ea\!\!\right.$ & \multicolumn{2}{c|}{$
                \ba{c}\eta=0,\ M\ne0,\
                \theta=B_2=0,\\  B_3\mu \ne0,\, H_7=0,\, H_{10}>0 \ea $} & Config.\ 4.11\\
\hline\rule[0mm]{0mm}{7.0mm}
 \!\!\!({\it IV\!.12})\ $\left\{\!\!\ba{l} \dot x= l\big[(x+ b)^2-1\big],\\[-0.8mm]
                \dot y= (l-1)xy,\ l,b\in\!\R,\\[-0.8mm]
                l(l-1)(b^2-1)\ne0\ea\!\!\right.$ & \multicolumn{2}{c|}{$
                \ba{c}\eta=0,\ M\ne0,\
                \theta=B_3=0,\\  KH_6\ne0,\,\mu\!=\!H_7\!=\!0,\, H_{11}\!>\!0 \ea $} & Config.\ 4.12\\
 \hline\rule[0mm]{0mm}{9.0mm}
  \!\!\!({\it IV\!.13})\ $\left\{\!\!\ba{l} \dot x= l(x^2+1), \\[-0.8mm]
                \dot y=(y\!+\!b)[y\!+\!(l\!-\!1)x\!-\!b],\\[-0.0mm]
                  l,b\in\R,\ l(l-1)\ne0 \ea\!\!\right.$  &\!\!\!\! $
                 \ba{c}\eta>0,\,
                 \theta=B_2=0\\ \mu\, B_3\,H_4\ne0,\\H_7=0\ea \!\!\!\!$
                 &\!\!\!$\ba{c}   H_{10}N<0,\\[-0.9mm]
                 \!\!\!\rule[0.5mm]{0.16 \textwidth}{0.15mm}\!\!\!\\[-0.9mm] \!\!\!N\!=\!0,H_8\!<\!0 \!\!\!\ea $\!\!\!
                      & Config.\ 4.13 \\
\hline\rule[0mm]{0mm}{7.0mm}
 \!\!\!({\it IV\!.14})\ $\left\{\!\!\ba{l} \dot x= x^2+ xy,\ l\in\!\R  \\[-0.8mm]
                \dot y= (y+l)^2+1\ea\!\!\right.$ & \multicolumn{2}{c|}{$
                \ba{c}\eta=0,\ M\ne0,\
                \theta=B_2=0,\\  B_3\mu \ne0,\, H_7=0,\, H_{10}<0 \ea $} & Config.\ 4.14\\
\hline
\end{tabular}
\end{center}
\end{table}
\begin{table}[!htb]
\begin{center}
\begin{tabular}{|l|c|c|c|}
\multicolumn{4}{r}{\bf Table 2}{\it(continued)}\\[1mm]
\hline \raisebox{-0.7em}[0pt][0pt]{\qquad Orbit representative}  &
\multicolumn{2}{c|}{Necessary and sufficient}
           & \!\!\raisebox{-0.7em}[0pt][0pt]{Configuration}\!\! \\[-1mm]
                             & \multicolumn{2}{c|}{conditions} & \\
 \hline
\hline\rule[0mm]{0mm}{8.0mm}
 \!\!\!({\it IV\!.15})\ $\left\{\!\!\ba{l} \dot x= l\big[(x+ b)^2+1\big],\\[-0.8mm]
                \dot y= (l-1)xy,\ l,b\in\!\R,\\[-0.8mm]
                l(l-1)\ne0\ea\!\!\right.$ & \multicolumn{2}{c|}{$
                \ba{c}\eta=0,\ M\ne0,\
                \theta=B_3=0,\\  KH_6\ne0,\,\mu\!=\!H_7\!=\!0,\, H_{11}\!<\!0 \ea $} & Config.\ 4.15\\
\hline\rule[0mm]{0mm}{6.0mm}
 \!\!\!({\it IV\!.16})\ $\left\{\!\!\ba{l} \dot x= l+x,\ \ l\in\R, \\[-0.8mm]
                \dot y=y(y-x),\ l(l\!-\!1)\ne0  \ea\!\!\right.$ & \multicolumn{2}{c|}{$ \ba{c}\eta>0,\
                \theta= B_2=\mu=0, \\B_3\ne0,\,  H_7=0,\,H_9\ne0\ea $} & Config.\ 4.16\\
\hline\rule[0mm]{0mm}{6.0mm}
 \!\!\!({\it IV\!.17})\ $\left\{\!\!\ba{l} \dot x= x,\\[-0.8mm]
                \dot y=y(y-x)  \ea\!\!\right.$ & \multicolumn{2}{c|}{$ \ba{c}\eta>0,\
                \theta= B_2=\mu=0, \\B_3\ne0,\,  H_7=H_9=0,\,H_{10}\ne0\ea $} & Config.\ 4.17\\
\hline\rule[0mm]{0mm}{6.0mm}
\! \!\!\!({\it IV\!.18})\ $\left\{\!\!\ba{l} \dot x= l(l\!+\!1)+lx+y,\ l\in\!\R, \\[-0.8mm]
                \dot y=y(y-x),\ l(l+1)\ne0  \ea\!\!\right.$ & \multicolumn{2}{c|}{$ \ba{c}\eta>0,\
                \theta= B_3=\mu=0, \\  H_7\ne0\ea $} & Config.\ 4.18\\
\hline\rule[0mm]{0mm}{6.0mm}
 \!\!\!({\it IV\!.19})\ $\left\{\!\!\ba{l} \dot x= l+x,\quad l\in\!\R,\\[-0.8mm]
                \dot y= -xy,\ \   l(l-1)\ne0\!\ea\!\!\right.$ & \multicolumn{2}{c|}{$
                \ba{c}\eta=0,\ M\ne0,\
                \theta=B_3=K=0,\\  NH_6\ne0,\,\mu\!=\!H_7\!=\!0,\, H_{11}\!\ne\!0 \!\!\ea $} & Config.\ 4.19\\
\hline\rule[0mm]{0mm}{6.0mm}
 \!\!\!({\it IV\!.20})\ $\left\{\!\!\ba{l} \dot x= lx^2+ xy,\ l(l\!-\!1)\!\ne\!0  \\[-0.8mm]
                \dot y= (l-1)xy+y^2,\ l\in\!\R\ea\!\!\right.$ & \multicolumn{2}{c|}{$
                \ba{c}\eta=0,\ M\ne0,\
                \theta\ne0,\\  B_3=H_7= D=0 \ea $} & Config.\ 4.20\\
\hline\rule[0mm]{0mm}{6.0mm}
 \!\!\!({\it IV\!.21})\ $\left\{\!\!\ba{l} \dot x= lx^2+ xy,\ l(l\!-\!1)\!\ne\!0  \\[-0.8mm]
                \dot y= (y+1)(lx-x+y),\ l\in\!\R\!\ea\!\!\right.$ & \multicolumn{2}{c|}{$
                \ba{c}\eta=0,\ M\ne0,\
                \theta\ne0,\\  B_3=H_7=0,\,D\ne0,\,\mu\ne0\ea $} & Config.\ 4.21\\
\hline\rule[0mm]{0mm}{6.0mm}
  \!\!\!({\it IV\!.22})\ $\left\{\!\!\ba{l} \dot x= lx^2,\  l\!\in\!\R,\ l(l\!-\!1)\!\ne\!0  \\[-0.8mm]
                \dot y=(y\!+\!1)[y\!+\!(l\!-\!1)x\!-\!1],\\[-0.0mm]
                 \ea\!\!\right.$ &\!\!\!\! $
                 \ba{c}\eta>0,\,
                 \theta=B_2=0\\ \mu\, B_3\,H_4\ne0,\\H_7=0\ea\!\!\!\! $
                 &\!\!\!$\ba{c}\!\! N\!\ne\!0,  H_{\!10}\!=\!0\!\!\\[-0.9mm]
                 \!\!\!\rule[0.5mm]{0.16 \textwidth}{0.15mm}\!\!\!\\[-0.9mm] \!\!\!N\!=\!0,H_8\!=\!0 \!\!\!\ea $\!\!\!
                      & Config.\ 4.22 \\
\hline\rule[0mm]{0mm}{6.0mm}
 \!\!\!({\it IV\!.23})\ $\left\{\!\!\ba{l} \dot x= x^2+ xy,\\[-0.8mm]
                \dot y= (y+1)^2\ea\!\!\right.$ & \multicolumn{2}{c|}{$
                \ba{c}\eta=0,\ M\ne0,\
                \theta=B_2=0,\\  B_3\mu\!\ne\!0,\, H_7\!= \!0,\, H_{10}\!=\!0 \!\!\ea $} & Config.\ 4.23\\
\hline\rule[0mm]{0mm}{6.0mm}
 \!\!\!({\it IV\!.24})\ $\left\{\!\!\ba{l} \dot x= l (x+ 1)^2,\ l\in\!\R,\\[-0.8mm]
                \dot y= (l-1)xy,\
                l(l\!-\!1)\!\ne0\!\ea\!\!\right.$ & \multicolumn{2}{c|}{$
                \ba{c}\eta=0,\ M\ne0,\
                \theta=B_3=0,\\  KH_6\!\ne\!0,\,\mu\!=\!H_7\!=\!0,\, H_{11}\!=\!0\!\! \ea $} & Config.\ 4.24\\
\hline\rule[0mm]{0mm}{6.0mm}
 \!\!\!({\it IV\!.25})\ $\left\{\!\!\ba{l} \dot x= lx^2+ xy,\ l\in\!\R,\ l(l\!-\!1)\!\ne\!0  \\[-0.8mm]
                \dot y= y+(l-1)xy+y^2\ea\!\!\right.$ & \multicolumn{2}{c|}{$
                \ba{c}\eta=0,\ M\ne0,\
                \theta\ne0,\\  B_3=0,\, H_7\ne0\ea $} & Config.\ 4.25\\
\hline\rule[0mm]{0mm}{6.0mm}
 \!\!\!({\it IV\!.26})\ $\left\{\!\!\ba{l} \dot x=  xy,\ \\[-0.8mm]
                \dot y= (y+1)(y-x)\ea\!\!\right.$ &  \multicolumn{2}{c|}{$
                \ba{c}\eta=0,\ M\ne0,\
                \theta\ne0,\\  B_3=H_7=0,\,D\ne0,\,\mu=0\ea $} & Config.\ 4.26\\
\hline\rule[0mm]{0mm}{6.0mm}
 \!\!\!({\it IV\!.27})\ $\left\{\!\!\ba{l} \dot x= 2lx+ 2y,\quad l\in\!\R  \\[-0.8mm]
                \dot y= l^2+1-x^2-y^2\ea\!\!\right.$ & \multicolumn{2}{c|}{$ \ba{c}\eta<0,\
                \theta=0,\  B_3=0,\\ N\ne0,\ H_7\ne0\ea $} & Config.\ 4.27\\
\hline\rule[0mm]{0mm}{6.0mm}
 \!\!\!({\it IV\!.28})\ $\left\{\!\!\ba{l} \dot x= x^2-1, \quad l\in\!\R,  \\[-0.8mm]
                \dot y= x+ly,\    l(l^2-4)\ne0 \!\ea\!\!\right.$ & \multicolumn{2}{c|}{$
                \ba{c}\eta=0,\, M\ne0,\,
                \theta\!=\!\mu\!=\!N\!=\!B_3\!=\!0,\\[-0.8mm]
                N_1N_2\ne0,\,  K\!=\!0,\, N_5\!>\!0,\,D\!\ne\!0\ea $} & Config.\ 4.28\\
\hline\rule[0mm]{0mm}{6.0mm}
 \!\!\!({\it IV\!.29})\ $\left\{\!\!\ba{l} \dot x= x^2-1, \quad l\in\!\R,  \\[-0.8mm]
                \dot y= l+x,\quad   l\ne\pm1 \!\ea\!\!\right.$ & \multicolumn{2}{c|}{$
                \ba{c}\eta=0,\, M\ne0,\,
                \theta\!=\!\mu\!=\!N\!=\!B_3\!=\!0,\\[-0.8mm]
                N_1N_2\ne0,\,  K\!=\!0,\, N_5\!>\!0,\,D\!=\!0\ea $} & Config.\ 4.29\\
\hline\rule[0mm]{0mm}{6.0mm}
 \!\!\!({\it IV\!.30})\ $\left\{\!\!\ba{l} \dot x= (1+x)(1+lx),\  \\[-0.8mm]
                \dot y= 1+(l-1)xy,\\[-0.8mm]
                 l\in\!\R,\ l(l^2-1)\ne0\!\ea\!\!\right.$ & \multicolumn{2}{c|}{$
                \ba{c}\eta=0,\ M\ne0,\
                \theta=H_6=0,\\  NB_3\ne0,\,\mu\!=0,\,K\!\ne\!0,\, H_{11}\!\ne\!0 \ea $} & Config.\ 4.30\\
\hline
 \!\!\!({\it IV\!.31})\ $\left\{\!\!\ba{l} \dot x= x +x^2, \quad  l\in\R,  \\[-0.8mm]
                \dot y= l\!-\!x^2\!+\!xy,\  l(l\!+\!1)\!\ne\!0 \!\ea\!\!\right.$ & \multicolumn{2}{c|}{$
                \ba{c}\eta=M=0,\
                \theta=B_3=0,\\[-0.8mm] N_6N\ne0,\ H_{11}\ne0\ea $} & Config.\ 4.31\\
\hline\rule[0mm]{0mm}{6.0mm}
 \!\!\!({\it IV\!.32})\ $\left\{\!\!\ba{l} \dot x= x^2+1, \quad l\in\!\R,  \\[-0.8mm]
                \dot y= x+ly,\quad l\ne0 \!\ea\!\!\right.$ & \multicolumn{2}{c|}{$
                \ba{c}\eta=0,\, M\ne0,\,
                \theta\!=\!\mu\!=\!N\!=\!B_3\!=\!0,\\[-0.8mm]
                N_1N_2\ne0,\,  K\!=\!0,\, N_5\!<\!0,\,D\!\ne\!0\ea $} & Config.\ 4.32\\
\hline\rule[0mm]{0mm}{6.0mm}
 \!\!\!({\it IV\!.33})\ $\left\{\!\!\ba{l} \dot x= x^2+1,   \\[-0.8mm]
                \dot y= l+x,\quad l\in\!\R\!\ea\!\!\right.$ & \multicolumn{2}{c|}{$
                \ba{c}\eta=0,\, M\ne0,\,
                \theta\!=\!\mu\!=\!N\!=\!B_3\!=\!0,\\[-0.8mm]
                N_1N_2\ne0,\,  K\!=\!0,\, N_5\!<\!0,\,D\!=\!0\ea $} & Config.\ 4.33\\
\hline\rule[0mm]{0mm}{6.0mm}
 \!\!\!({\it IV\!.34})\ $\left\{\!\!\ba{l} \dot x= l,\ l\in\{-1,1\} \\[-0.8mm]
                \dot y=y(y-x) \ea\!\!\right.$ & \multicolumn{2}{c|}{$ \ba{c}\eta>0,\
                \theta= B_2=\mu=0, \\B_3\ne0,\,  H_7=H_9=H_{10}=0\ea $} & Config.\ 4.34\\
\hline
\end{tabular}
\end{center}
\end{table}
\begin{table}[!htb]
\begin{center}
\begin{tabular}{|l|c|c|}
\multicolumn{3}{r}{\bf Table 2}{\it(continued)} \\[1mm]
\hline \raisebox{-0.7em}[0pt][0pt]{\qquad Orbit representative}  &
Necessary and sufficient
           & \raisebox{-0.7em}[0pt][0pt]{Configuration} \\[-1mm]
                             & conditions & \\
 \hline\hline
 \!\!\!({\it IV\!.35})\ $\left\{\!\!\ba{l} \dot x= l+ y,\ l\in\!\R,\ l\ne0 \\[-0.8mm]
                \dot y= -xy\ea\!\!\right.$ & $
                \ba{c}\eta=0,\ M\ne0,\
                \theta=B_3=0,\\  N\ne0,\,\mu\!=\!0,\, H_7\ne0 \ea $ & Config.\ 4.35\\
\hline\rule[0mm]{0mm}{6.0mm}
 \!\!\!({\it IV\!.36})\ $\left\{\!\!\ba{l} \dot x= l,\quad l\in\{-1,1\},\\[-0.8mm]
                \dot y= -xy\!\ea\!\!\right.$ & $
                \ba{c}\eta=0,\ M\ne0,\
                \theta=B_3=K=0,\\  NH_6\ne0,\,\mu\!=\!H_7\!=\!0,\, H_{11}\!=\!0 \ea $ & Config.\ 4.36\\
\hline\rule[0mm]{0mm}{6.0mm}
 \!\!\!({\it IV\!.37})\ $\left\{\!\!\ba{l} \dot x= l+x, \ b\!\in\!\R,\ l\!\in\!\{0,1\},\!\!  \\[-0.8mm]
                \dot y= by\!-\!x^2,\ \ b(b^2\!-\!1)\!\ne\!0 \!\!\ea\!\!\right.$ & $
                \ba{c}\eta=M=0,\
                \theta=B_3=N=0,\\[-0.8mm] N_3D_1\ne0,\, N_6\ne0,\, D\ne0\ea $ & Config.\ 4.37\\
\hline\rule[0mm]{0mm}{6.0mm}
 \!\!\!({\it IV\!.38})\ $\left\{\!\!\ba{l} \dot x= l+x, \ b\!\in\!\R,\ l\!\in\!\{0,1\},\!\!  \\[-0.8mm]
                \dot y= b\!-\!x^2,\ \ b-l^2\ne\!0 \!\!\ea\!\!\right.$ & $
                \ba{c}\eta=M=0,\
                \theta=B_3=N=0,\\[-0.8mm] N_3D_1\ne0,\, N_6\ne0,\, D=0\ea $ & Config.\ 4.38\\
\hline\rule[0mm]{0mm}{6.0mm}
 \!\!\!({\it IV\!.39})\ $\left\{\!\!\ba{l} \dot x= x^2,    \\[-0.8mm]
                \dot y= x+y  \!\ea\!\!\right.$ & $
                \ba{c}\eta=0,\, M\ne0,\,
                \theta\!=\!\mu\!=\!N\!=\!B_3\!=\!0,\\[-0.8mm]
                N_1N_2\ne0,\,  K\!=\!0,\, N_5\!=\!0\ea $ & Config.\ 4.39\\
\hline\rule[0mm]{0mm}{6.0mm}
 \!\!\!({\it IV\!.40})\ $\left\{\!\!\ba{l} \dot x= 1+x,   \\[-0.8mm]
                \dot y= 1-xy\!\ea\!\!\right.$ & $
                \ba{c}\eta=0,\ M\ne0,\
                \theta=H_6=0,\\  NB_3\ne0,\,\mu\!=0,\,K=0 \ea $ & Config.\ 4.40\\
\hline\rule[0mm]{0mm}{6.0mm}
 \!\!\!({\it IV\!.41})\ $\left\{\!\!\ba{l} \dot x= l xy, \quad l\in\{-1,1\},  \\[-0.8mm]
                \dot y= y-x^2+lxy, \!\ea\!\!\right.$ & $
                \ba{c}\eta=M=0,\
                \theta\ne0,\\[-0.8mm] B_3=0,\ H_7=0,\
               D\ne0\ea $ & Config.\ 4.41\\
\hline\rule[0mm]{0mm}{6.0mm}
 \!\!\!({\it IV\!.42})\ $\left\{\!\!\ba{l} \dot x= l xy, \quad l\in\{-1,1\},  \\[-0.8mm]
                \dot y= -x^2+lxy, \!\ea\!\!\right.$ & $
                \ba{c}\eta=M=0,\
                \theta\ne0,\\[-0.8mm] B_3=0,\ H_7=D=0\ea $ & Config.\ 4.42\\
\hline\rule[0mm]{0mm}{6.0mm}
 \!\!\!({\it IV\!.43})\ $\left\{\!\!\ba{l} \dot x= lx^2,\ l\!\in\!\R,\, l(l^2\!-\!1)\!\ne\!0 \\[-0.8mm]
                \dot y= 1+(l-1)xy,\!\ea\!\!\right.$ & $
                \ba{c}\eta=0,\ M\ne0,\
                \theta=H_6=0,\\  NB_3\ne0,\,\mu\!=0,\,K\!\ne\!0,\, H_{11}\!=\!0 \ea $ & Config.\ 4.43\\
\hline\rule[0mm]{0mm}{6.0mm}
 \!\!\!({\it IV\!.44})\ $\left\{\!\!\ba{l} \dot x= x^2, \quad l\in\{-1,1\},  \\[-0.8mm]
                \dot y= l\!-\!x^2\!+\!xy \!\ea\!\!\right.$ & $
                \ba{c}\eta=M=0,\
                \theta=B_3=0,\\[-0.8mm] N_6N\ne0,\ H_{11}=0\ea $ & Config.\ 4.44\\
\hline\rule[0mm]{0mm}{6.0mm}
 \!\!\!({\it IV\!.45})\ $\left\{\!\!\ba{l} \dot x= l xy, \quad l\in\{-1,1\},  \\[-0.8mm]
                \dot y= x-x^2+lxy, \!\ea\!\!\right.$ & $
                \ba{c}\eta=M=0,\
                \theta\ne0,\\[-0.8mm] B_3=0,\
               H_7\ne0\ea $ & Config.\ 4.45\\
\hline\rule[0mm]{0mm}{6.0mm}
 \!\!\!({\it IV\!.46})\ $\left\{\!\!\ba{l} \dot x= 1,\\[-0.8mm]
                \dot y=y-x^2 \!\!\ea\!\!\right.$ & $
                \ba{c}\eta=M=0,\
                \theta=B_3=N=0,\\[-0.8mm] N_3D_1\ne0,\, N_6=0\ea $ & Config.\ 4.46\\
\hline
\end{tabular}
\end{center}

\vspace{5mm}
 \begin{center}
\begin{tabular}{|l|c|}
\multicolumn{2}{r}{\bf Table 3}\\[2mm]
\hline
 \hfil Perturbed systems\hfil & \hfil Invariant straight lines \hfil \\
\hline
\hline({\it IV.16${}_\varepsilon$})\,:$\ba{l} \dot x=
(l+x)(\varepsilon x+1),\ \
  \dot y= y(y-x)\ea$ &
  $ \ba{c}y=0,\, x=- l,\, \varepsilon x=-1 \ea $ \!\!\\[0mm]
\hline ({\it IV.17${}_\varepsilon$})\,:$\ba{l} \dot x=
x(\varepsilon x+1),\ \
  \dot y= y(y-x)\ea$ &
  $ \ba{c}y=0,\, x=0,\, \varepsilon x=-1 \ea $ \!\!\\[0mm]
\hline ({\it IV.18${}_\varepsilon$})\,:$\left\{\!\!\ba{l} \dot x=
(l^2+l+lx+y)(\varepsilon x+1),\\[-0.8mm]
  \dot y= y[(y-x)-\varepsilon y\big[l(l+1)(\varepsilon+1)+1-y\big]\!\! \ea\right.$ &
  $ \ba{c}y=0,\,  \varepsilon x+1=0,\\[-0.8mm]
      y-x(1\!-\!l\varepsilon)\!=\!(l\!+\!1)(l\varepsilon+1)\!\!  \ea $ \!\!\\[0mm]
\hline ({\it IV.19${}_\varepsilon$})\,:$\ba{l} \dot x=
            (l+x)(\varepsilon x+1),\quad
           \dot y= -xy\!\! \ea$ & $ \ba{c}y=0,\ x=-l,\  \varepsilon x=-1  \ea $ \\[0mm]
\hline ({\it IV.20${}_\varepsilon$})\,:$\ba{l} \dot x= l
x^2+(\varepsilon+1)xy,\
  \dot y= (l-1)xy+y^2\!\! \ea$ & $ \ba{c}x=0,\, y=0,\, x+ \varepsilon y=0  \ea $ \\[0mm]
\hline ({\it IV.21${}_\varepsilon$})\,:$\ba{l} \dot x=\varepsilon
x\!+\! l x^2\!+\!(\varepsilon\!+\!1)xy,\
  \dot y= (y\!+\!1)(lx\!-\!x\!+\!y)\!\! \ea$ & $ \ba{c}x=0,\, y=-1,\, x+ \varepsilon y=-\varepsilon  \ea $ \\[0mm]
\hline({\it IV.22${}_\varepsilon$})\,:$\ba{l} \dot x=
l(x^2-\varepsilon^2),\
  \dot y= (y+1)[y+(l-1)x-1]\!\! \ea$ & $ \ba{c}y+1=0,\, x=\pm  \varepsilon  \ea $ \\[0mm]
\hline ({\it IV.23${}_\varepsilon$})\,:$\ba{l} \dot x= x^2+
xy,\quad   \dot y= (y+l)^2-\varepsilon^2\!\! \ea$ & $ \ba{c}x=0,\quad y+l=\pm \varepsilon  \ea $ \\[0mm]
\hline ({\it IV.24${}_\varepsilon$})\,:$\ba{l} \dot x=
l(x+1)^2-l\varepsilon^2,\quad
           \dot y= (l-1)xy\!\! \ea$ & $ \ba{c}y=0,\quad x+1=\pm \varepsilon  \ea $ \\[0mm]
\hline ({\it IV.25${}_\varepsilon$})\,:$\left\{\!\!\ba{l} \dot
x=\varepsilon
lx+ l x^2+(\varepsilon+1)xy,\\[-0.8mm]
  \dot y=y+ (l-1)xy+y^2\!\! \ea\right.$ & $ \ba{c}x=0,\, y=0,\, x+ \varepsilon y=-\varepsilon  \ea $ \\[0mm]
\hline ({\it IV.26${}_\varepsilon$})\,:$\ba{l} \dot x=\varepsilon
x+(\varepsilon\!+\!1)xy,\quad
  \dot y= (y\!+\!1)(y-x)\!\! \ea$ & $ \ba{c}x=0,\, y=-1,\, x+ \varepsilon y=-\varepsilon  \ea $ \\[0mm]
\hline
\end{tabular}
\end{center}
\end{table}

\begin{table}[!htb]
\begin{center}
\begin{tabular}{|l|c|}
\multicolumn{2}{r}{\bf Table 3}{\it(continued)}\\[2mm]
\hline
 \hfil Perturbed systems\hfil & \hfil Invariant straight lines \hfil \\
\hline
\hline
({\it IV.27${}_\varepsilon$})\,:$\left\{\!\!\ba{l}
\dot x=
2(1-2\varepsilon)(lx+y)(1+\varepsilon\,x),\\[-0.8mm]
  \dot y= l^2+1 +2(l^2+1)\varepsilon\,x \\[-0.8mm]
  \ +(1-2\varepsilon)\left[-x^2
  +2l\varepsilon\,xy -(1-2\varepsilon)y^2\right]\!\! \ea\right.$ &
  $ \ba{c}\varepsilon x+1=0,\\ (1-2 \varepsilon) (x\pm i y)=1\mp il\!\!  \ea $ \!\!\\[0mm]
\hline ({\it IV.28${}_\varepsilon$})\,:$\ba{l} \dot x=
            x^2-1,\quad
           \dot y= (x+ly)(1+\varepsilon y)\!\! \ea$ & $\ba{c}x=\pm1,\quad \varepsilon y=-1 \ea $ \\[0mm]
\hline ({\it IV.29${}_\varepsilon$})\,:$\ba{l} \dot x=
            x^2-1,\quad
           \dot y= (l+x)(1+\varepsilon y)\!\! \ea$ & $\ba{c}x=\pm1,\quad \varepsilon y=-1 \ea $ \\[0mm]
\hline ({\it IV.30${}_\varepsilon$})\,:$\left\{\!\!\ba{l} \dot x=
            (1+x)(1+lx)-\varepsilon,\\[-0.8mm]
           \dot y= 1+(l-1)xy-\varepsilon y^2\!\! \ea\right.$ & $ \ba{c}x+ \varepsilon y+1=0\\[-0.8mm]
                                                  lx^2+(l+1)x+1=\varepsilon  \ea $ \\[0mm]
\hline
  ({\it IV.31${}_\varepsilon$})\,:$\left\{\!\!\ba{l} \dot x=
           -l\varepsilon +(1+\varepsilon)x+ (1+\varepsilon)x^2,\\[-0.8mm]
           \dot y= l-x^2+ xy\!\! \ea\right.$
            & $ \!\!\!\ba{c} x+ \varepsilon y=-1-\varepsilon,\\[-0.8mm]
             (1+\varepsilon)x^2+(1+\varepsilon)x=l\varepsilon \!\!\! \ea \!$ \\[0mm]
\hline ({\it IV.32${}_\varepsilon$})\,:$\ba{l} \dot x=
            x^2+1,\quad
           \dot y= (x+ly)(1+\varepsilon y)\!\! \ea$ & $\ba{c}x=\pm i,\quad \varepsilon y=-1 \ea $ \\[0mm]
\hline ({\it IV.33${}_\varepsilon$})\,:$\ba{l} \dot x=
            x^2+1,\quad
           \dot y= (l+x)(1+\varepsilon y)\!\! \ea$ & $\ba{c}x=\pm i,\quad \varepsilon y=-1 \ea $ \\[0mm]
\hline ({\it IV.34${}_\varepsilon$})\,:$\ba{l} \dot x=
l(1-\varepsilon^2 x^2),\ \
  \dot y= y(y-x)\ea$ &
  $ \ba{c}y=0,\ \  \varepsilon x=\pm1 \ea $ \!\!\\[0mm]
\hline ({\it IV.35${}_\varepsilon$})\,:$\left\{\!\!\ba{l} \dot
x=(\varepsilon x+1)[(\varepsilon+1)(y+l)+l\varepsilon x],\\[-0.8mm]
  \dot y=l\varepsilon^2y+ (l\varepsilon^3-1)xy+\varepsilon^2y^2\!\! \ea\right.$ &
                           $ \ba{c} y=0,\, \varepsilon x=-1,\\[-0.8mm]
                             \varepsilon^2(x+\varepsilon y)+\varepsilon=-1    \ea $ \\[0mm]
\hline ({\it IV.36${}_\varepsilon$})\,:$\ba{l} \dot x=
            l(1-\varepsilon^2 x^2),\quad
           \dot y= -xy\!\! \ea$ & $ \ba{c}y=0,\quad   \varepsilon x=\pm1  \ea $ \\[0mm]
\hline
  ({\it IV.37${}_\varepsilon$})\,:$\left\{\!\!\ba{l} \dot x=l+x+
            \varepsilon (2-l\varepsilon)x^2,\\[-0.8mm]
           \dot y= by-x^2+\varepsilon(1+b-2l\varepsilon) xy+\varepsilon^2(b-l\varepsilon)y^2\!\! \ea\right.$
            & $ \!\!\!\ba{c} \varepsilon x+ \varepsilon^2 y+1=0,\\[-0.8mm]
             \varepsilon(2-l\varepsilon)x^2+x+l=0\!\!\! \ea \!$ \\[0mm]
\hline
  ({\it IV.38${}_\varepsilon$})\,:$\left\{\!\!\ba{l} \dot x=l+x+
            \varepsilon (2-l\varepsilon+b\varepsilon^2)x^2,\\[-0.8mm]
           \dot y= b-x^2+\varepsilon(1-2l\varepsilon-2b\varepsilon^2) xy-\varepsilon^3( l+b\varepsilon)y^2\!\! \ea\right.$
            & $ \!\!\!\ba{c} \varepsilon x+ \varepsilon^2 y+1=0,\\[-0.8mm]
             \varepsilon(2-l\varepsilon+b\varepsilon^2)x^2+x=-l\!\!\! \ea \!$ \\[0mm]
\hline ({\it IV.39${}_\varepsilon$})\,:$\ba{l} \dot x=
            x^2-\varepsilon^2,\quad
           \dot y= (x+ly)(1+\varepsilon y)\!\! \ea$ & $\ba{c}x=\pm \varepsilon,\quad \varepsilon y=-1 \ea $ \\[0mm]
\hline ({\it IV.40${}_\varepsilon$})\,:$\left\{\!\!\ba{l} \dot x=
            (1+x)(1+\varepsilon x)-\varepsilon,\\[-0.8mm]
           \dot y= 1+(\varepsilon-1)xy-\varepsilon y^2\!\! \ea\right.$ & $ \ba{c}x+ \varepsilon y+1=0\\[-0.8mm]
                                                  \varepsilon x^2+(\varepsilon+1)x+1=\varepsilon  \ea $ \\[0mm]
\hline
  ({\it IV.41${}_\varepsilon$})\,:$\left\{\!\!\ba{l} \dot x=
           l\varepsilon x^2/2+lxy,\\[-0.8mm]
           \dot y= \varepsilon^2\!+\!2\varepsilon x\!+\!(1\!+\!2l\varepsilon^2)y\!-\!x^2\!+\!
           2l\varepsilon xy\!+\!l(1\!+\!l\varepsilon^2)y^2\!\! \ea\right.$
            & $ \ba{c}x=0, \ 2x+ l\varepsilon y=-\varepsilon,\\[-0.8mm]
              x-2l\varepsilon y=2\varepsilon \ea $ \\[0mm]
\hline
  ({\it IV.42${}_\varepsilon$})\,:$\left\{\!\!\ba{l} \dot x=
           \varepsilon x^2/2+lxy,\\[-0.8mm]
           \dot y= -x^2+
           2\varepsilon xy+(l+\varepsilon^2)y^2\!\! \ea\right.$
            & $ \ba{c}x=0, \ 2x+ \varepsilon y=0,\\[-0.8mm]
              x-2\varepsilon y=0 \ea $ \\[0mm]
\hline ({\it IV.43${}_\varepsilon$})\,:$\ba{l} \dot x=
            l(x^2-\varepsilon^2),\quad
           \dot y= 1+(l-1)xy-l\varepsilon^2 y^2\!\! \ea$ & $ \ba{c}x=\pm\varepsilon,\
                                                  x+l\varepsilon^2y=0  \ea $ \\[0mm]
\hline
  ({\it IV.44${}_\varepsilon$})\,:$\left\{\!\!\ba{l} \dot x=
           -l\varepsilon +\varepsilon(1+\varepsilon)x+ (1+\varepsilon)x^2,\\[-0.8mm]
           \dot y= l-x^2+ xy\!\! \ea\right.$
            & $ \!\!\!\ba{c} x+ \varepsilon y=-\varepsilon(1+\varepsilon),\\[-0.8mm]
             (1+\varepsilon)x^2+\varepsilon(1+\varepsilon)x=l\varepsilon \!\!\! \ea \!$ \\[0mm]
\hline
 ({\it IV.45${}_\varepsilon$})\,:$\left\{\!\!\ba{l} \dot x=
            \varepsilon(\varepsilon^2-l)x/l+\varepsilon x^2+lxy,\\[-0.8mm]
           \dot y= x\!+\!\varepsilon(2\varepsilon^2\!-\!l)y/l\!-\!x^2\!-
           \!2\varepsilon xy\!+\!(l\!-\!2\varepsilon^2)y^2\!\! \ea\right.$
            & $ \ba{c}x=0, \ x+ \varepsilon y=0,\\[-0.8mm]
              x+2\varepsilon y=2\varepsilon^2/l  \ea $ \\[0mm]
\hline
  ({\it IV.46${}_\varepsilon$})\,:$\left\{\!\!\ba{l} \dot x=1+\varepsilon x+
            \varepsilon x^2,\\[-0.8mm]
           \dot y= y-x^2+\varepsilon(1-\varepsilon) xy+\varepsilon^2( l-\varepsilon)y^2\!\! \ea\right.$
            & $ \ba{c} \varepsilon x+ \varepsilon^2 y+1=0,\\[-0.8mm]
             \varepsilon x^2+\varepsilon x+1=0  \ea \!$ \\[0mm]
\hline
\end{tabular}
\end{center}
\end{table}

\noindent {\it Proof of the Main Theorem:}  Since  we only discuss
the case $C_2\ne0$, in what follows it suffices to consider only
the canonical forms $(\SSS_I)$ to $(\SSS_{I\,V})$.  The idea of
the proof is to perform a case by case discussion for each one of
these canonical forms, for which   according   to
Lemma~\ref{lm_NB2}  we must examine two subcases: $(i)\
\theta\ne0,$\ $B_3=0$ and  $ (ii)\ \theta=B_2=0$. Each one of
these conditions yields specific conditions on the parameters. The
discussion proceeds further by breaking these cases in more
subcases and then by constructing new invariants or T-comitants to
fit the conditions on parameters.


\subsection{Systems with the divisor\ $D_S(C,Z)=1\cdot \omega_1+1\cdot \omega_2+1\cdot \omega_3$}

 For this case we shall later  need the following $T$-comitants.
\bnot\label{not_H4,5} Let us denote
$$
\bal
&H_7(a)=(N,C_1)^{(2)},\quad  H_9(\ab)=
-\Big(\big((D,D)^{(2)},D,\big)^{(1)}D\Big)^{(3)},\\
 & H_8(\ab)= 9\Big((C_2,D)^{(2)},(D,D_2)^{(1)}\Big)^{(2)}+
2\Big[(C_2,D)^{(3)}\Big]^2,\\
& H_{10}(\ab)= \big((N,D)^{(2)},\ D_2\big)^{(1)}.
\eal
$$
\enot
 According to Lemma \ref{lm_3:2} the systems with this type
of  divisor can be brought by linear transformations  to the
canonical form $(\SSS_I)$ for which we have:
\be\label{theta:I}
\theta(\ab,x,y)= -8(g-1)(h-1)(g+h).
\ee
\vspace{-5mm}
\subsubsection{The case $\theta\ne0,$  $B_3=0$}

The condition $\theta\ne0$ yields $(g-1)(h-1)\ne0$ and in
$(\SSS_I)$  we may assume $d=e=0$ via the translation $x\to
x+d/(1-h)$ and $y\to y+e/(1-g)$. Thus  we obtain the  systems
\be\label{S1:1}
  \dot x=k + cx + gx^2+(h-1)xy,\quad
  \dot y=l + fy + (g-1)xy+ hy^2,
\ee
for which we  calculate
$$
\bal
B_3=7& 3l(g-1)^2x^3(2y-x)+ 3k(h-1)^2y^3(y-2x) +\\
&\ 3[(c-f)(fg+ch) - k(1 + g)(-1 + g + 2 h) + l(1 + h)(-1 + 2 g + h)]x^2y^2.\\
\eal
$$
Then by $\theta\ne0$ the condition $B_3=0$ yields $k=l=0$ and
$(c-f)(fg+ch)=0$. Hence $c-f=0$ or $fg+ch=0$.  An invariant which
capture the condition $c-f=0$ is $H_7$. Indeed, for the systems
(\ref{S1:1}) we have $H_7=4(f-c)(g-1)(h-1)$ and since $\theta\ne0$
 the condition $f-c=0$ is equivalent with $H_7=0$.

\textbf{Subcase $H_7=0$.} Then $f=c$  and we obtain the systems
\be\label{CS1:1}
  \dot x= cx + gx^2+(h-1)xy,\qquad
  \dot y= cy + (g-1)xy+ hy^2,
\ee
for which calculations yield  $ {\cal H}=\gcd\left({\cal
E}_1,{\cal E}_2\right)= 2XY(X-Y)$ in the ring $\R[c,g,h,X,Y,Z]$
which means that for a concrete system corresponding to
 ($\mathbf{c ,  g, h}$),\   $M_{{}_\IL}\ge4$ and by
Lemma \ref{gcd:3} $M_{{}_\IL}$ cannot be 5.

Let us examine the singularities of systems \eqref{CS1:1}. Clearly
for $c=0$ the point $(0,0)$ will be of the multiplicity four and
since for these systems we have $H_1(\ab)=576 c^2$ the condition
$c=0$ is equivalent to $H_1=0$. Hence, for $H_1=0$ we obtain
Config. 4.5.

Consider now $H_1\ne0$ (i.e. $c\ne0$). By Remark \ref{rem:transf}
( $\gamma=c,\ s=1$) we may assume $c=1$ and then for
$gh(g+h-1)\ne0$ the systems
\eqref{CS1:1} possess the following finite singular points:
$$
(0,\,0); \quad (0,\ -1/h); \quad (-1/g,\ 0); \quad
\big(1/(1-g-h),\ 1/(1-g-h)\big).
$$
On the other hand for systems \eqref{CS1:1} we have
$\mu=32gh(g+h-1)$ and hence for $\mu\ne0$ we get the Config. 4.3.

Consider now $\mu=0$ for which we obtain $gh(g+h-1)=0$ and without
loss of generality we may assume $h=0$. Indeed, if $g=0$
(respectively, $g+h-1=0$) we can apply the  linear transformation
which will replace the straight line $x=0$ with $ y=0$
(respectively, $x=0 $ with $y=x$) which reduces the case $g=0$
(respectively, $g+h-1=0$) to the case $h=0$. In this case from
\eqref{theta:I} we have $\theta=8g(g-1)\ne0$ and hence we obtain Config.
4.4.

\smallskip
\textbf{Subcase $H_7\ne0$.} Then $c-f\ne0$  and hence, the
condition $B_3=0$ yields  $fg+ch=0$. If    $g=0$ then   $ch=0$. In
this case $\theta=8h^2\ne0$ hence we must have   $c=0$   yielding
degenerate  systems. So, $g\ne0$ and by introducing a new
parameter $u$, $c=gu$ we obtain $f=-hu$. Then the condition
$H_7=-4(g-1)(h-1)(g+h)u\ne0$ implies $u\ne0$ and we may assume
$u=1$ via Remark \ref{rem:transf} ($\gamma=u$, $s=1$). This leads
to the systems:
\be\label{CS1:2}
  \dot x= g x + gx^2+(h-1)xy,\qquad
  \dot y= -h y + (g-1)xy+ hy^2,
\ee
for which calculations yield  \ $ {\cal H}=\gcd\left({\cal
E}_1,{\cal E}_2\right)= 2XY(X-Y+Z)$  in the ring $\R[g,h,X,Y,Z]$
which means that for a concrete system corresponding to
$\mathbf{g}$, $\mathbf{h}$,   $M_{{}_\IL}\ge4$ and by Lemma
\ref{gcd:3} $M_{{}_\IL}$ cannot be 5. The singularities of systems
\eqref{CS1:2} are
$$
(0,0);\quad (0,1);\quad (-1,0);\quad (-h,g).
$$
The point $(-h,g)\not\in  {\cal G}$ (otherwise the systems
\eqref{CS1:2} become degenerate). So we obtain the Config 4.4.

\vspace{-5mm}
\subsubsection{The case $\theta=B_2=0$}

According to \eqref{theta:I} the condition $\theta=0$ yields
$(g-1)(h-1)(g+h)=0$ and without loss of generality we can consider
$h=1$. Indeed, if $g=1$ (respectively, $g+h=0$) we can apply the
linear transformation which will replace the straight line $x=0$
with $ y=0$ (respectively, $x=0 $ with $y=x$) reducing this case
to $h=1$.  Assuming $h=1$ for the systems~$(\SSS_I)$ we calculate
$N=(g^2-1)x^2$ (see Lemma \ref{lm4}) and we shall examine two
subcases: $N\ne0$ and $N=0$.

\fbox{\textbf{Subcase $N\ne0$.}} Then $(g-1)(g+1)\ne0$. For
systems $(\SSS_I)$ with $h=1$ we have $\mu=32g^2$ and we shall
consider two subcases: $\mu\ne0$ and $\mu=0$.

\textbf{1)} If $\mu\ne0$ then  $g\ne0$ and  we may assume $c=f=0$
via the translation $x\to x-c/(2g)$,\ \ $y\to
y+\big[c(g-1)-2fg\big]/(4g)$.
 Thus we obtain the systems
\be\label{S1:3}
  \dot x=k + dy + gx^2,\quad  \dot y=l + ex+ (g-1)xy +y^2,
\ee
for which $ H_7=4d(g^2-1)$. We claim that a given system $S(\ab)$
of the form
\eqref{S1:3} to belong to the class $\QSL_{\bf4}$ $\mathbf{ d}=0$ (i.e. $H_7(\ab)=0$)
is necessary. Indeed,  for systems
\eqref{S1:3} we calculate
$$
C_2=xy(x-y),\quad H=-x\big[(g-1)^2x+4gy]
$$
and hence there exist 3 directions $(u,v)$  for invariant lines
$ux+vy+w=0$, and namely: $(1,0),$ $(0,1)$ and $(1,-1)$. Moreover,
according to Lemma \ref{lem:H} the parallel lines could be only in
the directions given by the $T$-comitant $H$, i.e. $(1,0)$ or
$((g-1)^2,4g)$. However since $N\ne0$ (i.e. $g^2-1\ne0$) we obtain
$g-1\ne0$ and $(g-1)^2\ne-4g$ (otherwise
$0=(g-1)^2+4g=(g+1)^2\ne0$). Thus  the parallel lines could only
be in the direction $(1,0)$ and  we conclude that to have
$M_{{}_\IL}=4$ for a system  $S(\ab)$,  the existence of at least
one line in the direction $(1,0)$ is necessary, i.e. of a line
with equation $x+\alpha=0$. Calculations yield $d=0$ in this case
and our claim is proved.

Thus we assume $H_7=0$ (i.e. $d=0$). Then for systems
\eqref{S1:3} we calculate:
$$
B_2=-648\big[e^2+l(g-1)^2\big]\big[e^2+(l-k)(g+1)^2\big]
$$
and the condition $B_2=0$ yields either
$$
(i)\quad e^2+l(g-1)^2 =0\qquad \text{or}\qquad (ii)\quad
e^2+(l-k)(g+1)^2=0.
$$
We claim that the case $(i)$ can be reduced by a linear
transformation and time rescaling to the case $(ii)$ and
viceversa. Indeed, via the transformation $x_1=x$, $y_1=x-y$ and
$t_1=-t$  systems \eqref{S1:3} with $d=0$ keep  the same form
$$
  \dot x_1=\tilde k + \tilde gx_1^2,\quad
  \dot y=\tilde l + \tilde e x_1+ (\tilde g-1)x_1y_1 +y_1^2
$$
but with new parameters: $ \tilde k=-k,\ \tilde g=-g,\ \tilde
l=l-k,$\ and\ $\tilde e=e$.\ Then obviously we have:
$$
\bal
 & \tilde e^2+\tilde l(\tilde g-1)^2 = e^2+(l-k)(g+1)^2,\qquad
   \tilde e^2+(\tilde l-\tilde k)(\tilde g+1)^2 = e^2+l(g-1)^2
\eal
$$
and this proves our claim.

In what follows we assume that the condition $(i)$ holds. Since
$g-1\ne0$ we may set $e=u(g-1)$ (where $u$ is a new parameter) and
then we obtain $l=-u^2$. So we get the systems
\be\label{S1:4}
  \dot x=k + gx^2,\quad  \dot y=(y+u)\big[y+(g-1)x-u\big]
\ee
which posses three invariant straight lines: $y+u=0$ and
$gx^2+k=0$. So for a concrete system in this family
$M_{{}_\IL}\ge4$. We shall find the conditions on the parameters,
for systems \eqref{S1:4} to be  in the class $\QSL_{\bf4}$, i.e.
to possess exactly 4 invariant lines, including the line at
infinity and including multiplicities. Calculations yield
$$
\bal
&{\cal E}_1=2\big[(g-1) X^2 - (g-3)XY - 2Y^2 -u(g+1)XZ +2u Y Z - k(g-1) Z^2\big]{\cal H},\\
&{\cal E}_2= \big[(g-1) X + Y-u Z\big]\big[g X^2 +(1-g) X Y  - Y^2 + u(1-g) X Z + (k+u^2) Z^2\big]{\cal H},\\
\eal
$$
where $ {\cal H}=  (Y+uZ)(g X^2 + k Z^2)\in\R[g,k,u, X,Y,Z]$. The
condition on the parameters $g,k,u$ so as to have an additional
common factor of ${\cal E}_1$, ${\cal E}_2$ according to Lemma
\ref{Trudi:2} is
$$
\Res_X({\cal E}_1/{\cal F},\ {\cal E}_2/{\cal F})=
8(1-g)\big[4gu^2 +k(g + 1)^2 \big] \big[g Y^2 -2gu  Y Z  +
(k(g-1)^2+gu^2) Z^2\big]^2Z^2 \equiv0
$$
in $\R[X,Y,Z]$.
 Thus, since $g(g-1)\ne0$  the condition
$\Res_X \equiv0$ is equivalent to the condition $4gu^2 +k(g + 1)^2
=0$. So the condition for $(S)\in \QSL_{\bf4}$ is $4gu^2 +k(g +
1)^2\ne0$ which in view of
$$
B_3=-3[4gu^2 +k(g + 1)^2] x^2 y^2
$$
means $B_3\ne0$.

It remain to examine in more details the invariant lines
configuration of systems~\eqref{S1:4}. We observe that the  lines
depends on  $gk$. For these systems we have
$$
H_{10}=-32gk(g^2-1),\qquad N=(g^2-1)x^2.
$$
So $H_{10}N=-32gk(g^2-1)^2x^2$. Since $g^2-1\ne0$ for a real value
$x\ne0$, $H_{10}N=-gk \Lambda$ with $\Lambda>0$.

 \textbf{a)} Assume $H_{10}N>0$. Then $gk<0$ and since
$\mu\ne0$ (i.e. $g\ne0$) we may set $k=-gv^2$, where $v\ne0$ is a
new parameter. We may assume $v=1$ via Remark \ref{rem:transf}
($\gamma=v$,$s=1$) and  systems~\eqref{S1:4} become
\be\label{S1:5}
  \dot x= g(x^2-1),\quad  \dot y=(y+u)\big[y+(g-1)x-u\big]
\ee
with the singularities:
$$
(1,-u),\quad (-1,-u),\quad (1,u-g+1),\quad (-1,u+g-1).
$$
So among these there can be double points if and only if either
$g-1=2u$ or $g-1=-2u$. On the other hand for systems~(\ref{S1:5})
we have $H_9=2^{10}3^2g^6\big[(g-1)^2-4u^2\big]^2$ and hence, for
$H_9\ne0$ we get Config. 4.9 whereas for  $H_9=0$ we obtain Config
4.6. In the last case we may choose the relation $g-1=2u$ as
$g-1=-2u$ leads to this relation via the change $x\to-x$, $y\to
-y$ and $t\to -t$.

 \textbf{b)} Suppose now  $H_{10}N<0$, i.e. $gk>0$. The
 systems \eqref{S1:4} possess 2 imaginary affine parallel invariant
 lines.
 As above we may set $k=gv^2$  and assuming $v=1$ by
Remark \ref{rem:transf} ($\gamma=v$,$s=1$), we obtain the systems
\be\label{S1:5a}
  \dot x= g(x^2+1),\quad  \dot y=(y+u)\big[y+(g-1)x-u\big]
\ee
without real finite  singularities. Thus we get the Config. 4.13.

 \textbf{c)} For $H_{10}=0$, since $g\ne0$ we obtain $k=0$ and for
 systems \eqref{S1:4} the invariant line
 $x=0$ is a double line. Since in this case $B_3=-12gu^2x^2y^2\ne0$
 we can assume $u=1$ via the Remark \ref{rem:transf}
 ($\gamma=u$,$s=1$). Thus we obtain the systems
$$
  \dot x= gx^2,\quad  \dot y=(y+1)\big[y+(g-1)x-1\big]
$$
and this leads to the Config. 4.22.

\smallskip
\textbf{2)} Assume $\mu=0$.  Then  $g=0$ and  we may consider
$e=f=0$ via the translation $x\to x+2e+f$,\  $y\to y+e$. So the
systems~$(\SSS_I)$ with $\theta=\mu=0$ (i.e. $h=1$, $g=0$) become
\be\label{S1:6}
  \dot x=k + cx+ dy,\qquad  \dot y=l  -xy +y^2
\ee
for which \ $ B_2=-648 l(c^2 + cd - k + l)x^4=0$.\ Hence the
condition $B_2=0$ yields either $(i)\ l=0$ or $ (ii)\ c^2 + cd - k
+ l=0$.

We claim that the case $(ii)$ can be reduced by an affine
transformation and time rescaling to the case $(i)$ and viceversa.
Indeed, via the transformation
$$
x_1=x+2c+d,\qquad y_1=x-y+c+d\qquad t_1=-t
$$
the systems \eqref{S1:6} keep the same form
$$
  \dot x_1=\tilde k + \tilde cx_1+ \tilde dy_1,\qquad  \dot y_1=\tilde l  -x_1y_1 +y_1^2
$$
but with new parameters:
$$
\bal
 \tilde k=2c(c+d)-k,\quad \tilde c=-(c+d),\quad \tilde
l=l-k+c(c+d),\quad \tilde d=d.
\eal
$$
Then obviously we have:
$$
\bal
 & \tilde l=c^2 + cd - k + l,\qquad
   \tilde c^2 + \tilde c\tilde d - \tilde k + \tilde l = l
\eal
$$
and this proves our claim.

 Hence we only need to consider  $l=0$. For  systems
(\ref{S1:6}) we have   $H_7=-4d$ and we shall consider two
subcases: $H_7\ne0$ and $H_7=0$.

\textbf{a)} Consider first the case  $H_7\ne0$. Then $d\ne0$ and
we may assume $d=1$ via Remark \ref{rem:transf}
 ($\gamma=d$, $s=1$). For systems (\ref{S1:6}) with
$l=0$ and $d=1$ calculations yield:
$$
\bal
 &{\cal E}_1=\big[-2Y^3+Y^2F_1(X,Z)+YF_2(X,Z)+F_3(X,Z)\big]{\cal H},\\
&{\cal E}_2=(Y-X)(c X + Y + k Z)(X Y -Y^2 +c X Z +  Y Z+k
Z^2){\cal H},
\eal
$$
where ${\cal H}=YZ\in \R[c,k,X,Y,Z]$ and $F_i\in\R[c,k,X,Z]\
(i=1,2,3)$ are homogeneous in $X,Z$ of degree $i$. Hence
$\deg\,{\cal H}=2$ and for a system $S(\ab)$ to have an additional
common factor of ${\cal E}_1$ and ${\cal E}_2$   according to
Lemma \ref{Trudi:2}  it is necessary
 that  the following identity holds in $\R[X,Z]$:
$$
 \Res_Y({\cal E}_1/{\cal H},\ {\cal E}_1/{\cal H})=
  -16(c^2+c  -k)(c X + k Z)^2(c X +  X + k Z)^6 Z^4 =0.
$$
So,  we obtain the condition $(c^2 + c  -
k)(c^2+k^2)((c+1)^2+k^2)=0$. However, for $k=c(c+1)=0$ we get a
degenerate system. On the other hand for systems (\ref{S1:6}) with
$l=0$ and $d=1$ we calculate $B_3=3(c^2 + c  - k)x^2y^2$ and
hence, to have an additional common factor of ${\cal E}_1$ and
${\cal E}_2$  it is necessary that $B_3=0$. Then $k=c(c + 1)$ and
we get the systems
\be\label{S1:7}
  \dot x=c(c+1) + cx+y,\qquad  \dot y= -xy +y^2
\ee
for which we obtain\ ${\cal H}=\gcd\left({\cal E}_1,{\cal
E}_2\right)=  Y Z ( X - Y + c Z +  Z) $, i.e.  $\deg {\cal H}=3$.

We claim that these systems  belong to class $\QSL_{\bf4}$ or they
are  degenerate. Indeed, for systems \eqref{S1:7} calculations
yield:
$$
\bal
&{\cal E}_1=\big(-c X^2 + 2 c X Y +  Y^2 + 2 c^2 Y Z + 2 c  Y Z
+     c^3 Z^2 + c^2  Z^2){\cal H},\\
&{\cal E}_2= -3(X - Y)(Y + c Z)(c X +  Y + c^2 Z + c  Z){\cal H}.
\eal
$$
Hence by  Lemma \ref{Trudi:2} to have an additional common factor
of ${\cal E}_1$ and ${\cal E}_2$    the condition
$$
\Res_Y({\cal E}_1/{\cal H},\ {\cal E}_2/{\cal H})= 9 c^2  (c +
1)^2 (X + c Z)^6 =0
$$
must hold in $\R[X,Z]$. Therefore   $c=0$ or $c=-1$. Both cases
lead to  degenerate systems~(\ref{S1:7}). Our claim is proved.

We observe that  systems~(\ref{S1:7})  possess the two invariant
affine straight lines:\  $ y=0$ and $ y=x+c+1.$ Taking into
account that  $Z\mid{\cal H}$, we have  by Corollary
\ref{Mult:Z=0} that the line $Z=0$ could be of multiplicity 2.
This is confirmed by the perturbations ({\it
IV.18${}_\varepsilon$}) from Table 3. On the other hand the
systems~(\ref{S1:7}) have two finite singularities: $(-c-1,0)$ and
$(-c,-c)$ with $c\ne0$ (otherwise we get a degenerate system).
 Thus, we get Config. 4.18.

\textbf{b)} Assume now  $H_7=0$. Then $d=0$ and we get the systems
\be\label{S1:8}
  \dot x=k + cx,\qquad  \dot y= -xy +y^2
\ee
for which we obtain
$$
\bal
&{\cal E}_1=2(-X^2 + 3 X Y - 2 Y^2 + c Y Z + k Z^2){\cal H},\quad
 {\cal E}_2= (Y-X)(X Y - Y^2 + c X Z + k Z^2){\cal H}.
\eal
$$
where  ${\cal H}=  YZ(cX+kZ) =\gcd\left({\cal E}_1,{\cal
E}_2\right)$\ and due to Lemma \ref{Trudi:2} in order to have for
a specific system an additional nontrivial common factor of ${\cal
E}_1$ and ${\cal E}_2$  it is necessary that the following
condition holds in $\R[Y,Z]$:
$$
\Res_X({\cal E}_1/{\cal H},\ {\cal E}_2/{\cal H})= 4(c^2 -k ) ( c
Y+ k Z)^2Z^4 =0,
$$
i.e. $c^2 -k=0$ or $c=k=0$. The second condition leads to
degenerate system.  Since for systems~\eqref{S1:8} we have
$B_3=3(c^2-k)x^2y^2$ evidently the condition $\Res_X\not\equiv0$
is equivalent to $B_3\ne0$. We observe, that these  systems
possess 2 invariant affine lines $y=0$ and $cx+k=0$ for $c\ne0$.
Taking into account the value of ${\cal H}$, we have by Corollary
\ref{Mult:Z=0} that the line $Z=0$ could be of multiplicity two
and even three if $c=0$. This is confirmed by the perturbations
({\it IV.16${}_\varepsilon$}) and ({\it IV.34${}_\varepsilon$})
from Table 3. It remains to note that in the case $c\ne0$ the
systems~\eqref{S1:8} possess the finite singular points $(-k/c,0)$
and $(-k/c,-k/c)$. Moreover, these points are distinct for $k\ne0$
and they coincide for $k=0$. These points lie on the invariant
lines $y=0$ and $cx+k=0$.

Thus, since for the systems~(\ref{S1:8}) $H_{10}=-8c^2$ and
$H_{9}=-576c^4k^2$, we obtain:  Config. 4.16 for $H_9\ne0$,
 Config. 4.17   for $H_9=0$, $H_{10}\ne0$   and
Config. 4.34 for $H_9=H_{10}=0$. We note that for $c\ne0$
(respectively, $c=0,$ $k\ne0$)  we may assume $c=1$ (respectively,
$k\in\{-1,1\}$)  via Remark \ref{rem:transf} ($\gamma=c$,$s=1$)
(respectively, $\gamma=|k|$,$s=1/2$).

\fbox{\textbf{Subcase $N=0$.}} Then $(g-1)(g+1)=0$ and  we may
assume $g=1$ (otherwise  the transformation $x\to -x$ and $y\to
y-x$ can be applied). So, $g=1$ and by translation of the origin
of coordinates to the point ($-c/2,-f/2$)   we obtain the systems
\be\label{S1:2}
  \dot x=k + dy + x^2,\quad  \dot y=l + ex +y^2.
\ee
For this  systems we have $\  H_4=96(d^2+e^2)$\ and
$$
 B_2=648\big[e^2(4k-4l-e^2)x^4+2d^2e^2(2x^2-3xy+2y^2)
-d^2(4k-4l+d^2)y^4\big].
$$
Therefore, the condition $B_2=0$ yields
$de=e(4k-4l-e^2)=d(4k-4l+d^2)=0$ and we may assume $d=0$, since
for $d\ne0$, $e=0$ the substitution  $x\leftrightarrow y$,
$d\leftrightarrow e$ and $k\leftrightarrow l$  can be applied.
Then the condition $B_2=0$ yields $e(4k-4l-e^2)=0$.

We claim that for $M_{{}_\IL}<4$, $e\ne0$  (i.e. $H_4\ne0$)  is
necessary. Indeed, let us suppose $e=0$. Then for
systems~\eqref{S1:2} with $d=e=0$ calculations yield:
$$
{\cal H}=\gcd\left({\cal E}_1,{\cal E}_2\right)= (X^2 + k Z^2)(Y^2
+ lZ^2),
$$
i.e. \ $\deg\,{\cal H}>3$ and our claim is proved. Thus, $H_4\ne0$
(i.e. $e\ne0$) and replacing $e$ by $2e$ the condition $B_2=0$
yields $l=k-e^2$.  Hence, we get the systems
\be\label{CS1:3}
  \dot x=k +  x^2,\quad  \dot y=k-e^2 + 2ex +y^2
\ee
for which calculations yield:
$$
\bal
&{\cal E}_1= -4(Y^2 - eYZ + eZX + kZ^2){\cal H},\quad
{\cal E}_2=-(X + Y - eZ)(Y^2 + 2 eX Z - e^2Z^2 + k Z^2){\cal H},\\
\eal
$$
where $ {\cal H}= (Y-X+eZ)(X^2+kZ^2)$. Thus, $\deg{\cal H}=3$ and
we shall show that for all values given to the parameters $k$ and
$e$ ($e\ne0$) the degree of $\gcd({\cal E}_1,{\cal E}_2)$ remains
3.   Indeed, since the common factor of the polynomials ${\cal
E}_1/{\cal H}$ and ${\cal E}_2/{\cal H}$ must depend on $Y$,
according to Lemma \ref{Trudi:2} it is sufficient to observe that
$ \Res_Y({\cal E}_1/{\cal H},\ {\cal E}_2/{\cal H})=-64e^2Z^2(X^2
+ k Z^2)^2\ne0$ since $e\ne0$. The systems (\ref{CS1:3}) possess
the invariant lines $x-y=e$, and $x^2+k=0$. Since for these
systems $H_8=-2^93^3e^2k$ we obtain that for $k\ne0$
$\sign(k)$=$-\sign(H_8)$.

\textbf{1)} If $H_8>0$ then $k<0$ and we may assume $k=-1$  via
the Remark \ref{rem:transf} ($\gamma=-k$, $s=1/2$). Then systems
(\ref{CS1:3}) possess the following real singular points:
$(1,\pm(e-1))$ and $(-1,\pm(e+1)),$ which are distinct if
$e\ne\pm1$. On the other hand for systems (\ref{CS1:3}) with
$k=-1$ we calculate $H_9=-2^{14}3^2(e^2-1)^2$ and hence  we obtain
Config 4.5 for $H_9\ne0$ and Config 4.6 for $H_9=0$. In the last
case we may consider $e=1$, otherwise the substitution $x\to -x$,
$y\to -y$ and $t\to -t$ can be applied.

\textbf{2)} Assume $H_8<0$. Then $k>0$ and we may assume $k=1$ via
the Remark \ref{rem:transf} ($\gamma=k$, $s=1/2$). Then the
invariant lines $x^2+1=0$ are imaginary and, moreover,  systems
(\ref{CS1:3}) do not possess  any real singular point. Thus we get
Config. 4.13.

\textbf{3)} For $H_8=0$ since $e\ne0$, we obtain $k=0$. We may
assume $e=1$ via the Remark \ref{rem:transf} ($\gamma=e$, $s=1$).
This leads to the system
\be\label{CS1:3a}
  \dot x= x^2,\quad  \dot y=-1 + 2x +y^2
\ee
for which  $ {\cal H}= X^2(Y-X+Z)$. According to Lemma \ref{lm3}
the line $x=0$ could  be of multiplicity two and the perturbations
({\it 4.22${}_\varepsilon$}) from Table 3 show  this. Taking into
consideration that both singular points $(0,\pm 1)$ of system
(\ref{CS1:3a}) are double points we get Config. 4.22.

\medskip
It remains to note  that via the transformation $\ x_1=-x$ and
$y_1=y-x$ the systems (\ref{CS1:3}) can be brought to the systems
\be\label{S:N0}
  \dot x_1=-k -  x_1^2,\quad  \dot y=-e^2 -2ex_1-2x_1y_1 +y_1^2.
\ee
On the other hand we observe that systems \eqref{S:N0} can be
obtained from systems \eqref{S1:4} by setting $g=-1$   and
changing $u\to e$ and  $k\to-k$. For systems (\ref{CS1:3}) we have
$H_7=0$ and $B_3=-12e^2x^2(x-y)^2\ne0$,  and   for systems
\eqref{S1:4} we have $H_4=48(g-1)\big[k(g+1)^2+4gu^2\big]$. Since
for systems \eqref{S1:4} in $\QSL_{\bf4}$
$B_3=-3\big[k(g+1)^2+4gu^2\big]x^2y^2\ne0$  we obtain $H_4\ne0$.
Therefore we conclude, that the representatives of the orbits in
the case $N=0$ for systems \eqref{S:N0}  can be included as
respective particular cases of the representatives corresponding
to systems \eqref{S1:4} for which  $N\ne0$ as indicated in Table 2
for Config. 4.10, 4.13 and 4.22.

\subsection{Systems with the divisor\ $D_S(C,Z)=1\cdot \omega_1^c+1\cdot \omega_2^c+1\cdot \omega_3$}

According to Lemma \ref{lm_3:2} in this case we shall consider the
canonical systems $(\SSS_{I\!I})$ for which we have:
\be\label{val:Kappa-N}
\ba{l}
 \theta= 8(h+1)[(h-1)^2+g^2],\\ N=(g^2-2h+2)x^2 +2g(h+1)xy+(h^2-1)y^2.
\ea
\ee
\brm\label{rm:B_3=0}
We note that  two infinite points of the systems~$(\SSS_{I\!I})$
are imaginary. Therefore by Lemma \ref{lm:BGI} we conclude that a
system of this class could belong to $\QSL_{\bf4}$ only if for
this system the condition $B_3=0$ holds.
\erm
In what follows we shall assume that for a system $(\SSS_{I\!I})$
the condition $B_3=0$ is fulfilled.

\subsubsection{The case
$\theta\ne0$}

The condition $\theta\ne0$ yields $(h+1)\ne0$ and we may assume
$c=d=0$ via the translation $x\to x-d/(h+1)$ and $y\to y+(2dg-
c(h+1))/(h+1)^2$. Thus we obtain the  systems
\be\label{S2:1}
  \dot x=k +  gx^2+(h+1)xy,\quad
  \dot y=l + ex+fy -x^2+ gxy+ hy^2,
\ee
for which we  have: Coefficient$[B_3,\ y^4]=-3k(h+1)^2$. Since
$h+1\ne0$ the condition $B_3=0$ implies $k=0$ and then we obtain
\beq
B_3=
3[ef(h+1)+2gl(h-1)-f^2g]x^2(x^2-y^2)+6[f^2+efg-e^2h+l(h-1)^2-g^2l]x^3y.
\eeq
Hence, the condition $B_3=0$ yields the following system of
equations
\be\label{Eq:1}
\ba{l} Eq_1\equiv ef(h+1)+2gl(h-1)-f^2g=0,\\
Eq_2\equiv f^2+efg-e^2h+l\big[(h-1)^2-g^2]=0.
\ea
\ee
Because $\theta\ne0$ the conditions $gl(h-1)=(h-1)^2-g^2=0$ are
impossible. Hence   we calculate
$$
\Res_l(Eq_1,\,Eq_2)=(e-eh+ f g)[2egh+ f(h^2-1- g^2)].
$$
On the other hand for  systems (\ref{S2:1}) with $k=0$ we obtain:\
$ H_7=4(h+1)[e(1-h)+fg] $\ and we shall consider two subcases:
$H_7\ne0$ and $H_7=0$.

\textbf{1)} If $H_7\ne0$ then the equality $\Res_l(Eq_1,\,Eq_2)=0$
yields $2egh= f(g^2+1-h^2)$. Since  $\theta\ne0$  from
\eqref{val:Kappa-N} we have $(gh)^2+(g^2+1-h^2)^2\ne0$ then
without loss of generality we may set: $e=(g^2+1-h^2)u$ and
$f=2ghu$ where $u$ is a new parameter. Therefore from \eqref{Eq:1}
we obtain
$$
g(h-1)\big[l-hu^2(g^2+(h+1)^2)\big]=0=\big[(h-1)^2-g^2\big]\big[l-hu^2(g^2+(h+1)^2)\big]
$$
and hence, $l=hu^2[g^2+(h+1)^2]$ which leads to the systems:
\be\label{S2:2}
\bal
 & \dot x=  gx^2+(h+1)xy,\\
 & \dot y=hu^2[g^2+(h+1)^2] + u(g^2+1-h^2)x+2ghuy -x^2+ gxy+ hy^2.
\eal
\ee
For these systems we have $H_7=4u(h+1)^2[g^2+(h-1)^2]\ne0$. Then
$u\ne0$ and we may assume $u=1$ via Remark \ref{rem:transf}
($\gamma=u$,$s=1$).

For systems \eqref{S2:2} with $u=1$ calculations yield:
$$
{\cal H}=\gcd\left({\cal E}_1,{\cal E}_2\right)=X
[X^2+Y^2-2(h+1)XZ +2g\,YZ+(g^2+(1+h)^2)Z^2].
$$
Thus  $\deg\,{\cal H}=3$ and for every specific system in the
family \eqref{S2:2} with $u=1$ we have  $M_{{}_\IL}=4$. Since
$\theta\ne0$, by Lemma~\ref{gcd:3} $M_{{}_\IL}$ cannot be 5. In
this case systems
\eqref{S2:2} possess three invariant lines: $x=0$ , $y\pm
i\big[x-(h+1)\big]+g=0$ and  four singularities:
$$
\big(0,\,-g\pm i(h+1)\big),\quad (h+1,\,-g),\quad
\big(-h(h+1),\,gh \big).
$$
Due to the condition $\theta\ne0$ (i.e. $h+1\ne0$) all points are
distinct and the intersection point $(h+1,\, -g)$ of the imaginary
lines  is not placed on the line $x=0$. This leads to the  Config.
4.2.

\textbf{2)} Assume  $H_7=0$. Then  $e(h-1)=fg$\  and since the
condition $\theta\ne0$ yields $g^2+(h-1)^2\ne0$, we may assume
$e=gu$ and $f=(h-1)u$, where $u$ is a new parameter. Then from
(\ref{Eq:1}) we have
$$
 g(h-1)(l+u^2)=0=[(h-1)^2-g^2](l+u^2).
$$
By $g^2+(h-1)^2\ne0$ evidently we obtain $l=-u^2$ and this leads
to the  systems:
\be\label{S2:3}
  \dot x=  gx^2+(h+1)xy,\quad
  \dot y=-u^2 + gu\,x+u(h-1)y -x^2+ gxy+ hy^2.
\ee
For these systems  we have: $ {\cal H}=\gcd\left({\cal E}_1,{\cal
E}_2\right)= 2 X \big[X^2 +(Y+uZ)^2\big].$ Then for a specific
systems in this family $M_{{}_\IL}\ge4$, and by Lemma~\ref{gcd:3}
due to the condition $\theta\ne0$, $M_{{}_\IL}$ cannot be 5. The
systems (\ref{S2:3}) possess three invariant lines: $x=0$ and
$x\pm i(y+u)=0$ and the following singular points:
$$
(0,-u),\quad (0,u/h), \quad \Big(\frac{u(h+1)}{g\pm i(h+1)},\,
-\frac{gu}{g\pm i(h+1)}\Big)
$$
We observe that all three lines have the common point of
intersection: $(0,-u)$. Moreover if $u=0$ this point becomes of
multiplicity 4 and the point $(0,u/h)$ tends to infinity when $h$
tends to zero. On the other hand for systems (\ref{S2:3}) we have
$\mu=-32h\big[g^2+(h+1)^2\big]$ and $H_{9}=2^83^2u^8(h+1)^8$. Thus
in the case $H_7=0$ we obtain:  Config.~4.6 for $\mu H_9\ne0$,\
Config.~4.8 for $\mu\ne0$, $H_9=0$\ and Config.~4.7  for $\mu=0$
(in this case $u\ne0$ otherwise we get  degenerate systems). It
remains to note that if  $u\ne0$  we may assume $u=1$ via Remark
\ref{rem:transf} ($\gamma=u,$ $s=1$).

\subsubsection{The case $\theta=0$}

 According to (\ref{val:Kappa-N}) we have
$(h+1)[(h-1)^2+g^2]=0$ and we shall consider two subcases $N\ne0$
and $N=0$.

\textbf{Subcase $N\ne0$.} Then by (\ref{val:Kappa-N}) the
condition $\theta=0$ yields $h=-1$ and  in addition we may assume
$f=0$ due to the translation: $x\to x$ and $y\to y+f/2$. Hence, we
obtain the systems
\be\label{S2:4}
  \dot x=k +cx+dy+ gx^2,\quad
  \dot y=l + ex -x^2+ gxy-y^2,
\ee
for which by Remark \ref{rm:B_3=0} the condition $B_3=0$ must be
satisfied. Calculation yields
$$
\text{Coefficient}[B_3,y^4]=-3\,d^2g,\quad H_7=4d(g^2+4).
$$
So the condition $B_3=0$ implies  $dg=0$. We claim that for
systems \eqref{S2:4} to be in the class $\QSL_{\bf4}$ the
condition $d\ne0$ must be fulfilled. Indeed, suppose $d=0$. Then
we have:
$$
B_3=3(2ce-4gl+4k-g^2k)x^2(x^2-y^2)-6(c^2-e^2-4l+lg^2-4gk)x^3y
$$
and hence the condition $B_3=0$ yields a   system of two linear
equations with respect to the parameters $k$ and $l$.
 As its determinant equals
$(g^2+4)^2\ne0$  we obtain the unique solution:
\be\label{val:param-kl}
\!\!\ba{l}
k=2(2 c +  e g)(c g-2e)/(g^2+4)^2,\\[2mm]
l=(4 c^2  - 4 e^2 +  8 c e g - c^2 g^2 + e^2 g^2)/(g^2+4)^2.
\ea
\ee
Then  for    systems~(\ref{S2:4}) with  these  values of the
parameters $l$ and $k$ and with $d=0$ calculations yield:
\beq
&&{\cal H}= 2 \Big[(4 + g^2)X - (2 e - c g) Z\big] \Big[g(g^2+4)X
+ 2(2 c  +
           e g) Z\Big]\times\\
&&  \Big[(4 + g^2)( X^2 +  Y^2)
-2(2e-cg)XZ+2(2c+eg)YZ+(c^2+e^2)Z^2\Big],
\eeq
i.e. in this case the family of systems \eqref{S2:4} possess four
affine lines, which implies $M_{{}_\IL}\ge5$. Our claim is proved.

Let us assume  $d\ne0$, i.e. $H_7\ne0$. Then  $g=0$ and we may
assume $e=0$ via the translation: $x\to x+e/2$, $y\to y$.  After
that for systems
\eqref{S2:4} calculations yield:
$$
B_3=12k x^2(x^2-y^2)-6(c^2-4l+d^2)x^3y.
$$
Therefore the condition $B_3=0$ yields  $k=0$ and $l=(c^2+d^2)/4$.
Then by replacing  $c$ with $2c$ and   $d$ with $ 2d$ and assuming
$d=1$ (due to the Remark \ref{rem:transf} ($\gamma=d,$ $s=1$) we
get the systems:
\be\label{S2:5}
  \dot x= 2cx+2y,\quad
  \dot y=c^2+1  -x^2-y^2,
\ee
for which we calculate
\beq
 &&{\cal E}_1=2\big[ X^2 - 2 c X Y -  Y^2 +  (c^2+ 1)(2 X+ Z)Z\big]{\cal H},\\
 &&{\cal E}_2= (c X +  Y) [X^2 + Y^2 + 2  X Z - 2 c Y Z +
        (c^2  + 1) Z^2]{\cal H},\\
&&{\cal H}=\gcd\left({\cal E}_1,{\cal E}_2\right)=2 Z [(X-Z)^2 +(Y
+ c Z)^2].
\eeq
Thus, for this family of systems we obtain $\deg\,{\cal H}=3$.
According to Lemma~\ref{Trudi:2}  to have an additional common
factor of ${\cal E}_1$ and ${\cal E}_2$ it is necessary
 that the following  holds in $\R[X,Z]$:
$$
\Res_Y({\cal E}_1/{\cal H},\ {\cal E}_2/{\cal H})= 32  (c^2 + 1)^2
(X + d Z)^6 =0.
$$
However this is impossible and hence the systems (\ref{S2:5})
belong to the class $\QSL_{\bf4}$ for any value of the parameter
$c$. Note that these systems possess the invariant lines $y+c=\pm
i(x-1)$. Since  $Z\mid{\cal H}$,  by Corollary \ref{Mult:Z=0} we
conclude that the line $Z=0$ could be of multiplicity two and this
is confirmed by the perturbations ({\it IV.27${}_\varepsilon$})
from Table 3.

 The systems~\eqref{S1:8} possess two finite distinct
singular points: $(-1,\,c)$ and $(1,-c)$ and the second one is the
point of intersection of the imaginary lines. This leads to the
Config. 4.27.

\textbf{Subcase $N=0$.} Then from (\ref{val:Kappa-N}) we have
$g=h-1=0$  and without loss of generality we may assume $c=d=0$
via the translation $x\to x-d/2$, $y\to y-c/2$. Hence  we obtain
the systems
$$
\dot x=k + 2xy,\quad
  \dot y=l + ex +fy -x^2-y^2,
$$
for which calculations yield:
$$
 B_3=6\,\left[(ef-2k)x^4+(f^2-e^2)x^3y-(4k+ef)x^2y^2-2ky^4\right].
$$
Therefore, the condition $B_3=0$ yields $k=e=f=0$ and this leads
to the following  systems
$$  \dot x= 2xy,\qquad
   \dot y= l -x^2 + y^2,
$$
for which  we obtain:\ ${\cal H}=\gcd\left({\cal E}_1,{\cal
E}_2\right)= 2 X \left[X^4 +
        2 X^2 (Y^2 - l Z^2) + (Y^2 + l Z^2)^2\right].
$\ Hence, for $N=0=B_3$ we have $\deg\,{\cal H}=5$ and  systems
$(\SSS_{II})$ cannot belong to the class $\QSL_{\bf4}$.

\subsection{Systems with\ $D_S(C,Z)=2\cdot \omega_1+1\cdot \omega_2$}

 For this case we shall later  need the following $T$-comitant.
\bnot\label{not_H11} Let us denote
$\quad H_{11}(a,x,y)=8
H\big[(C_2,D)^{(2)}+(D,D_2)^{(1)}\big]+3H_2^2.$
\enot
By Lemma \ref{lm_3:2} systems (\ref{2l1}) can be brought via a
linear transformation to the canonical form $(\SSS_{I\!I\!I})$ for
which we have:
\be\label{val: kappa-N-3}
\theta =  -8h^2(g-1),\quad  \mu=32gh^2,\quad N=(g^2-1)x^2
+2h(g-1)xy+h^2y^2.
\ee

\vspace{-5mm}
\subsubsection{The case $\theta\ne0, B_3=0$}

Then $h(g-1)\ne0$ and  we may assume $h=1$ due to the substitution
$y\to y/h$. Then  via the translation $x\to x-d$ and $y\to y+2dg-
c$ we may assume $c=d=0$. Thus we obtain the systems
\be\label{S3:1a}
  \dot x=k +  gx^2+ x y,\quad
  \dot y=l + ex+fy + (g-1)xy+ y^2,
\ee
for which   calculation yields
$$
B_3=-3\big[l(g-1)^2+ef(1-g)+e^2\big]x^4 +3(l-2k-2gk)x^2y^2-6kxy^3.
$$
Hence  the condition $B_3=0$ implies $k=l=0$ and $e(f-f g+e)=0$.
On the other hand for the systems~(\ref{S3:1a}) we have $ H_7=-4(f
- f g + e)$ and we shall consider two subcases: $H_7\ne0$ and
$H_7=0$.

\textbf{Subcase $H_7\ne0$.} Then $f-fg+e\ne0$ and the condition
$B_3=0$ yields $e=0$. Therefore the condition $H_7=4f(g-1)\ne0$
implies $f\ne0$ and we may assume $f=1$ via Remark
\ref{rem:transf} ($\gamma=f,$ $s=1$).  This leads to the systems:
\be\label{S3:1}
  \dot x=  gx^2+ x y,\quad
  \dot y= y + (g-1)xy+ y^2,
\ee
for which we calculate:
\beq
 &&{\cal E}_1= 2\big[g(g-1) X^2 + 2 g X Y +  Y^2 + 2 g X Z -
    g Z^2\big]{\cal H},\\
 &&{\cal E}_2= (g X +  Y)^2 (g X-X  +  Y +  Z){\cal H},\quad
 {\cal H}=\gcd\left({\cal E}_1,{\cal E}_2\right)=X^2Y.
\eeq
Thus, we obtain $\deg\,{\cal H}=3$ and we claim that for systems
\eqref{S3:1} we cannot have an additional common factor of ${\cal E}_1$
and ${\cal E}_1$ for any choice of value for $g\not\in\{0,1\}$ .
Indeed, by Lemma~\ref{Trudi:2} to have such a common factor  it is
necessary that the following holds in $\R[X,Z]$:
$$
\Res_Y({\cal E}_1/{\cal H},\ {\cal E}_2/{\cal H})= 8 (1 - g) g^2
(X -  Z)^6 =0.
$$
As $g\not\in\{0,1\}$ this is impossible.

We observe that  the systems \eqref{S3:1} possess the invariant
affine lines $x=0$ and $y=0$ and the line $x=0$ is a double one.
This is confirmed by the perturbations   ({\it
IV.25${}_\varepsilon$})   from Table 3.

Note that  the systems \eqref{S3:1} have three distinct  finite
singular points: $(0,0)$ (which is double),\ $ (0,\,-1)$ and
$(1,-g)$.  Therefore we obtain Config. 4.25.

\textbf{Subcase $H_7=0$.} Then $e=f(g-1)$ and we obtain the
systems:
\be\label{S3:2}
  \dot x=  g\,x^2+ x y,\quad
  \dot y= f(g-1)\,x+  f\, y + (g-1)\,xy+ y^2,
\ee
for which calculations yield:
\beq
 &&{\cal E}_1= 2\big[g(g-1) X^2 + 2 g  X Y +  Y^2
 +f(1- g) X Z - f Y Z)\big]{\cal H},\quad\
{\cal H}=X^2(Y+fZ),\\
 &&{\cal E}_2=3 (g X +  Y)^2 (gX-X   +  Y ){\cal H},\quad
\Res_Y({\cal E}_1/{\cal H},\ {\cal E}_2/{\cal H})= 8 (1 - g)
 (gX - f Z)^2X^4.
\eeq
Thus, we obtain $\deg\,{\cal H}=3$ and since $\Res_Y({\cal
E}_1/{\cal H},\ {\cal E}_2/{\cal H})\not\equiv0$, according to
Lemma~\ref{Trudi:2} we cannot  have an additional common factor
 of ${\cal E}_1$ and ${\cal E}_1$ for any specific value  of the
 parameters.

Note that the systems  \eqref{S3:2} have the finite singular
points: $(0,0)$ (which is double),\ $ (0,\,-f)$ and $(f/g,-f)$.
For $f\ne0$ (then we can assume $f=1$ due to the Remark
\ref{rem:transf} ($\gamma=f,$ $s=1$)) these points are distinct
and finite if $g\ne0$. When $g\to0$ the point $(f/g,-f)$ tends to
infinity. For $f=0$ (then $g\ne0$, otherwise we get a degenerate
system) the point $(0,0)$ is of   multiplicity four.

On the other hand for the systems \eqref{S3:2} we have: $D=-fx^2y$
and $\mu=32g$, i.e. the conditions $f=0$ and $g=0$ are captured by
the $T$-comitants $D$ and $\mu$, respectively.

We observe that  the systems \eqref{S3:2} possess the invariant
affine lines $x=0$ and $y+f=0$ and the line $x=0$ is a double one.
This is confirmed in the case $f\ne0$ (respectively  $f=0$) by the
perturbations  ({\it IV.21${}_\varepsilon$}) and ({\it
IV.26${}_\varepsilon$}) (respectively  ({\it
IV.20${}_\varepsilon$})) from Table 3.

Thus for $\theta\ne0$, $B_3=0$ and $H_7=0$, if   $D\ne0$ and
$\mu\ne0$  (respectively $D\ne0$ and $\mu=0$;\ $D=0$) we obtain
Config. 4.21 (respectively  Config. 4.26; Config. 4.20).

\subsubsection{The case $\theta=0=B_2$}

\fbox{\textbf{Subcase $\mu\ne0$.}} From (\ref{val: kappa-N-3}) we
obtain $h\ne0$, $g=1$  and then  we may assume $h=1$ due to the
change $y\to y/h$. Moreover,  we may assume $c=d=0$ via the
translation $x\to x-d$ and $y\to y+2d- c$. So, we obtain the
canonical systems
\be\label{S3:3}
  \dot x=k +  x^2+ x y,\quad
  \dot y=l + ex+fy +  y^2,
\ee
for which calculation yields:
$$
B_2=-648 \Big[e^2(l-4k)x^4 - 4 e^2kx^3y+  k^2 y^4\Big], \quad
H_7=-4e.
$$
Hence, the condition $B_2=0$ yields $k=el=0$. We claim, that if
$e\ne0$ (i.e. $H_7\ne0$) then no    system of the family
 \eqref{S3:3} can belong to  $\QSL_{\bf4}$.  Indeed, supposing $e\ne0$
(then we may assume $e=1$ due to the Remark \ref{rem:transf}
($\gamma=e,$ $s=1$)), we have $l=0$ and we obtain the systems
$$
  \dot x=  x^2+ x y,\quad
  \dot y= x+fy +  y^2
$$
 for which we calculate:
\beq
 &&{\cal E}_1=2[X^2Z+2 XY^2 +Y^3+(2f-1)X Y Z +(1-f)XZ^2+f(1-f)YZ^2]{\cal H},\\
 &&{\cal E}_2= (X +  Y)^2 ( Y^2 +  X Z + f Y Z){\cal H},\quad
{\cal H}=X^2.
\eeq
So, $\deg\,{\cal H}=2$ and according to Lemma~\ref{Trudi:2} to
have an additional common factor of ${\cal E}_1$ and ${\cal E}_2$
for some specific value of $f$ the condition   $ \Res_Y({\cal
E}_1/{\cal F},\ {\cal E}_2/{\cal F})= 16 X^4 Z^2 (X - f Z + Z)^6=0
$\ must hold in $\R[X,Z]$ and this is  impossible.

 Assuming $e=0$ (then  $H_7=0$ )  we get the following
systems:
\be\label{S3:4}
  \dot x= x^2+ x y,\quad
  \dot y=l + fy +  y^2
\ee
for which calculations yield:
\beq
 &&{\cal E}_1=2(2 X^2 +  X Y - f X Z +  l Z^2){\cal H},\quad
{\cal E}_2=  X (X +  Y)^2{\cal H},
\eeq
where ${\cal H}=X(Y^2+fYZ+lZ^2)$. Thus, we obtain $\deg\,{\cal
H}=3$ and to have an additional common factor of ${\cal E}_1$ and
${\cal E}_2$ for some specific value of $(f,l)$, according to
Lemma~\ref{Trudi:2} it is necessary and sufficient that the
following condition holds:
$$
\Res_X({\cal E}_1/{\cal F},\ {\cal E}_2/{\cal F})= 8 l Z^2 ( Y^2 +
f Y Z + l Z^2)^2 \equiv0,
$$
which is equivalent to $l=0$. On the other hand for the
systems~\eqref{S3:4} we have $B_3=3lx^2y^2$ and hence   for these
systems to belong to $\QSL_{\bf4}$ we must have $B_3\ne0$.

Assume  $l\ne0$ (i.e. $B_3\ne0$). Then the systems~\eqref{S3:4}
possess three invariant lines: $x=0$ and $y^2+fy+l=0$ and the
position of these two lines depends on the values of the affine
invariant $H_{10}=8(f^2-4l)$.

 \textbf{1)} Assume $H_{10}>0$. Then we may set $f^2-4l=4u^2\ne0$
($u$ is a new parameter) and replacing $f$ by $2f$ we have
$l=f^2-u^2$. Assuming $u=1$ via Remark \ref{rem:transf}
($\gamma=u,$ $s=1$) we get the systems
\be\label{S3:5}
  \dot x= x^2+ x y,\quad
  \dot y=(y+f)^2-1.
\ee
  These systems possess three real distinct invariant lines
($x=0$ and $y+f=\pm1$) and the following singular points:
$$
(0,1-f),\quad (0,-1-f),\quad (-1+f,1-f),\quad (1+f,-1-f).
$$
Since for systems~\eqref{S3:5}  we have $B_3=3(f^2-1)x^2y^2)\ne0$,
we obtain $|f|\ne1$ and, hence, all singular points  are distinct.
Thus we  get Config. 4.11.

\textbf{2)} If $H_{10}<0$ then  as above assuming
$f^2-4l=-4u^2\ne0$ and replacing $f$ by $2f$  and $u=1$ we get the
systems
\be\label{S3:6}
  \dot x= x^2+ x y,\quad
  \dot y=(y+f)^2+1
\ee
which possess one real  ($x=0$) and two imaginary ( $y+f=\pm i$)
 invariant lines. Since all singular points are imaginary we
 obtain Config 4.23.

\textbf{3)} Assume now   $H_{10}=0$, i.e. $f^2=4l$. Then replacing
$f$ by $2f$ we obtain the systems
\be\label{S3:7}
  \dot x= x^2+ x y,\quad
  \dot y=(y+f)^2
\ee
for which we have $B_3=3f^2x^2y^2\ne0$. Then $f\ne0$ and we may
assume $f=1$ due to the Remark \ref{rem:transf} ($\gamma=f,$
$s=1$).  We observe that the system~\eqref{S3:7} has  a simple
invariant line ($x=0$) and a double one ($y=-1$). Moreover, these
systems have 2 double singular points: $(0,-1)$ and $(1,-1)$. Thus
we   obtain Config. 4.23.

\smallskip

 \fbox{\textbf{Subcase $\mu=0$.}} Since $\theta=0$  from
(\ref{val: kappa-N-3}) we obtain $h=0$ and  for the  systems
$(\SSS_{I\!I\!I})$ we calculate
$$
 N=(g^2-1)x^2, \quad H_7=4d(g^2-1).
$$

\fbox{\textbf{I.}} If $N\ne0$  then $g-1\ne0$ and  we may assume
$e=f=0$ via the translation $x\to x+f/(1-g)$ and $y\to y+e/(1-g)$.
This leads to the  systems
\be\label{S3:7a}
  \dot x=k + cx + dy + gx^2,\quad
  \dot y=l +(g-1)xy,
\ee
for which  calculations yield:
$$
B_2=648 d\Big[cl(g-1)^3x^4 + 4 dlg(g-1)^2x^3y- d^3 g^2 y^4\Big]
$$
and we shall consider two subcases:  $H_7\ne0$ and $H_7=0$.

\textbf{1)} Suppose  $H_7\ne0$. Then $d\ne0$ and we may assume
$d=1$ due to the substitution $y\to y/d$. Then by $g-1\ne0$ the
condition $B_2=0$ yields $g=cl=0$.

We claim that the systems
\eqref{S3:7a} with $d=1$ could   be in the class $\QSL_{\bf4}$ only if the condition
$c=0$ is fulfilled. Indeed, supposing  $c\ne0$ from $B_2=0$ we
obtain $g=l=0$ and this leads to the systems
$$
  \dot x=k + cx + y,\quad
  \dot y= -xy.
$$
 For these  systems   calculations yield:
\beq
 &&{\cal E}_1=-2(c X^3 + k X^2 Z - c X Y Z -  Y^2 Z -
    c k X Z^2 - 2  k Y Z^2 - k^2 Z^3){\cal H},\\
 &&{\cal E}_2=-X Z ( c X   +  Y+k Z)^2{\cal H},\quad {\cal H}=YZ.
\eeq
 So, $\deg\,{\cal H}=2$ and since $c\ne0$, to have an additional
 common factor of ${\cal E}_1$ and
${\cal E}_2$ for some specific value of $(c,k)$, according to
Lemma~\ref{Trudi:2} it is necessary    that the following  holds
in $\R[Y,Z]$:
$$
\Res_X({\cal E}_1/{\cal H},\ {\cal E}_2/{\cal H})= -8 c^2 Y^2 Z^4
( Y + k Z)^6=0,
$$
which is impossible since $c\ne0$.

 Let us now assume  $c=0$. Then we obtain
the  systems:
$$
  \dot x=k +  y,\quad
  \dot y= l-xy,
$$
 for which we calculate:
\beq
 &&{\cal E}_1=2(-k X^2 Y +  Y^3 -  l X Y Z + 2  k Y^2 Z +
    k l X Z^2 + k^2 Y Z^2 + l^2 Z^3){\cal H},\quad {\cal H}=Z^2,\\
 &&{\cal E}_2=-( Y + k Z)^2 (X Y - l Z^2){\cal H},\quad
 \Res_Y({\cal E}_1/{\cal H},\ {\cal E}_2/{\cal H})=-8  l Z^6
(k X +  l Z)^6.
\eeq
Thus, we have $\deg\,{\cal H}=2$ and   to have an additional
common factor of ${\cal E}_1$ and ${\cal E}_2$ for some specific
value of $(k,l)$, according to Lemma~\ref{Trudi:2} it is necessary
 $\Res_Y({\cal
E}_1/{\cal H},\ {\cal E}_2/{\cal H})\equiv0$. This condition is
equivalent to $l=0$ in which case the additional common factor is
$Y$. Thus  the systems~\eqref{S3:7a} with  $d=1$ and $g=cl=0$
belong  to the class $\QSL_{\bf4}$ if and only if $c=l=0$. Since
for these systems we have $B_3=-3x^2(lx^2-cy^2)$ we conclude that
the condition $c=l=0$ is equivalent to $B_3=0$.

Therefore,  for $\theta=\mu=0$, $N\ne0$, $H_7\ne0$ and $B_3=0$
(then $B_2=0$) we obtain the systems:
\be\label{S3:8}
  \dot x=k + y,\quad
  \dot y= -xy,
\ee
for which  we have\ ${\cal H}= \gcd\left({\cal E}_1,{\cal
E}_2\right)=YZ^2$\ and
\beq
 &&{\cal E}_1=(-k X^2 + d^2 Y^2 + 2 d k Y Z + k^2 Z^2){\cal H},\quad
{\cal E}_2=-3 X (d Y + k Z)^2{\cal H}.
\eeq
Hence  $\deg\,{\cal H}=3$ and according to Lemma~\ref{Trudi:2} to
have an additional common factor of ${\cal E}_1$ and ${\cal E}_2$
for some specific value of $k$ it is necessary to have $
\Res_Y({\cal E}_1/{\cal H},\ {\cal E}_2/{\cal H})=8 k^2
X^6\equiv0$, i.e. $k=0$. This case is ruled out as the
corresponding system (\ref{S3:8}) is degenerate.

The systems~(\ref{S3:8}) possess the invariant line $y=0$ and due
to Corollary \ref{Mult:Z=0} since $Z^2\mid{\cal H}$ the line $Z=0$
could be of the multiplicity three. And this is confirmed by the
perturbations ({\it IV.35${}_\varepsilon$}) from Table 3. So in
this case we obtain Config. 4.35.

\textbf{2)} Assume $H_7=0$. Then $d=0$ and  this implies $B_2=0$.
Hence we get the systems
\be\label{S3:9}
  \dot x=k + cx +  gx^2,\quad
  \dot y=l +(g-1)xy,
\ee
 for which  calculations yield:
\beq
 &&{\cal E}_1=2\big[(g-1)X^2 Y  + l(g+1) X Z^2  + k(1-g) Y Z^2 + c l Z^3\big]{\cal H},\\
 &&{\cal E}_2= (g X^2 + c X Z + k Z^2) (-X Y + g X Y + l
 Z^2){\cal H},\quad
{\cal H}=  gX^2+cXZ+kZ^2,\
\eeq
Thus, we have $\deg\,{\cal H}=2$, but for systems~(\ref{S3:9}) to
be in the class $\QSL_{\bf4}$ we need an additional common factor
of  of ${\cal E}_1$ and ${\cal E}_2$. Since $g-1\ne0$ according to
Lemma~\ref{Trudi:2} we can obtain such a common factor if at least
one of the following two identities holds:
$$
\bal
& \Res_X({\cal E}_1/{\cal F},\ {\cal E}_2/{\cal
F})=8(g\!-\!1)\Big[c^2-k(g+1)^2\Big]\Big[k(g-1)^2 Y^2\!+ \!c\,
l(1\!\!-\!g) Y Z \! + \!g\, l^2 Z^2\Big]^2 Y Z^6
 \equiv0,\\
&\Res_Y({\cal E}_1/{\cal F},\ {\cal E}_2/{\cal F})=2\,l (g-1) (g
X^2 + c XZ + k Z^2)^2 Z^2  \equiv0.
\eal
$$
So  to have an additional common factor of  ${\cal E}_1$ and
${\cal E}_2$ one of the following   three conditions must be
fulfilled:
$$
(i)\,\ l=0;\qquad (ii)\,\ c^2-k(g+1)^2=0;\qquad (iii)\,\
k=c\,l=g\,l=0.
$$
Since for $l\ne0$ the condition $(iii)$ leads to the degenerate
systems~(\ref{S3:9}), only the conditions $(i)$ and $(ii)$ remain
to be examine.

On the other hand, for the systems~(\ref{S3:9}) we have
$$
B_3=-3\,l(g-1)^2x^4,\qquad H_6=128(g-1)^4[k(g+1)^2-c^2]x^6
$$
and we shall consider two subcases: $B_3=0$ and $B_3\ne0$,
$H_6=0$.

\textbf{a)} If $B_3=0$ (i.e. $l=0$)  we get the systems:
\be\label{S3:10}
  \dot x=k + cx+gx^2,\quad
  \dot y=(g-1)xy
\ee
with $k\ne0$, otherwise they are degenerate.   For these systems
calculation yields:
$$
\bal
 &{\cal E}_1=2(X^2-k Z^2){\cal H},\qquad
{\cal H}= \gcd\left({\cal E}_1,{\cal E}_2\right)=(g-1)(gX^2+cXZ+kZ^2)Y,\\
 &{\cal E}_2=X (g X^2 + c X Z + k Z^2){\cal H},\qquad
 \Res_X({\cal E}_1/{\cal H},\ {\cal E}_2/{\cal H})=8 k^2 [c^2
-k(g+1)^2] Z^6.
\eal
$$
Thus, we have $\deg\,{\cal H}=3$ and since $k\ne0$ according to
Lemma~\ref{Trudi:2} if $c^2 -k(g+1)^2\ne0$ we cannot have an
additional common factor of of ${\cal E}_1$ and ${\cal E}_2$. The
last condition is equivalent to $H_6\ne0$.

We observe, that the systems~\eqref{S3:10} have the invariant
lines $y=0$ and $gx^2+cx+k=0$. The positions of the last two lines
depend on the values of the $T$-comitants $H_{11}$ and $K$,
because for these systems we have
$$
H_{11}= 48(g-1)^4(c^2-4gk)x^4,\quad K=2g(g-1)x^2
$$

\textbf{a.1)}  Assume $K\ne0$. In this case $g\ne0$ and both
invariant lines $gx^2+cx+k=0$ are affine.

$\alpha)$ If $H_{11}>0$ then without loss of generality we may
introduce two new parameters $u$ and $v$ as follows:
$c^2-4gk=4g^2u^2>0$ and $c=2gv$. Then $k=g(v^2-u^2)$ and and since
$u\ne0$  we may assume $u=1$ via Remark \ref{rem:transf}
($\gamma=u,$ $s=1$). Then we obtain the systems
\be\label{S3:11}
  \dot x=g\big[(x+v)^2-1\big],\quad
  \dot y=(g-1)xy,
\ee
 which possess three real distinct
invariant lines ($y=0$ and $x+v=\pm1$) and two finite  singular
points: $\ (-v\pm1,\, 0)$. We note that $v\ne\pm1$, otherwise we
get degenerate systems \eqref{S3:11}. Thus for $K\ne0$ and
$H_{11}>0$ we obtain Config. 4.12.

$\beta)$ For $H_{11}<0$  as above we may set $c^2-4gk=-4g^2u^2<0$
and $c=2gv$. Then $k=g(v^2+u^2)$ and assuming $u=1$ due to Remark
\ref{rem:transf} ($\gamma=u,$ $s=1$) we obtain the systems
$$
 \dot x=g\big[(x+v)^2+1\big],\quad
  \dot y=(g-1)xy
$$
with invariant lines $y=0$ and $x+v=\pm i$. These systems have no
real singularities and we get Config. 4.15.

$\gamma)$ Assume now that $H_{11}=0$, i.e. $c^2=4gk$. Setting
$c=2gv$ we obtain $k=gv^2$ and this leads to the systems
$$
\dot x=g(x+v)^2,\quad
  \dot y=(g-1)xy.
$$
We observe that the line $x=-v$ is of multiplicity two (it cannot
be of multiplicity three in view of Lemma \ref{lm3}, as ${\cal
H}=Y(X+vZ)^2$). Since $v\ne0$ (otherwise systems become
degenerate) we may assume $v=1$ via Remark \ref{rem:transf}
($\gamma=v,$ $s=1$) and in this case we obtain Config. 4.24.

\textbf{a.2)}  Let us consider the case  $K=0$. Then $g=0$ and the
systems \eqref{S3:10} become
\be\label{S3:12}
  \dot x=k + cx,\quad
  \dot y=-xy
\ee
for which ${\cal H}=YZ(cX+kZ)$. Since in this case
$H_{11}=48c^2x^4$ we shall consider two cases: $H_{11}\ne0$ and
$H_{11}=0$.

$\alpha)$ If $H_{11}\ne0$ then $c\ne0$ and we may assume $c=1$ via
 Remark \ref{rem:transf} ($\gamma=c,$ $s=1$). So the systems
 \eqref{S3:12} with $c=1$ possess  invariant affine lines $y=0$
 and $x=-k$. Moreover, the line $Z=0$ is double as it is confirmed
 by the perturbations ({\it
IV.19${}_\varepsilon$})  from Table 3. Taking into account the
existence of a unique finite  singular point $(-k,0)$ we obtain
Config. 4.19.

$\beta)$ For $H_{11}=0$ we have $c=0$ and then we may assume
$k\in\{-1,1\}$ via Remark \ref{rem:transf} ($\gamma=|k|,$
$s=1/2$). In these cases the systems
 \eqref{S3:12} with $c=0$ possess  only one invariant affine line
 ($y=0$), and   the line $Z=0$ is triple, as it is confirmed
 by the perturbations ({\it
IV.36${}_\varepsilon$})   from Table 3. Therefore we get Config.
4.36.

\textbf{b)} Assume now $B_3\ne0$ and $H_6=0$. Then $l\ne0$ and the
condition $c^2 -k(g+1)^2=0$ is fulfilled. Since $g+1\ne0$ (due to
$N\ne0$) we may use a new parameter $u$: $c=u(g+1)$ and then
$k=u^2$. This leads to the systems
\be\label{S3:13}
  \dot x=(u+x)(u+ gx),\quad
  \dot y=l +(g-1)xy,
\ee
for which calculation yields:
$$
\bal
 &{\cal E}_1=2\big[(g-1) X Y + u(1-g) Y Z  + l(g+1) Z^2\big]{\cal
 H},\\
 &{\cal E}_2= (g X + cu Z) (-X Y + g X Y + l Z^2){\cal H},\\
&{\cal H}= \gcd\left({\cal E}_1,{\cal E}_2\right)=
 (X + u Z)^2 (g X + u Z),\\
& \Res_Z({\cal E}_1/{\cal H},\ {\cal E}_2/{\cal H})= 8
\,l\,(g-1)^2(g+1)\big[lg^2 X + u^2(g-1) Y\big]^2 X^2Y^2.
\eal
$$
So, we have $\deg\,{\cal H}=3$ and according to
Lemma~\ref{Trudi:2} to have  an additional common factor of ${\cal
E}_1$ and ${\cal E}_2$ for some  specific value of $(g,l,u)$ we
must have $\Res_Z({\cal E}_1/{\cal H},\ {\cal E}_2/{\cal H})=0$ in
$\R[X,Y]$. Therefore, since  $l(g^2-1)\ne0$, we obtain $g=u=0$,
however this condition leads to  degenerate systems
\eqref{S3:13}. Hence for $l(g^2-1)\ne0$ (i.e. $B_3N\ne0$) the  systems
\eqref{S3:13} belong to the class $\QSL_{\bf4}$.

 We observe that the  systems~\eqref{S3:13} possess invariant lines
$x+u=0$ and $gx+u=0$ and by Lemma \ref{lm3} and  the expression
for ${\cal H}$   the first line could be a double one. Clearly,
the values of the parameters $g$ and $u$ govern the position of
the invariant lines in the configuration. On the other hand for
systems~\eqref{S3:13} we have:
$$
K=2g(g-1)x^2,\quad H_{11}=48u^2(g-1)^6x^4.
$$

\textbf{b.1)} Assume   $KH_{11}\ne0$. Then $gu\ne0$ and we may set
$u=1$ via Remark \ref{rem:transf} ($\gamma=u,$ $s=1$). Moreover
due to the additional change $y\to ly$ we  may also assume $l=1$.
In this case the systems~\eqref{S3:13} possess two parallel
invariant affine lines $x=-1$ and $gx=-1$. The first line is a
double one as it is confirmed by the perturbations ({\it
IV.30${}_\varepsilon$}) from Table 3. Taking into account the
singular points $\big(-1,\,1/(g-1)\big)$ and $\big(g/(g-1),\,
-1/g\big)$ which are finite and distinct since $g(g-1)\ne0$, we
get Config. 4.30.

\textbf{b.2)} For $K\ne0$ and $H_{11}=0$ we have $g\ne0$, $u=0$
and the systems~\eqref{S3:13} possess the  invariant line $x=0$
which is of the multiplicity three. This is confirmed by the
perturbations   ({\it IV.43${}_\varepsilon$})   from Table 3.
Since $l\ne0$ (then we may assume $l=1$ via the change $y\to ly$)
these systems do not have  finite singular points. So we obtain
Config. 4.43.

\textbf{b.3)} Finally,  assume  that $K=0$, i.e.  $g=0$. Then
$u\ne0$ (otherwise we get degenerate systems) and as above we may
assume $u=1$ and $l=1$. So we obtain the system
\be\label{S3:14}
  \dot x=1+ x,\quad
  \dot y=1 -xy,
\ee
for which ${\cal H}=Z(X+Z)^2$, i.e. the line $x+1=0$ as well as
the line $Z=0$ are of the multiplicity two. This is confirmed by
the perturbations ({\it IV.40${}_\varepsilon$})   from Table 3.
Considering the existence of a simple finite singular point
$(-1,-1)$ we obtain Config. 4.40.

\fbox{\textbf{II.}} Assume $N=0$. According to \eqref{val:
kappa-N-3} the condition $\theta=\mu=0$ yields $h=0$ for systems
$(\SSS_{I\!I\!I})$ and then for these systems we have:
$$
  N=(g^2-1)x^2,\qquad K=2g(g-1)x^2.
$$
So, the condition $N=0$ implies either $g=1$ or $g=-1$ and we
shall consider two subcases: $K\ne0$ and $K=0$.

\textbf{1)} For $K\ne0$ we obtain $g=-1$ and  we may assume
$e=f=0$ via the  translation  $x\to x+f/2$ and $y\to y+e/2$. Then
calculations yield:
$$
B_2=-648\, d (8 c l x^4 + 16 d l x^3 y + d^3 y^4).
$$
Hence, the condition $B_2=0$ implies $d=0$  and we obtain the
systems:
\be\label{S3:15}
  \dot x=k  +cx - x^2,\quad
  \dot y=l-2xy.
\ee
for which calculations yield:
$$
\bal
 &{\cal E}_1=2(-2 X^2 Y + 2 k Y Z^2 + c l Z^3){\cal H},\quad
  {\cal E}_2= (X^2 - c X Z - k Z^2) (2 X Y - l Z^2){\cal H},\\
\eal
$$
where ${\cal H}= -X^2 + c X Z + k Z^2.$\ Thus, we have
$\deg\,{\cal H}=2$ and we observe that an additional common factor
of the polynomials ${\cal E}_1$ and ${\cal E}_2 $ must depend on
$X$ or/and on $Y$. According to Lemma~\ref{Trudi:2} this occurs
 if and only if at least one of the following two
identities holds:
\beq
&& \Res_X({\cal E}_1/{\cal H},\ {\cal E}_2/{\cal H})=
 -16 c^2 Y Z^6 (4 k Y^2 + 2c\,l Y Z  - l^2 Z^2)^2 \equiv0, \\
&&\Res_Y({\cal E}_1/{\cal H},\ {\cal E}_2/{\cal H})=-4\, l Z^2
(-X^2 + c X Z + k Z^2)^2 \equiv0.
\eeq
Since $k=l=0$ yield degenerate systems we obtain the condition
$cl=0$.

We claim that for the systems \eqref{S3:15} to be in the class
$\QSL_{\bf4}$ it is necessary to have $c\ne0$ and $l=0$. Indeed,
supposing that  $c=0$ we get the systems
$$
  \dot x=k - x^2,\quad
  \dot y=l-2xy,
$$
for which  calculations yield: ${\cal H}= \gcd\left({\cal
E}_1,{\cal E}_2\right)=2(X^2-kZ)^2$. Hence, $\deg\,{\cal H}=4$,
i.e. these systems belong to the class $\QSL_{\bf5}$ and not to
$\QSL_{\bf4}$.

Since   for systems \eqref{S3:15} we have  $B_3=-12lx^4$ and
$H_6=-2^{11}c^2 x^6$,  in what follows we shall consider $B_3=0$
and $H_6\ne0$ (i.e. $l=0$ and $c\ne0$). In this case we get the
systems
\be\label{S3:16}
  \dot x=k  +cx - x^2,\quad
  \dot y=-2xy,
\ee
which in fact  are  particular case of systems
\eqref{S3:10} (see page \pageref{S3:10}). More exactly, when
$g=-1$ from \eqref{S3:10} we obtain \eqref{S3:16}. On the other
hand  when $g=-1$ for systems \eqref{S3:10}  we have
$$
H_{11}= 768(c^2+4k)x^4,\quad K=4x^2,\quad H_6=-2^{11}c^2 x^6,
$$
i.e. $H_6K\ne0$. Then as it was proved on the page \pageref{S3:10}
if $H_{11}>0$ (respectively  $H_{11}<0$; $H_{11}=0$) we obtain
Config. 4.12 (respectively  Config. 4.15; Config. 4.24).

\textbf{2)} For $K=0$ and $N=0$ we obtain $g=1$ and via the
translation $x\to x-c/2$ and $y\to y$ the
systems~$(\SSS_{I\!I\!I})$ can be brought to the systems
$$
  \dot x=k +  dy + x^2,\quad
  \dot y=l+ex+fy,
$$
for which
$$
B_2=-648 d^4y^4,\quad B_3=6d(fx-dy)xy^2.
$$
Hence, the condition $B_2=0$ yields $d=0$ which implies $B_3=0$
and we obtain the systems:
\be\label{S3:17}
  \dot x=k  + x^2,\quad
  \dot y=l+ex+fy.
\ee
Calculations yield:
\beq
 &&{\cal E}_1=2\big[e X^2 + 2 f X Y + (2\,l- e f) X Z  - f^2 Y Z -
    (e k + f l) Z^2\big]{\cal H},\\
 &&{\cal E}_2= (e X + f Y + l Z) (X^2 + k Z^2){\cal H},\quad
{\cal H}= \gcd\left({\cal E}_1,{\cal E}_2\right)=Z(X^2+kZ^2).
\eeq
So, we have $\deg\,{\cal H}=3$ and
\beq
&& \Res_X({\cal E}_1/{\cal H},\ {\cal E}_2/{\cal H})= -8 e (f^2
+ 4 k) Z^2 \Big[f^2 Y^2 + 2 f l Y Z + (e^2 k + l^2)Z^2\Big]^2,\\
&& \Res_Y({\cal E}_1/{\cal H},\ {\cal E}_2/{\cal H})=-2ef(X^2+k^2Z^2),\\
 && \Res_Z({\cal E}_1/{\cal H},\ {\cal E}_2/{\cal H})=-8\, e(f^2 +4 k)
  \Big[(e^2 k + l^2) X^2 + 2 e f k X Y + f^2 k Y^2\Big]^2X^2.
\eeq
We note, that the conditions $f=e^2k+l^2=0$ yield degenerate
systems~(\ref{S3:17}), and  for $f=0$  both polynomials ${\cal
E}_1/{\cal H}$ and ${\cal E}_2/{\cal H}$ do not depend on $Y$.
Hence, according to Lemma~\ref{Trudi:2},  to have an additional
common factor of the polynomials ${\cal E}_1$ and ${\cal E}_2 $
the condition $e(f^2+4k)=0$ has to be satisfied. Therefore the
systems~(\ref{S3:17}) will belong to $\QSL_{\bf4}$ only if
$e(f^2+4k)\ne0$. On the other hand, for systems~(\ref{S3:17})  we
have:
$$
N_1=8ex^4,\quad  N_2=4(f^2+4k)x,\quad  N_5=-64kx^2
$$
and hence, the condition $e(f^2+4k)\ne0$ is equivalent to
$N_1N_2\ne0$.

The systems~(\ref{S3:17}) possess invariant affine lines $x^2+k=0$
which  could  be real, imaginary or could coincide  depending on
the value of the parameter $k$ (i.e. of the value  of $N_5$).
Moreover, since  $Z\mid{\cal H}$   the line $Z=0$ could be of the
multiplicity two and this is confirmed by the following
perturbations:
$$  \dot x=k  + x^2,\quad
  \dot y=(l+ex+fy)(1+\varepsilon\,y).
$$

\textbf{a)} Assume $N_5>0$. Then $k<0$ and since $e\ne0$  we may
consider $k=-1$ and $e=1$ via the transformation
$x\to(-k)^{1/2}x$, $y\to ey $ and $t\to(-k)^{-1/2}t$. So we get
the systems
\be\label{S3:18}
  \dot x= x^2-1,\quad
  \dot y=l+x+fy,
\ee
with $f\ne\pm2$. These systems    possess singular points
$\big(-1,\,(1-l)/f\big)$ and $\big(1,\,-(1+l)/f  \big)$ which tend
to infinity when $f\to 0$. We note that for the
systems~(\ref{S3:18}) we have $D=-f^2x^2y$.

Therefore, for $D\ne0$ (i.e. $f\ne0$) we may assume  $l=0$ via the
transformation $x\to x$ and $y\to y-l/f$ and taking into account
that $Z=0$ is double  we get Config. 4.28.

For $D=0$ we obtain $f=0$ and then $l\ne\pm1$, otherwise we get a
degenerate system \eqref{S3:18}. In this case we obtain Config.
4.29.

\textbf{b)} If $N_5<0$ then $k>0$ and since $e\ne0$  we may
consider $k=1$ and $e=1$ via the transformation $x\to k^{1/2}x$,
$y\to ey $ and $t\to k^{-1/2}t$. So we get the systems
\be\label{S3:19}
  \dot x= x^2+1,\quad
  \dot y=l+x+fy,
\ee
with $l,f\in \R$. These systems   have two imaginary invariant
lines ($x=\pm i$) and imaginary  singular points
$\big(-i,\,(i-l)/f \big)$ and $\big(i,\,-(i+l)/f \big)$.

Therefore, for $D=-f^2x^2y\ne0$ (i.e. $f\ne0$) we may assume $l=0$
via a transformation (as above) and taking into account that $Z=0$
is double  we obtain Config. 4.32.

For $D=0$ we obtain $f=0$ and then  we get Config. 4.33.

\textbf{c)} Assume now that  $N_5=0$, i.e.  $k=0$. Since $e\ne0$
we may assume $e=1$ via the change $y\to ey$. Then we obtain the
systems
\be\label{S3:20}
  \dot x= x^2,\quad
  \dot y=l+x+fy,
\ee
for which the condition $N_2=4f^2x\ne0$ yields $f\ne0$. Then via
the transformation $x_1=x/f$, $y_1=y+l/f$ and $t_1=ft$  we may
assume $l=0$ and $f=1$. Systems
\eqref{S3:20} possess the invariant line $x=0$ on which the double
point $(0,0)$ is placed. Taking into account that both the lines
$x=0$ and $Z=0$ are double (this is confirmed
 by the perturbations ({\it
IV.39${}_\varepsilon$})   from Table 3) we obtain Config. 4.39.

\subsection{Systems with\ $D_S(C,Z)=3\cdot \omega$}

In this subsection we shall consider the canonical system
$(\SSS_{IV})$ for which we have: $\theta =  8h^3$.

\subsubsection{The case $\theta\ne0$,\ $B_3=0$}

Then $h\ne0$ and  we can assume  $c=d=g=0$ via the affine
transformation:
$$
 x_1=x-\frac{d}{h},\quad y_1= \frac{g}{h}x+y+\frac{ch-dg}{h^2}.
$$
 Thus we obtain the canonical systems after returning to the same notations
 for the variables:
\be\label{S4:1}
  \dot x=k +  h x y,\quad
  \dot y=l + ex+fy -x^2+ hy^2.
\ee
We  calculate
\beq
&&B_3=3\,h x^2 \left[(ef-4k) x^2 +    2 h l x y - 3 h k
y^2\right],\quad
 H_7=-4eh^2,
\eeq
and hence,  by $\theta\ne0$ (i.e. $h\ne0$) the condition $B_3=0$
yields $k=l=fe=0$. We shall consider two subcases: $H_7\ne0$ and
$H_7=0$.

 \textbf{Subcase $H_7\ne0$.} Then $e\ne0$ and the
condition $B_3=0$ yields $k=l=f=0$. This leads to the systems
\be\label{S4:2}
  \dot x=  h x y,\quad
  \dot y= ex -x^2+ hy^2
\ee
for which we may assume $e=1$ and $h\in\{-1,1\}$ due to the
substitution: $x\to ex$, $y\to e|h|^{-1/2}y$ and $t\to
e^{-1}|h|^{-1/2}t$. Then for systems \eqref{S4:2} calculation
yields:
\beq
 &&{\cal E}_1=2(X^2  + h Y^2 - 2  X Z + Z^2){\cal H},\quad
 {\cal H}= \gcd\left({\cal E}_1,{\cal E}_2\right)=hX^3,\\
&& {\cal E}_2=  Y ^3{\cal H},\quad   \Res_Y({\cal E}_1/{\cal H},\
{\cal E}_2/{\cal H})= 8   (X -  Z)^6.
\eeq
So, we obtain $\deg\,{\cal H}=3$ and since\ $\Res_Y({\cal
E}_1/{\cal H},\ {\cal E}_2/{\cal H})\not\equiv0$,  according to
Lemma~\ref{Trudi:2} we cannot have additional nontrivial factors
of ${\cal E}_1$ and ${\cal E}_2$.

The invariant line $x=0$ is of the multiplicity 3 and this is
confirmed by the perturbations ({\it IV.45${}_\varepsilon$}) from
Table 3. Taking into account that systems \eqref{S4:2} possess the
singular points   $(0,0)$ (triple) and $(1,0)$ (simple), and the
last point is not located on the invariant line $x=0$, we obtain
Config. 4.45.

\textbf{Subcase $H_7=0$.} In this case the condition $B_3=0$
yields $k=l=e=0$ and we obtain  the systems
\be\label{S4:3}
  \dot x=  h x y,\quad
  \dot y= fy -x^2+ hy^2.
\ee
Calculations yield:
\beq
 &&{\cal E}_1=2(X^2 + h Y^2 - f Y Z){\cal H},\quad {\cal H}= \gcd\left({\cal E}_1,{\cal
E}_2\right)=\,hX^3,\\
 && {\cal E}_2=   h Y ^3{\cal H},\quad
\Res_Y({\cal E}_1/{\cal H},\ {\cal E}_2/{\cal H})= 8 h^4 X^6.
\eeq
So,  $\deg\,{\cal H}=3$ and from $h\ne0$ we have $\Res_Y({\cal
E}_1/{\cal H},\ {\cal E}_2/{\cal H})\not\equiv0$, i.e. we cannot
have an additional nontrivial factor  of ${\cal E}_1$ and ${\cal
E}_2$. The invariant line $x=0$ is of the multiplicity three and
this is confirmed by the perturbations ({\it
IV.41${}_\varepsilon$}) (for $f\ne0$) and ({\it
IV.42${}_\varepsilon$}) (for $f=0$)  from Table 3. We observe that
systems \eqref{S4:3} possess the singular points $(0,0)$ (triple)
and $(0,-f/h)$ (simple), and for these systems $D=-f^2x^3$.

Therefore, if $D\ne0$ then  $f\ne0$ and  we may assume $f=1$,
$h\in\{-1,1\}$  via the substitution $x\to f|h|^{-1/2}x$, $y\to
f|h|^{-1}y$ and $t\to f^{-1} t$. Since  in this case the singular
points above are distinct we obtain  Config. 4.41.

 Assume $D=0$, i.e. $f=0$. Then systems \eqref{S4:3} possess one point $(0,0)$
 of  multiplicity four and we get Config. 4.42. Note that in this case
we may assume $h\in\{-1,1\}$  via the substitution $x\to
|h|^{-1/2}x$, $y\to |h|^{-1}y$.

\subsubsection{The case $\theta=0=B_2$}

The condition $\theta=0$ yields $h=0$ and then
$$
B_2=-648d^2\big[\big((d+cg-fg)^2+g^2(2f^2-cf-gk)\big)x^2+4dg(d+cg-3fg)xy-6d^2g^2y^2\big]x^2.
$$
The condition $B_2=0$ yields $d=0$ and calculations yield:
\be\label{N-B_3new}
N=g^2x^2,\qquad B_3=3g(cf-f^2-gk)x^4
\ee
and we shall consider two subcases: $N\ne0$ and $N=0$.

\fbox{\textbf{Subcase $N\ne0$.}} Then $g\ne0$ and   we may assume
$g=1$ and $e=f=0$ via the transformation
$$
x_1= x+f/g,\quad y_1= gy+(2f+eg)/g,\quad t_1=g\,t.
$$
 So  keeping the previous notations we obtain the  systems
\be\label{S4:4}
  \dot x=  k+ c x + x^2,\quad
  \dot y= l -x^2+  xy
\ee
for which we calculate:
 \beq
 &&{\cal E}_1=2\big[X^3 + c X^2 Z + (2 k + l) X Z^2 - k Y Z^2 + c l Z^3\big]{\cal H},\\
&& {\cal E}_2=  (X^2 + c X Z + k Z^2)^2{\cal H},\quad
{\cal H}=  (X^2 + c X Z + k Z^2),\\
&& \Res_X({\cal E}_1/{\cal H},\ {\cal E}_2/{\cal H})= 16 k^2 Z^8
\left[k Y^2 + c(k  - l) Y Z + (k^2 - c^2l + 2 k l  + l^2)
Z^2\right]^2.
\eeq
Thus, we have $\deg\,{\cal H}=2$ and  according to
Lemma~\ref{Trudi:2} to have an additional   factor  of ${\cal
E}_1$ and ${\cal E}_2$  we must have $\Res_X({\cal E}_1/{\cal H},\
{\cal E}_2/{\cal H})=0$ in $\R[Y,Z]$. Therefore, we obtain the
condition $k=0$ which is equivalent to $B_3= 0$ since
$B_3=-3ky^4$. So, we get the systems
\be\label{S4:5}
  \dot x=   c x + x^2,\quad
  \dot y= l -x^2+  xy
\ee
for which we calculate again:
 \beq
&&{\cal E}_1=2(X^2 + l Z^2){\cal H},\quad
{\cal H}= \gcd\left({\cal E}_1,{\cal E}_2\right)= X (X + c Z)^2,\\
&& {\cal E}_2=  X^2 (X + c Z){\cal H},\quad
  \Res_X({\cal E}_1/{\cal H},\ {\cal E}_2/{\cal H})=8 l^2 (c^2
+ l) Z^6.
\eeq
The condition $l=0$ yields degenerate systems~(\ref{S4:5}). Hence,
in view of  Lemma~\ref{Trudi:2} in order to have  $\Res_X({\cal
E}_1/{\cal H},\ {\cal E}_2/{\cal H})\not\equiv0$ the condition
$c^2+l\ne0$ has to be satisfied. On the other hand for
systems~(\ref{S4:5}) we have $N_6=8(c^2+l)x^3$. So, the condition
$c^2+l\ne0$ is equivalent to $N_6\ne0$. We observe that
systems~(\ref{S4:5}) possess invariant lines: $x=0$ and $x+c$=0.
Moreover, the last line is of the multiplicity two, and for $c=0$
the line $x=0$ is of the multiplicity three. This is confirmed by
the perturbations ({\it IV.31${}_\varepsilon$}) (for $c\ne0$) and
({\it IV.44${}_\varepsilon$}) (for $c=0$) from Table 3.

On the other hand for systems~(\ref{S4:5}) we have
$H_{11}=48c^2x^4$ and hence the condition $c=0$ is equivalent to
$H_{11}=0$.  So, for $H_{11}\ne0$ we obtain $c\ne0$ and we may
assume $c=1$ via Remark \ref{rem:transf} ( $\gamma=c,\ s=1$). This
leads to Config. 4.31.

Assuming $c=0$ since $l\ne0$ we may consider $l\in\{-1,1\}$ via
Remark \ref{rem:transf} ( $\gamma=|l|,\ s=1/2$) and we obtain
Config. 4.44.

\fbox{\textbf{Subcase $N=0$.}} Then $g=0$ and from
\eqref{N-B_3new} we obtain $B_3=0$. We may assume $e=0$ via the
translation $x\to x+e/2$ and $y\to y$ and therefore we obtain the
systems
\be\label{S4:6}
  \dot x=  k+ c x,\quad
  \dot y= l + f y  -x^2
\ee
for which calculations yield:
 \beq
 &&{\cal E}_1=2\big[(c+f) X^2 +2k X Z  + f(c - f) Y Z +
    l(c-f) Z^2\big]{\cal H},\\
&& {\cal E}_2=Z (cX + k Z)^2{\cal H},\quad {\cal H}=
\gcd\left({\cal E}_1,{\cal E}_2\right)=Z^2 (cX + k Z).
\eeq
Thus,  $\deg\,{\cal H}=3$ and, since the polynomial ${\cal E}_2$
does not depend on $Y$, according to Lemma~\ref{Trudi:2} we could
have an additional common factor of ${\cal E}_1$ and $ {\cal E}_2$
   only if at least one of the following
conditions holds:
\beq
&& \Res_X({\cal E}_1/{\cal H},\ {\cal E}_2/{\cal H})= 4 (c - f)^2
Z^4\left[c^2 f Y + (c^2l  - k^2 ) Z\right]^2\equiv0,\\
&& \Res_Z({\cal E}_1/{\cal H},\ {\cal E}_2/{\cal H})= 8 (c - f)^2
(c + f) X^4 \left[(c^2l   - k^2 ) X - c f k Y\right]^2\equiv0,
\eeq
which amount to at least  one of the  following   conditions:
\beq
(i)\ (c-f)(c+f)=0;\ \ (ii)\ cf=c^2l- k^2=0;\ \ (iii)\
cfk=c^2l-k^2=0.
\eeq
We observe, that the conditions $(ii)$ yield either $c=k=0$ or
$f=c^2l-k^2=0$ and both these cases lead to  degenerate systems
(\ref{S4:6}). If $k\ne0$ the conditions $(iii)$ are equivalent to
$(ii)$, and for $k=0$ we obtain $l=0$. However the conditions
$k=l=0$ do  not imply the existence of an additional  common
factor for the polynomials ${\cal E}_i/{\cal H}$ $(i=1,2)$ unless
$c+f=0$ which falls in the case $(i)$.

Thus, we do not have an additional  common factor for the
polynomials ${\cal E}_1$ and  ${\cal E}_2$ if $c^2-f^2\ne0$.

We observe that the systems~(\ref{S4:6}) possess invariant affine
lines $cx+k=0$ only for $c\ne0$. On the other hand, for these
systems we have
$$
N_3=3(c-f)x^3,\quad D_1=c+f,\quad N_6=8c(c-f)x^3.
$$
Hence, the condition $c^2-f^2\ne0$ is equivalent to $D_1N_3\ne0$
and for $N_3\ne0$ the condition $c=0$ is equivalent to $N_6=0$.

\textbf{1)} Assume firstly $N_6\ne0$, i.e. $c\ne0$. Then we may
consider $c=1$ via Remark \ref{rem:transf} ( $\gamma=c,\ s=1$). In
this case the systems~(\ref{S4:6}) possess the invariant straight
lines $x+k=0$ and the singular point $\big(-k, (k^2-l)/f \big)$
which is finite for $f\ne0$ and it tends to infinity  when $f$
tends to zero. Since for systems~(\ref{S4:6}) $D=-f^2x^3$ we
conclude that the condition $f\ne0$ is captured by the
$T$-comitant $D$.

  Assuming $D\ne0$ we obtain $f\ne0$ and then we may assume
$l=0$ due to a translation. Thus we obtain the systems
\be\label{S4:7}
  \dot x=  k+  x,\quad
  \dot y=  f y  -x^2
\ee
with $f(f^2-1)\ne0$. Moreover we may assume $k\in\{0,1\}$ via the
rescaling $x\to kx$ and $y\to k^2y$ (for $k\ne0$).  Taking into
account that the line $Z=0$ is triple (and this is confirmed by
the perturbations  ({\it IV.37${}_\varepsilon$}) from Table 3) we
get Config. 4.37.

If $D=0$ we have $f=0$ and this leads to the systems
\be\label{S4:7a}
  \dot x=  k+  x,\quad
  \dot y=  l - x^2
\ee
with $l-k^2\ne0$ and as above we may assume $k\in\{0,1\}$. Taking
into account that these systems do not have  finite singular
points and   the line $Z=0$ is triple ( this is confirmed by the
perturbations ({\it IV.38${}_\varepsilon$}) from Table 3) we get
Config. 4.38.

\textbf{2)} Assume now $N_6=0$, i.e. $c=0$. Then the condition
$D_1N_3\ne0$ yields $f\ne0$ and since $k\ne0$  for systems
\eqref{S4:6} we may assume $k=f=1$ and $l=0$ via the
transformation
$$
x_1=fk^{-1}x,\quad y_1=f^3k^{-2}y+lf^2k^{-2},\quad t_1=ft.
$$
Hence we obtain the system
\be\label{S4:8}
  \dot x=  1,\quad
  \dot y=  y - x^2
\ee
for which the line $Z=0$ is of the multiplicity four as it is
confirmed by the perturbations ({\it IV.46${}_\varepsilon$}) from
Table 3. Thus in this case  we get Config. 4.46.

\medskip
All the cases in  Theorem \ref{th_mil_4} are thus examined. To
finish the proof of the Theorem \ref{th_mil_4} it remains to show
that the conditions  occurring in the middle column of Table 2 are
affinely invariant. This follows from the proof of Lemma
\ref{Table:Propreties}.
  \EProof

\blm\label{Table:Propreties}
   The  polynomials which are used
   in  Theorem \ref{th_mil_4}   have the properties
indicated in the Table 4. In the last column are indicated the
algebraic sets on which the  $GL$-comitants on the left are
$CT$-comitants.
\elm
\begin{table}[!htb]
\begin{center}
\begin{tabular}{|c|c|c|c|c|c|}
\multicolumn{6}{r}{\bf Table 4}\\[2mm]
\hline \raisebox{-0.7em}[0pt][0pt]{Case} &
\raisebox{-0.7em}[0pt][0pt]{$GL$-comitants}
 & \multicolumn{2}{c|}{Degree in }  & \raisebox{-0.7em}[0pt][0pt]{Weight} &  Algebraic subset \\
\cline{3-4}
  &  & $\ \ a\ \ $ & $\!x$ and $y\!$ &    &  $V(*)$   \\
\hline\hline \rule{0pt}{2ex}
 $1$ & $\eta(a)$,\ $\mu(a)$,\ $\theta(a)$ & $4$ &  $0$ & $ 2$  &  $V(0)$ \\[0.5mm]
\hline\rule{0pt}{2ex} $2$  & $C_2(a,x,y)$   & $1$ &  $3$ & $-1$  & $V(0)$\\[0.5mm]
\hline\rule{0pt}{2ex} $3$ & $H(a,x,y),\ K(a,x,y)$  & $2$ &  $2$ & $ 0$   &  $V(0)$ \\[0.5mm]
\hline\rule{0pt}{2ex} $4$ & $M(a,x,y),\ N(a,x,y)$   & $2$ &  $2$ & $ 0$   &  $V(0)$ \\[0.5mm]
\hline\rule{0pt}{2ex} $5$ & $D(a,x,y)$  & $3$ &  $3$ & $-1$   &  $V(0)$ \\[0.5mm]
\hline\rule{0pt}{2ex} $6$ & $B_1(a)$  & $12$ &  $0$ & $3$   &  $V(0)$ \\[0.5mm]
\hline\rule{0pt}{2ex} $7$ & $B_2(a,x,y)$  & $8$ &  $4$ & $0$   &  $V(0)$ \\[0.5mm]
\hline\rule{0pt}{2ex} $8$ & $B_3(a,x,y)$  & $4$ &  $4$ & $-1$   &  $V(0)$ \\[0.5mm]
\hline\rule{0pt}{2ex} $9$ & $H_1(a)$  & $6$ &  $0$ & $2$   &  $V(0)$ \\[0.5mm]
\hline\rule{0pt}{2ex} $10$ & $H_2(a,x,y)$  & $3$ &  $2$ & $0$   &  $V(0)$ \\[0.5mm]
\hline\rule{0pt}{2ex} $11$ & $H_3(a,x,y)$  & $4$ &  $2$ & $0$   &  $V(0)$ \\[0.5mm]
\hline\rule{0pt}{2ex} $12$ & $H_4(a)$,\ $H_{10}(a)$  & $6$ &  $0$ & $2$   &  $V(0)$ \\[0.5mm]
\hline\rule{0pt}{2ex} $13$ & $H_5(a)$,\ $H_8(a)$  & $8$ &  $0$ & $2$   &  $V(0)$ \\[0.5mm]
\hline\rule{0pt}{2ex} $14$ & $H_6(a,x,y)$   & $8$ &  $6$ & $0$   &  $V(0)$ \\[0.5mm]
\hline\rule{0pt}{2ex} $15$ & $H_7(a)$  & $3$ &  $0$ & $1$   &  $V(0)$ \\[0.5mm]
\hline\rule{0pt}{2ex} $16$ & $H_9(a))$  & $12$ &  $0$ & $2$   &  $V(0)$ \\[0.5mm]
\hline\rule{0pt}{2ex} $17$ & $H_{11}(a,x,y))$  & $6$ &  $4$ & $0$   &  $V(0)$ \\[0.5mm]
\hline \rule{0pt}{2ex} $18$ & $N_1(a,x,y)$  & $3$ &  $4$ & $-1$   & $V(\eta,N,K)$ \\[0.5mm]
\hline\rule{0pt}{2ex}  $19$ & $N_2(a,x,y)$  & $3$ &  $1$ & $0$   &  $V(\eta,N,K,B_3)$ \\[0.5mm]
\hline\rule{0pt}{2ex}  $20$ & $N_3(a,x,y)$  & $2$ &  $3$ & $-1$   &  $V(M,N)$ \\[0.5mm]
\hline\rule{0pt}{2ex}  $21$ & $N_4(a,x,y)$  & $2$ &  $2$ & $-1$   &  $V(M,N,N_3)$ \\[0.5mm]
\hline\rule{0pt}{2ex}  $22$ & $N_5(a,x,y)$  & $4$ &  $2$ & $0$   &  $V(\eta,N,K,B_3)$ \\[0.5mm]
\hline\rule{0pt}{2ex}  $23$ & $N_6(a,x,y)$  & $3$ &  $3$ & $-1$   &  $V(M,\theta,B_3)$ \\[0.5mm]
\hline\rule{0pt}{2ex} $24$ & $D_1(a)$  & $1$ &  $0$ & $0$   &  $V(M,N)$ \\[0.5mm]
\hline
\end{tabular}
\end{center}
\end{table}

\BProof  {\it I. Cases  1--17}.  Assume that $\ab\in \R^{12}$ corresponds to an arbitrarily
given  system \eqref{2l1} and assume  $\tilde\ab\in \R^{12}$
corresponds to a system in the orbit of the given system
\eqref{2l1} under the action of the translation group, i.e. if
$\tau:\  x = \tilde x+x_0,$ $y=\tilde y+y_0$ then
$$
\tilde\ab:\quad \left\{\ba{l}\dot {\tilde x}=
P(\ab,x_0,y_0)+P_x(\ab,x_0,y_0)\tilde x
           +P_y(\ab,x_0,y_0)\tilde y+p_2(\ab,\tilde x,\tilde y), \\
\dot {\tilde y}= Q(\ab,x_0,y_0)+Q_x(\ab,x_0,y_0)\tilde
x+Q_y(\ab,x_0,y_0)\tilde y+q_2(\ab,\tilde x,\tilde y).\ea\right.
$$
Then for every $\ab\in\R^{12}$ and $(\tilde x,\tilde y)\in \R^2$
calculations yield:
$$
\bal
  & U (\tilde\ab)=U (\ab)\quad
  \text{for each}\quad U\in \{\eta,\mu,\theta,B_1,H_1,H_4,H_5,H_7,H_8,H_9,H_{10}\}, \\
 & W (\tilde\ab,\tilde x,\tilde y)=W (\ab,\tilde x,\tilde y)\quad
  \text{for each}\quad W\in \{C_2,K,H,M,N,D,B_2,B_3,H_2,H_3,H_6,H_{11}\}. \\
\eal
$$
Hence according to the definition of   $T$-comitants (see
\cite{Dana_Vlp1}) we conclude that the GL- comitants indicated in
the lines 1--17 of Table 4 are $T$-comitants for systems
\eqref{2l1}.

{\it II. Cases 18--24}. {\bf 1)} We consider firstly the
$GL$-comitants $N_1(a,x,y)$, $N_2(a,x,y)$ and $N_5(a,x,y)$ and we
shall  prove that $N_2$ and $N_5$ (respectively $N_1$) are
$CT$-comitants modulo $\langle\eta,N,K,B_3\rangle$ (respectively
modulo $\langle\eta,N,K \rangle$).  We shall examine the two
subcases: $M\ne0$ and $M=0$.

{\bf a)} For $\eta=0$ and $M\ne0$ we are in the  class of the
systems $(\SSS_{I\!I\!I})$, for which the conditions $
N=(g^2-1)x^2 +2h(g-1)xy+h^2y^2=0$ and $K=2g(g-1)x^2+4gh
xy+h^2y^2=0$ yield $h=g-1=0$. Hence applying the additional
translation  $x\to x-c/2$, $y\to y$ we obtain the systems
\be\label{S3_NH_0}
\dot x=k +dy+ x^2,\qquad \dot y= l+ex+fy.
\ee
 On the other hand  for any system
corresponding to a point $\tilde\ab\in \R^{12}$ in the orbit under
the translation group action of a system \eqref{S3_NH_0}
calculations yield:
$$
\bal
& N_1(\tilde \ab,\tilde x,\tilde y)=8\tilde x^2(e\tilde
x^2-2d\tilde y^2),\quad N_2(\tilde \ab,\tilde x,\tilde
y)=4(f^2+4k)\tilde x + 4df  \tilde y+8d (x_0\tilde y +2y_0\tilde x),\\
& N_5(\tilde \ab,\tilde x,\tilde y)=-16(4k \tilde x^2 - d^2
 \tilde y^2) + 64d\tilde x (x_0\tilde y -y_0\tilde x),\quad B_3(\tilde \ab,\tilde x,\tilde y)
  = 6d\tilde x\tilde y^2(f\tilde x-d\tilde y).
\eal
$$
 Hence the polynomial $N_1$
does not depend on the vector defining the translations and for
$B_3=0$ the same occurs for the polynomials $N_2$ and $N_5$.
Therefore we conclude that for $M\ne0$ the polynomial $N_1$  is a
$CT$-comitant modulo $\langle\eta,N,K,\rangle$, whereas the
polynomials $N_2$ and $N_5$ are $CT$-comitants modulo
$\langle\eta,N,K,B_3\rangle$.

{\bf b)} Assume now that $M=0$. Then we are in the class of the
systems $(\SSS_{I\!V})$, for which the condition
 $N=(g^2-2h)x^2+2ghxy+h^2y^2=0$
yields $h=g=0$ and then $K=0$. Then applying the additional
translation  $x\to x+e/2$, $y\to y$, we obtain the systems
\be\label{S4_N0}
\dot x=k+cx+dy,\qquad \dot y= l +fy-x^2.
\ee
 For any system corresponding to a point
$\tilde\ab\in \R^{12}$ in the orbit under the translation group
action of a system \eqref{S4_N0}  calculations yield:
$$
 N_1(\tilde \ab,\tilde x,\tilde y)=-24 d\tilde x^4,
 \quad N_2(\tilde \ab,\tilde x,\tilde
y)=12d(c+f)\tilde x,\quad N_5(\tilde \ab,\tilde x,\tilde
y)=0=B_3(\tilde \ab,\tilde x,\tilde y).
$$
Since the condition $M=0$ implies $\eta=0$, considering the case
{\bf 1\,a)} above, we conclude that independently of either
$M\ne0$ or $M=0$ the $GL$-comitant $N_1$  is a $CT$-comitant
modulo $\langle\eta,N,K\rangle$ and $N_2$ and $N_5$ are
$CT$-comitants modulo $\langle\eta,N,K,B_3\rangle$.

{\bf 2)} Let us now consider  the $GL$-comitants $N_3(a,x,y)$,
$N_4(a,x,y)$, $N_6(a,x,y)$ and $D_1(a)$. According to Table 4  we
only need to examine   the class of the systems $(\SSS_{I\!V})$
and we shall consider the two subcases: $N=0$ and $N\ne0$,
$\theta=0$.

{\bf a)} If for  a system $(\SSS_{I\!V})$  the condition $N=0$ is
fulfilled then   it was shown above  that this system can be
brought via a translation to the form \eqref{S4_N0}. For any
system corresponding to a point $\tilde\ab\in \R^{12}$ in the
orbit under the translation group action of a system
\eqref{S4_N0} calculations yield:
$$
\bal
& N_3(\tilde \ab,\tilde x,\tilde y)=3(c-f)\tilde x^3+2d\tilde
x^2\tilde y,\quad
B_3(\tilde \ab,\tilde x,\tilde y)=6d\tilde x^3(f\tilde x-d\tilde y),\\
& N_4(\tilde \ab,\tilde x,\tilde y)=12 k\tilde x^2 +
3(f^2-c^2)\tilde x\tilde y
               -3d(c+f)\tilde y^2+6\tilde x^2[(c-f)x_0+2dy_0],\\
& N_6(\tilde \ab,\tilde x,\tilde y)=8c(c-f)\tilde x^3 + 16 df
\tilde x^2 \tilde y-8d^2\tilde x\tilde y^2-48d x_0\tilde x^3,\quad
D_1(\tilde \ab)=c+f.
\eal
$$
These relations show us that: {\it (i)} the $GL$-comitants $N_3$
and $D_1$ are  $CT$-comitant modulo $\langle M,N\rangle$; {\it
(ii)} the $GL$-comitant $N_4$ is a $CT$-comitant modulo $\langle
M,N,N_3\rangle$; {\it (iii)} the $GL$-comitant $N_6$ is a
$CT$-comitant modulo $\langle M,N,B_3\rangle$.

{\bf b)} Assume that for the systems $(\SSS_{I\!V})$  the
conditions $\theta=8h^3=0$  and $N=(g^2-2h)x^2+2ghxy+h^2y^2\ne0$
are fulfilled.  Then $h=0$, $g\ne0$ and we may assume $g=1$ and
$e=f=0$ via the transformation $\ x_1= x+f/g,\ \ y_1=
gy+(2f+eg)/g,\quad t_1=g\,t.$\ Then we obtain the systems\\[2mm]
\centerline {$  \dot x=  k+ c x +dy + x^2,\quad
  \dot y= l -x^2+  xy$}\\[2mm]
for which  calculation yields: $
B_3=-3x^2\big[kx^2+2d(c+d)xy+3d^2y^2. $ Therefore the condition
$B_3=0$ yields $k=d=0$ and we obtain the following family of
systems which is  characterized by the conditions $M=\theta=B_3=0$
and $N\ne0$:
\be\label{S_IV_1g}
  \dot x=   c x   + x^2,\quad
  \dot y= l -x^2+  xy.
\ee
  For any system corresponding to a point
$\tilde\ab\in \R^{12}$ in the orbit under the translation group
action of a system \eqref{S_IV_1g} we have $N_6(\tilde \ab,\tilde
x,\tilde y)=8(l+3f^2)\tilde x^3$. Since the condition $N=0$
implies $\theta=0$, considering the case
 {\bf 2\,a)} above, we conclude in this case that independently of either
$N\ne0$ or $N=0$, the $GL$-comitant $N_6$  is a $CT$-comitant
modulo $\langle M,\theta,B_3\rangle$.

The Table 4 show us that all the conditions indicated in the
middle column of Table 2 are affinely invariant. Indeed, the
$CT$-comitants $N_i$, $i=1,\ldots,...,6$ and $D_1$ are used in
Table 2 only for the varieties indicated in the last column of
Table 4. This complete the proof of the Theorem \ref{th_mil_4}.
 \EProof

\medskip

 \centerline{\bf Acknowledgements}

\medskip

 {\it The second author is very thankful for the kind
hospitality of the Centre de Recherches Math\'ematiques and
D\'epartement de Math\'ematiques et de Statistique of the
Universit\'e de Montr\'eal. Special thanks are due to Yvan St
Aubin, Director of the D\'epartement de Math\'ematiques et de
Statistique of the Universit\'e de Montr\'eal.}



\end{document}